\documentclass[graybox]{SNmult}

\makeatletter
\def\l@title#1#2{}
\def\l@author#1#2{}
\def\l@section#1#2{\@dottedtocline{1}{\tocchpnum}{\tocsecnum}{\bfseries #1}{\bfseries #2}}
\def\tableofcontents{\section*{\contentsname}\@starttoc{toc}}
\providecommand{\authcount}[1]{}
\makeatother

\usepackage{etoolbox}
\makeatletter
\patchcmd{\@ssect}{\addcontentsline{toc}{section}{#5}}{}{}{}
\makeatother

\usepackage{amsmath}
\usepackage{type1cm} 
\usepackage{makeidx} 
\usepackage{graphicx} 
\usepackage{multicol} 
\usepackage[bottom]{footmisc}

\usepackage{newtxtext} 
\usepackage[varvw]{newtxmath}

\DeclareMathAlphabet{\mathbbe}{U}{bbold}{m}{n}
\newcommand{\bbOmega}{\mathbbe{\Omega}}

\usepackage{xcolor}

\graphicspath{{../}}
\usepackage{type1cm}        
\usepackage{makeidx}         
\usepackage{graphicx}        
\usepackage{multicol}        
\usepackage[bottom]{footmisc}

\usepackage{newtxtext}       %
\usepackage[varvw]{newtxmath}       

\makeindex             
                       
\usepackage{graphics}
\usepackage{xspace}
\usepackage{bm,mathrsfs}

\usepackage[colorlinks=true, linkcolor=blue!60!black, citecolor=blue!60!black, urlcolor=blue!60!black]{hyperref}

\numberwithin{equation}{section}

\newcommand{\qp}[1]{\ensuremath{\!\left({#1}\right)}}

\def\R{{\mathbb{R}}}

\def\EE{{\mathbb{E}}}

\def\d{{\mathrm d}}

\def\Th{{\mathcal{T} }_h }

\def\S{{\mathcal{S} }}

 \renewcommand{\P}{\mathbb{P}}

 \renewcommand{\L}{\mathcal{L}}

\def\<{{\langle }}
\def\>{{\rangle }}

\def\wU{\widehat U}

\def\PP{{\mathbb{P}}}

\begin{document}
\title*{Variational Principles and Machine Learning: a Numerical Analysis and PDE perspective}
\titlerunning{Variational Principles and Machine Learning}
\author{Charalambos G.\ Makridakis\orcidID{0009-0005-7430-1812} 
}
\institute{Charalambos G.\ Makridakis \at University of Sussex and IACM-FORTH \email{C.Makridakis@sussex.ac.uk}}
%
\maketitle
\renewcommand{\theequation}{\thesection.\arabic{equation}}
\abstract*{Variational principles have played a fundamental role in the development of modern mathematics, physics, and scientific computing. From the principles of least action in mechanics to the calculus of variations and the variational formulation of partial differential equations, they have provided a powerful framework for understanding complex problems and designing reliable numerical methods. In recent years, Machine Learning has introduced a new computational paradigm, in which highly expressive nonlinear function spaces and large-scale minimisation problems play a central role.
In these notes, we explore the connection between these two frameworks. We develop a variational and numerical-analysis perspective on Machine Learning algorithms, focusing in particular on their formulation as nonlinear minimisation problems. We investigate how classical concepts such as stability, consistency, approximability, and compactness acquire new interpretations in this setting, and how the variational structure of the underlying continuous problem can govern the behaviour of neural-network approximations. This perspective provides a mathematical framework for understanding existing ML algorithms while also suggesting new approaches to their design. In particular, alternative and strongly coercive variational principles can lead to more stable, reliable, and potentially more efficient computational methods.
}
\abstract{Variational principles have played a fundamental role in the development of modern mathematics, physics, and scientific computing. From the principles of least action in mechanics to the calculus of variations and the variational formulation of partial differential equations, they have provided a powerful framework for understanding complex problems and designing reliable numerical methods. In recent years, Machine Learning has introduced a new computational paradigm, in which highly expressive nonlinear function spaces and large-scale minimisation problems play a central role.
In these notes, we explore the connection between these two frameworks. We develop a variational and numerical-analysis perspective on Machine Learning algorithms, focusing in particular on their formulation as nonlinear minimisation problems. We investigate how classical concepts such as stability, consistency, approximability, and compactness acquire new interpretations in this setting, and how the variational structure of the underlying continuous problem can govern the behaviour of neural-network approximations. This perspective provides a mathematical framework for understanding existing ML algorithms while also suggesting new approaches to their design. In particular, alternative and strongly coercive variational principles can lead to more stable, reliable, and potentially more efficient computational methods.
}

\setcounter{tocdepth}{3}        
\newpage
{\small\tableofcontents}
\markboth{C.~G.~Makridakis}{Variational Principles and Machine Learning}


\section{Introduction}
\label{sec:1}

 \subsubsection*{Explainable artificial intelligence and simple ML algorithms}

In seeking to understand the main mechanisms that drive artificial intelligence (AI) algorithms, one encounters a lack of understanding even for the simplest machine learning (ML) processes, i.e., algorithms designed to perform basic tasks such as learning a function or solving simple differential equations. Such problems belong to the elementary curriculum in mathematics, engineering, and computer science. Classical methods for these tasks are known to be accurate and reliable, and it is therefore somewhat surprising that their machine learning counterparts are often complicated and lack interpretability. The corresponding machine learning algorithms are in an unfavourable position compared to standard approaches in terms of both computational cost and reliability, and this has led several scientists to question their use in scientific computing \cite{Grossmann_PINNs_vs_FEM:2024, McGreivy_Hakim_weak_baselines:2024}, acknowledging however that ML methods may offer advantages in high-dimensional settings where classical approaches become infeasible \cite{Han_Jentzen_E_PNAS:2018, e2017deep, Grohs_Jentzen_Schwab_highD_pde_DNN:2022}.\footnote{This statement may seem to contradict the remarkable achievements of AI systems reported at the time of writing, including contributions to long-standing open problems in mathematics. The lack of understanding referred to here concerns \emph{why} the algorithms work, not \emph{what} they can do. Moreover, in mathematics a proof can be certified independently of the process that produced it, by a proof assistant such as Lean.  In scientific computing no such universal certification mechanism exists, reliability must come from the method itself, and a method is judged by its computational cost as much as by its accuracy.}

However, a lack of understanding has historically been an important driving
force for scientific discovery. The rapid expansion of Machine Learning for
scientific applications was, to some extent, unexpected: it was enabled
primarily by computational and algorithmic advances — efficient optimisation
over large-scale nonlinear spaces, automatic differentiation, stochastic
methods, and hardware acceleration \cite{Baydin_AD_survey:2018,
Bottou_Curtis_Nocedal_SIREV:2018} — while the mathematical understanding of
why the resulting algorithms work remains incomplete. As a result, a
significant gap has opened between what these methods can achieve in
practice and what can currently be understood or guaranteed analytically
\cite{Grohs_Petersen_Review:2023, Bartlett_Montanari_Rakhlin_actaN:2021}.
This situation is not without precedent. In the early development of
scientific computing, too, practice ran ahead of theory: the reliability of
finite difference schemes came to be understood only through the stability
analysis pioneered by J.~von Neumann and the equivalence theorem of Lax and
Richtmyer \cite{Lax_Ricth_1956}; see \cite{Hou_Lax_1991} for a historical
account.

A first indication of what numerical analysis can offer is provided by a
phenomenon that has attracted considerable attention: the tendency of ML
algorithms to produce outputs that appear plausible but are incorrect, often
referred to as \emph{hallucinations}. In Section~\ref{SubS:hallucinations} we
argue that such behaviour is closely related to a deep phenomenon in the
numerical analysis of nonlinear PDEs: the non-uniqueness of weak solutions and
the resulting sensitivity to small perturbations. The analogue in machine
learning is the nonlinear structure and the selection of a local minimiser of a nonconvex loss among many.
Just as a physically irrelevant weak solution is a genuine solution of the
equations, a local minimiser is a genuine minimiser of the loss on a
neighbourhood, may achieve a small training loss, and may appear to be a
reasonable approximation, while being far from the solution sought. In both
cases reliability rests on additional structure: selection criteria, such as
entropy conditions for conservation laws, single out the physically relevant
weak solution, { while, as we shall see in Section~\ref{S:Analysis}, stability
of the loss guarantees that local minimisers with consistently small loss are
close to the true solution.} Section~\ref{S:typical} also shows that other
ingredients which may appear specific to ML have classical counterparts:
nonlinear approximation spaces arise already in finite element methods with
free nodes, and gradient descent, the workhorse of training, is the explicit
Euler discretisation of a gradient flow in parameter space.

In these notes, we focus on a key aspect of machine learning algorithms, namely their formulation as \emph{nonlinear} minimisation problems. Our aim is not only to analyse these problems as optimisation procedures, but also to understand them from the broader perspective of Numerical Analysis and the Calculus of Variations. Our goals are threefold: (i) to gain a deeper mathematical understanding of the behaviour of ML algorithms using tools from the Calculus of Variations and Numerical Analysis; (ii) to relate their behaviour to classical notions in the analysis of numerical algorithms, in particular stability, consistency, and approximability; and (iii) to demonstrate how analytical principles can guide the design of new algorithms with improved stability, reliability, and computational properties. These classical notions acquire new and meaningful interpretations in the ML setting. The nonlinear nature of neural-network approximation spaces, the non-convexity of the resulting minimisation problems, and the use of empirical or stochastic approximations introduce new challenges, but also reveal connections with well-established ideas in numerical analysis. Throughout, we use the word \emph{explainable} in the sense familiar from
numerical analysis and scientific computing: to explain an algorithm is to
identify the mechanisms that determine its behaviour, rather than to interpret
the output of a particular trained model.

Nonlinear minimisation problems have been studied in many fields of
mathematics. The calculus of variations provides a powerful machinery for
their analysis, and a central part of it is due to Ennio De Giorgi. His notion
of $\Gamma$-convergence, introduced with Franzoni \cite{DeGiorgi_Franzoni:1975},
was designed precisely to describe the limiting behaviour of sequences of
nonlinear minimisation problems, and in particular of their minimisers, see
\cite{DalMaso_Gamma:1993, braides2002gamma}. Its key feature, and the reason it
is so well suited to our purposes, is that it operates on the functionals in
function space: it sees only the energies and their minimisers, and neither
the parametrisation of the approximation spaces nor the procedure used to
compute the minimisers. Convergence of minimisers follows from a \emph{liminf}
inequality, combined with a suitable compactness property, and a
\emph{limsup} inequality, which expresses the existence of recovery
sequences. In Section~\ref{S:Abstract-liminf-limsup} we interpret this
framework from the perspective of numerical analysis, linking the liminf
inequality to stability and the limsup inequality to consistency and
approximability. Applied to our setting, it becomes yet another instance of
the principle that stability and consistency imply convergence, and it is
precisely what is needed to analyse minimisation over neural-network spaces,
since it bypasses the nonconvexity of the problem in parameter space.
$\Gamma$-convergence has been used in the numerical analysis of nonconvex and
singularly perturbed variational problems, see e.g.
\cite{bartels2017bilayer, grekas2022approximations}; machine learning offers a
natural new field of application. Although $\Gamma$-convergence is, in its standard form, a qualitative tool,
the liminf--limsup framework yields more in our setting. The liminf
inequality alone provides a posteriori error control for approximate local
minimisers through the computed losses, while, as the number of samples
grows, the accumulation points of global minimisers satisfy estimates in
terms of the best approximation error of the neural-network space, which
approximation results translate into rates.

We now come to the central theme of these notes, reflected in their title. A
given problem can typically be described by several variational principles,
equivalent at the continuous level but with very different mathematical
properties. Mechanics offers an early example: in contrast to Hamilton's
principle, which characterises the motion as a stationary point of the
action, Gauss's principle of least constraint is a genuine minimum principle,
formulated as a least-squares problem for Newton's equations
\cite{Gauss_least_constraint:1829}, see also \cite{Lanczos_VPM_1949} for a very interesting discussion. Once we
understand which properties, such as coercivity and stability, yield stable
and reliable ML discretisations, we can select, among the variational
principles describing a problem, one that possesses them, so that its
discretisations inherit them. For simple problems this selection may be
straightforward; for others, such as indefinite, saddle-point problems, or nonlinear phenomena it
becomes a quest for novel variational principles, see
Section~\ref{S:design}. Such a quest has the potential to lead to new
mathematics, compare for example to \cite{DeGiorgi_least_squares:1991}.

A simple example illustrates this theme. For a linear parabolic problem one
may construct time-discrete PINN losses based either on the implicit or on
the explicit Euler method. From the point of view of approximation the two
losses are equally natural and equally consistent, yet only the first
inherits the coercivity of the continuous problem. This is precisely what is
needed for \emph{discrete variational stability}, the stability notion for
discretised functionals introduced in
Section~\ref{S:Abstract-liminf-limsup}, and it allows us to prove convergence
of the corresponding neural-network minimisers. For the explicit loss no such
estimate appears to be available, and numerical experiments show that its
stability deteriorates as the number of spatial training points increases,
in close analogy with the conditional stability of the forward Euler method;
see Section~\ref{SubS:PINN-Parabolic}. This example, together with the problems discussed in
Section~\ref{S:design}, suggests two general principles, which guide the
design of the algorithms developed there. First, the stronger the coercivity properties of the continuous
functional, the better for the algorithm. Second, coercivity at the
continuous level is not sufficient on its own: the fully discrete functional,
including the discretisation in time and space and the training step, must
be designed so as to retain discrete variational stability. It is this
stability, in addition to  consistency, that determines the behaviour of the
algorithm.

Several of the results presented here have appeared in our recent work \cite{loulakis2023newap, GGM_IMA_2025, akrivis2025rungekuttaphysicsinformedneural, grekas2024deepritzfinite, MakrakisMMPP_2025variationalHelm, makr_pryer2024deepUz}; the purpose of these notes is to bring them together within a single framework. Besides the central theme discussed above, two further ideas run through the text. First, the stability of an ML algorithm is best understood as a property of the loss functional, including the training step (quadrature or sampling), rather than of the network parametrisation. Second, empirical measures and optimal transport provide a route to error control that avoids uniform bounds over the network class. Ultimately, this perspective provides concrete links between ML algorithms and their classical numerical counterparts, and suggests that insights from numerical analysis and the calculus of variations can play an important role not only in understanding existing methods, but also in designing new ones.

\subsubsection*{Overview of the notes}

These notes develop a \emph{variational and numerical-analysis perspective on machine learning algorithms}. We view neural networks primarily as nonlinear approximation spaces over which minimisation problems are solved, and investigate how the resulting approximations relate to the underlying continuous problems.

We begin in Sections~\ref{S:typical} and~\ref{S:3:Examples} with a discussion of basic machine learning algorithms and their formulation as discrete minimisation problems. This provides a common framework for several important methods, including function learning with deterministic and probabilistic training, the Deep Ritz method, Physics-Informed Neural Networks, time-discrete training, and generative models. Despite their different computational forms, these methods can often be interpreted through the common structure
$$
\text{continuous problem}
\;\longrightarrow\;
\text{discrete energy}
\;\longrightarrow\;
\text{minimisation over a nonlinear space}.
$$

We then establish, in Section~\ref{S:apriori_apost}, a connection with classical numerical analysis. In particular, we examine the roles of \emph{consistency, stability, approximability, compactness}, and a priori and a posteriori estimates. The central question is whether minimising an approximate energy over increasingly expressive neural-network spaces produces approximations that converge to solutions of the original continuous problem.

The core of the analysis, in Section~\ref{S:Analysis}, is based on a \emph{liminf--limsup framework}, inspired by the calculus of variations and $\Gamma$-convergence. The limsup condition provides an abstract form of consistency and approximability, while the liminf condition captures stability, compactness, and lower semicontinuity. These ideas lead to the notion of \emph{discrete variational stability}, which provides a general framework for establishing convergence of neural-network minimisers.

Having established the basic convergence theory, we turn in Section~\ref{S:reliability} to the question of \emph{reliability and error control}. We consider a posteriori estimates and alternative approaches based on empirical measures, Wasserstein distances, and optimal matching. These provide ways of relating computable losses to the error with respect to the underlying continuous problem.

Finally, in Section~\ref{S:design}, we use the preceding analysis to move from \emph{understanding existing algorithms to designing new ones}, guided by the central theme of the notes: identifying variational principles whose stability properties are inherited by their discretisations. We consider Runge--Kutta Physics-Informed Neural Networks, neural-network methods trained with finite element quadrature, strongly coercive variational formulations for the Helmholtz equation, and deep Uzawa methods, which replace saddle-point problems by sequences of strongly coercive minimisation problems. In particular, the Helmholtz equation illustrates how identifying strongly coercive variational principles can provide a basis for constructing neural-network and finite element methods with improved stability and robustness.
%

Our analysis concerns minimisers of the loss functionals rather than the
dynamics of the optimisation algorithms used to compute them. The
optimisation error is nevertheless partly addressed: for approximate local
minimisers, the error is controlled a posteriori by the computed loss. For
the optimisation dynamics in parameter space, including the convergence of
stochastic gradient descent and the neural tangent kernel and mean-field
regimes, we refer to the recent textbook
\cite{Spiliopoulos_Sowers_Sirignano:2025}. Approximation properties of
neural-network spaces enter as an ingredient and are quoted from the
literature rather than developed.  {We assume familiarity with Sobolev spaces, with the finite element method, and with the basic notions of measure theory and probability; the notions of weak convergence of measures that we need are recalled where they are used.}

\section{Typical machine learning algorithms}\label{S:typical}

\subsection{Machine Learning as Minimisation Problems}

A generic Machine Learning algorithm based on a deep neural network (DNN) can be formulated as a minimisation problem of the form
\begin{equation}\label{mm_nn:intro:abstract_theta}
\min  _ {v \in  V _{\mathcal{N}} } \mathcal{E}_h(v),
\end{equation}
where $V_{\mathcal{N}}$ denotes the class of functions represented by the chosen neural network architecture and $\mathcal{E}_h$ is the corresponding \emph{loss functional}. Here the subscript indicates that $\mathcal{E}_h$ (or $\mathcal{E}_N$ when sampling is used) is a computable approximation of an underlying continuous functional $\mathcal{E}$, obtained for instance by quadrature or sampling; see Section~\ref{S:3:Examples}.  A fundamental feature of this formulation is that $V_{\mathcal{N}}$ is generally \emph{not a linear space}.  
The neural network architecture imposes a nonlinear parametrisation of the admissible functions, and consequently the minimisation problem is, in general, non-convex and may admit multiple local minimisers. The practical training of the neural network amounts to finding an approximate minimiser of this problem, typically by applying an iterative optimisation algorithm to the neural network parameters.

 
%
\subsubsection*{{$V _{\mathcal{N}} \ $}: the discrete Neural Network function spaces }
These are  Finite Dimensional Spaces generated by Neural Networks:   
 The process
 \begin{equation*}
 \mathcal{C}_L:= C_L  \circ \sigma  \circ C_{L-1} \cdots \circ\sigma \circ C_{1} 
\end{equation*}
is in principle a map $\mathcal{C}_L : \mathbb{R} ^m \to \mathbb{R}^{m'} .$
Then the map $\mathcal{C}_L $ is a neural network with $L$ layers and activation function $\sigma.$ 
Any such map $\mathcal{C}_L$ is characterised by the intermediate (hidden) layers $C_k$, which are affine maps of the form 
\begin{equation}\label{C_k}
	C_k y = W_k y +b_k, \qquad \text{where }  W_k \in  \mathbb{R} ^ {d_{k+1}\times d_k}, b_k \in  \mathbb{R} ^ {d_{k+1}}.
\end{equation} 
The set of all networks $\mathcal{C}_L$ with a given structure (fixed $L, d_k,  k=1, \dotsc, L\,$)  
is called $\mathcal{N}.$ The total dimension (total number of degrees of freedom) of  $\mathcal {N} ,$ is $\text{dim } {\mathcal {N}}= \sum _{k=1} ^L d_{k+1} (d_k +1) \, .$ 
 
It is useful to distinguish this formulation from its equivalent representation in the parameter space $\Theta\subseteq\mathbb{R}^{\dim\mathcal{N}}$ (typically $\Theta=\mathbb{R}^{\dim\mathcal{N}}$), where $\theta\in\Theta$ represents collectively all the parameters of the network $\mathcal{C}_L$, namely $W_k, b_k$, $k=1,\dots,L$. The \emph{realisation map} $\mathcal{R}:\Theta\to V_{\mathcal{N}}$, $\theta\mapsto u_\theta$, assigns to each choice of parameters the function computed by the network, so that $V_{\mathcal{N}}=\mathcal{R}(\Theta)$. In the implementation of the algorithm, the actual minimisation which is taking place is with respect to the parameters $\theta$:

\begin{equation}\label{mm_nn:intro:parameter_theta}
	\min _ { \theta \in  \Theta } \mathcal{E}_h (u_\theta).
\end{equation}
This problem \emph{is non-convex} and admits in general more than one local {minimiser}. ML algorithms 
{search} for critical points in the parameter space, i.e., $\theta ^*$ such that 
$$ D_\theta \mathcal{E}_h(u _{\theta ^*}) =0\, .$$
It is to be noted that the corresponding to $u _{\theta ^*}$ Frech\'et derivative on the function space, i.e., the functional $D\mathcal{E}_h(u _{\theta ^*}) $  may not be zero.
Hence,  \cite{Petersen_Raslan_Voigtlaender:2021}, local minimisers in parameter space need not be local minimisers over $V_{\mathcal{N}}$.

Although \eqref{mm_nn:intro:parameter_theta} is the formulation used for computation, our primary perspective is to regard the neural network as defining a nonlinear approximation space and the algorithmic  procedure as an approximate method for solving the corresponding minimisation problem. This distinction is important: the nonlinearity and non-convexity of the parameterisation do not necessarily reflect the analytical properties of the underlying functional $\mathcal{E}_h$. In particular, even though the optimisation problem in parameter space is highly non-convex, the functional may possess useful variational properties, such as coercivity, compactness, or lower semicontinuity, which can be exploited to analyse the resulting approximation.

\subsection{Going from linear to nonlinear}

%
\subsubsection*{Other examples of spaces with finite number of parameters}
While neural network--based function spaces are the primary choice in machine learning algorithms, it is important to note that other choices are meaningful and, depending on the application, may be useful. Let us denote such a space, depending on a finite set of parameters $\Theta$, $V_\Theta$. The corresponding algorithms are then based on
\begin{equation}\label{mm_nn_V_Theta}
	\min_{v \in V_\Theta} \mathcal{E}_h(v).
\end{equation}
A standard choice for approximating physical problems that admit a variational principle is to consider \emph{linear} finite-dimensional subspaces, such as finite element spaces $V_{FE,h}$. In this case, problem~\eqref{mm_nn_V_Theta} with $V_\varTheta = V_{FE,h}$ and a given functional $\mathcal{E}(v)$ is relatively simple to analyse, provided that the underlying continuous problem
\begin{equation}\label{mm_nn_Cont}
	\min_{v \in V} \mathcal{E}(v)
\end{equation}
is well posed and admits a complete theoretical framework. Another class of spaces is to consider a \emph{union} of {(infinitely many)} linear finite element spaces, which for convenience we call \emph{Nonlinear Finite Element} 
spaces, $V_{FE,N}.$ These are nonlinear spaces and their analysis is related to nonlinear approximation theory, although minimisation problems posed over such nonlinear spaces have received comparatively little attention.

\subsubsection*{Discrete energy minimisation problem on parameter spaces:
an example}

 Let $f$ {be} a given function \emph{which at the moment we assume that we know everywhere}. An $L^2$-distance functional will be 
  \begin{equation*}\label{prob_E_error}
 {E} (u      ) =  \frac 1 2 \int_{\Omega}    |  u    - f|^2    \d x \, . \end{equation*}
Let $V_h$ {be} a finite element  (or other finite dimensional) \emph{subspace} of  $L^2(\Omega).$ Consider the approximation problem 
\begin{equation}\label{E-min-FE}
	\min  _ {v \in  V _{h}} E (v)\, .
\end{equation}

 \noindent
 \emph{$L^2$-  projection.\/}
At the minimum, $u_h \in V_h,$ the Frech\' et derivative will vanish, thus 
\begin{equation*}\label{ieE-minimize_NNQ_prob}
 \<  D E (u_h) , \phi \> = 0 \qquad \text {for all } \phi \in V_h\, .
\end{equation*}
But then
\begin{equation*}\label{ieE-minimize_NNQ_prob-d2}
 \<  D E (u_h) , \phi \> = \int_{\Omega}    ( u _h    - f ) \,   \phi \,  dx  \ \qquad \text {for all } \phi \in V_h\, , 
\end{equation*}
i.e. $u_h$ is the $L^2 $ projection of $f$ onto $V_h\, :$
\begin{equation*}\label{ieE-minimize_NNQ_prob-d3}
  \int_{\Omega}    u _h      \,   \phi \,  dx =    \int_{\Omega}     f   \,   \phi \,  dx \ \qquad \text {for all } \phi \in V_h\, . 
\end{equation*}

\noindent
 \emph{Error estimate.\/}
At the minimum, $u_h \in V_h,$ we have trivially:
  \begin{equation*}\label{prob_E_error-d2}
 {E} (u_h      ) \leq  E(\phi ) = \frac 1 2 \int_{\Omega}    |  f    - \phi |^2    \d x  \qquad \text {for all } \phi \in V_h\, . \end{equation*}

Thus, 
\begin{equation*}\label{ieE-minimize_NNQ_prob-d4}
 \| u_h - f \| _{L^2} \leq  \inf _{\phi \in V_h} \| \phi  - f \| _{L^2}   . 
\end{equation*}

\noindent
 \emph{Nonlinear Finite Elements.\/}
  Let $V_{FE, N}$ {be} a nonlinear finite element  \emph{subset} of  $L^2(\Omega).$ Consider the approximation problem 
\begin{equation}\label{E-min-FE-d2}
	\min  _ {v \in  V_{FE, N} } E (v)\, .
\end{equation}
then
    the space $V_{FE, N} $ is \emph{not  linear}    \emph{subspace} of  $L^2({\Omega})\, .$
  	This fact has several important consequences. First, the Euler--Lagrange equations with respect to the function variables do not necessarily hold at the minimisers, since the minimisation is performed over  a nonlinear function space. Second, the corresponding optimisation problem in the parameter space is generally non-convex and may therefore admit multiple local minimisers.
  	
  	As an illustrative example, consider $f \in L^2(0,1)$ and let $V_{FE,N}=V_{h,M}$ denote the space of   piecewise constant functions defined on partitions of $[0,1]$ consisting of $M$ elements. The only restriction on the partitions is the number of subintervals and not the exact location of the nodes that define the partition. In other words, $V_{h,M}$ is the union of all standard  piecewise constant finite element spaces associated with partitions of $[0,1]$ having $M$ elements. This space is clearly nonlinear and provides a prototypical example of a function space studied in the theory of nonlinear approximation.\\
  	

\noindent
 \emph{Nonlinear Finite Elements: error estimate.\/}
 Nevertheless,  
at a {minimiser} $u_{h, N} \in V_{FE, N}$  we still have trivially:
  \begin{equation*}\label{prob_E_error-d3}
 {E} (u_{h, N}     ) \leq  E(\phi ) = \frac 1 2 \int_{\Omega}    |  f    - \phi |^2    \d x  \qquad \text {for all } \phi \in V_{FE, N}\, . \end{equation*}

Thus, 
\begin{equation*}\label{ieE-minimize_NNQ_prob-d5}
 \|u_{h, N} - f \| _{L^2} \leq  \inf _{\phi \in V_{FE, N}} \| \phi  - f \| _{L^2}   . 
\end{equation*}
  
\noindent
 \emph{Do we gain anything?\/}
It is clear that departing from the linear framework introduces several additional computational and theoretical complications to the original approximation problem. In this simple setting, the best approximation result still holds, but the greater flexibility of nonlinear spaces offers the possibility of substantially improved approximation properties. In particular, one expects the nonlinear best approximation error to decay faster than its linear counterpart. Indeed, in the example considered above, free node finite element spaces  have a significantly greater capability to approximate functions with singularities. These issues are discussed in detail in \cite{DeVore_ActaNumer}.
As an illustrative example, consider  $f(x) = x^\alpha,$
 $0<\alpha <1\, .$ In this case the nonlinear best approximation error is
 \begin{equation*}\label{ieE-minimize_NNQ_prob-d6}
 \inf _{\phi \in {V_{h,M}}} \| \phi  - f \| _{ {L^\infty}}  \leq \frac  C {M} \, ,
\end{equation*}
while the same error over a corresponding linear space with $M$ equi-distributed elements is $O( M ^{-\alpha})\, .$ Thus clearly allowing free movement of partitions allows better approximations especially for singular functions.

Although nonlinear approximation spaces offer greater flexibility and potentially much better approximation properties, solving a minimisation problem over such a space is generally much more computationally demanding. A major computational and software development that enabled the rapid expansion of modern machine learning was the ability to solve, in practice, large-scale optimisation problems over nonlinear neural network spaces efficiently. This became possible through the combination of back-propagation and automatic differentiation, stochastic optimisation methods, and powerful hardware and software implementations {~\cite{Baydin_AD_survey:2018, Bottou_Curtis_Nocedal_SIREV:2018}}. As a result, despite the intrinsic nonlinearity and non-convexity of the underlying optimisation problems, neural networks can be trained efficiently even when they involve a very large number of parameters. This computational tractability, combined with the flexibility of nonlinear approximation spaces, makes ML methods particularly attractive in high-dimensional settings, where the performance of most classical approaches tends to deteriorate rapidly as the dimension increases.

\subsection{Training: gradient descent and its variants}\label{SubS:GD}

 The minimisation problems arising in neural-network methods are typically not solved exactly. Instead, one relies on iterative optimisation procedures, with  \emph{gradient descent} and its variations such as \emph{stochastic gradient descent} (SGD) being one of the most widely used. Gradient {descent} has the general form
 $$
\theta_{n+1}=\theta_n-\eta_n D_\theta \mathcal{E}_h(\theta_n),
$$
where $\eta_n>0$ is the learning rate. {It can be viewed as the \emph{explicit Euler} discretisation, with time step $\eta_n$, of the gradient flow $\dot\theta=-D_\theta\mathcal{E}_h(\theta)$ in the parameter space.}

The stationary points of this flow 
 are the critical points 
$$ D_\theta \mathcal{E}_h(u _{\theta ^*}) =0\, .$$
 
This interpretation is more than an analogy. As for the explicit Euler method applied to stiff problems, gradient descent is only conditionally stable: for a quadratic loss it is stable provided $\eta_n$ is small enough, e.g.,   $\eta_n<c /\lambda_{\max}$, where $\lambda_{\max}$ is the largest eigenvalue of the Hessian; compare the conditional stability of the forward Euler method in Section~\ref{SubS:aPaP:parabolic}. Remarkably, full-batch gradient descent on neural networks has been observed to operate typically at this stability limit, in a regime called the \emph{edge of stability} \cite{Cohen_EoS:2021}. Backward error analysis, a classical tool in the numerical analysis of ODE solvers, shows further that the iterates of gradient descent follow more closely the gradient flow of a \emph{modified loss}, containing an implicit regularisation term \cite{Barrett_Dherin_IGR:2021}.
The implicit Euler method, in contrast, can be recovered from  the \emph{proximal point iteration}
$$
\theta_{n+1}\in\operatorname*{argmin}_{\theta\in\Theta}\Big\{\mathcal{E}_h(u_\theta)+\frac{1}{2\eta_n}|\theta-\theta_n|^2\Big\},
$$
which is the parameter-space version of De Giorgi's \emph{minimizing movements}, see also {Section~\ref{SubS:Ex:measures}, where the scheme is used on the space of probability measures.}

In the case where we have a high number of data samples the model is probabilistic and  the loss functional has, in general, the form
$$
\mathcal{E}_N(v)=\frac{1}{N}\sum_{i=1}^N \mathcal{L}(v;X_i),
$$
where $X_1,\ldots,X_N$ are samples and $\mathcal{L}$ denotes the corresponding pointwise or sample-wise contribution to the loss,  see Section \ref{SubS:ModP:function-learning}. Its gradient with respect to the parameters is therefore an average of the gradients associated with the individual samples. Computing the full gradient at every iteration can be expensive, particularly when $N$ is very large. SGD replaces this full gradient by an estimate obtained from a randomly selected subset, or \emph{mini-batch}, of the samples.
More precisely, if $B_n\subset {\{}1,\ldots,N {\}}$ denotes a mini-batch at iteration $n$, the stochastic ``gradient'' is
$$
G_n(\theta_n)
=
\frac{1}{|B_n|}
\sum_{i\in B_n}
D_\theta \mathcal{L}(v_{\theta_n};X_i),
$$
and the parameters are updated according to
$$
\theta_{n+1}=\theta_n-\eta_n G_n(\theta_n),
$$
where $\eta_n>0$ is the learning rate. {If the mini-batches are drawn uniformly at random, $G_n$ is an unbiased estimate of the full empirical gradient.}
 Thus, SGD can be interpreted as a computationally cheaper approximation of gradient descent applied to the empirical variational problem. {It is an instance of stochastic approximation in the sense of 
 \cite{Robbins_Monro:1951}; see \cite{Bottou_Curtis_Nocedal_SIREV:2018} for a review in the context of machine learning.}

From our point of view, the stochastic nature of the algorithm introduces a
further level of approximation: the continuous functional is replaced by an
empirical or discrete one, the minimisation is restricted to a neural-network
space, and the resulting nonconvex problem is solved only approximately.
Separating these sources of error, known in statistical learning as the
estimation, approximation and optimisation errors, allows the properties of
the underlying variational problem to be studied independently of the
optimisation algorithm.

%
%

As mentioned,  convergence of gradient {descent}  should not be confused with convergence of the underlying variational approximation. Even if an optimisation procedure converges to a stationary point of the parameter-space objective, this does not by itself imply that the corresponding neural-network function is close to a minimiser of the continuous problem. Conversely, as we shall see in Section 
\ref{SubS:Conv:MCInterpolation},
the stability properties of the variational problem can provide information about the quality of approximate minimisers without requiring the optimisation problem to be solved exactly. This distinction is particularly important for neural-network methods, where the nonconvexity of the parameterisation makes global optimisation difficult and where practical algorithms generally produce only approximate local minimisers.
This suggests that, rather than seeking a complete theory of global optimisation in the highly nonconvex parameter space, one may instead seek conditions under which the loss produced by the optimisation algorithm is sufficiently small and the underlying variational problem is stable. Which approximate minimiser is reached may nevertheless depend on the algorithm: the \emph{implicit bias} of gradient descent towards particular minimisers is an active topic in machine learning \cite{Soudry_implicit_bias:2018}. We return to this selection mechanism in Section~\ref{SubS:hallucinations}.

\subsection{Hallucinations or numerical artefacts?}\label{SubS:hallucinations}

The term \emph{hallucination} has become widely used in the machine learning literature {~\cite{Ji_hallucination_survey:2023}} to describe outputs that appear plausible but are not supported by the underlying data or model. A related phenomenon can occur when neural networks are used as numerical approximations of mathematical problems: the computed solution may exhibit features that are not present in the exact solution and may nevertheless appear reasonably accurate. This is somewhat different from the standard numerical artefacts encountered in computational science, which are often recognisable from the outset as manifestations of numerical instability, insufficient resolution or noise. In the neural-network setting, however, the computed results may look entirely plausible and even achieve a small training loss, while still leading to unreliable predictions. It is therefore tempting to interpret such behaviour as a ``hallucination'' of the neural network.

However, similar behaviour has been observed for many years in the numerical approximation of nonlinear PDEs, particularly for problems that are sensitive to small perturbations. In such cases, numerical approximations may appear accurate while exhibiting qualitatively incorrect behaviour or converging to the wrong solution altogether. This phenomenon is well understood in the numerical analysis of nonlinear problems,  where it is attributed to the non-uniqueness of weak solutions, combined with instability and with the sensitivity of the underlying problem to perturbations in the data or the approximation scheme. 
 The neural-network setting is heavily nonlinear and further introduces additional sources of uncertainty — the non-convexity of the optimisation landscape and the use of stochastic methods — but the essential difficulty is qualitatively the same. Recognising this connection suggests that the mathematical tools developed for the analysis of nonlinear numerical methods may offer a more principled framework for understanding these phenomena, and it is from this perspective that we approach the discussion that follows.

\subsubsection*{Unreliable predictions and nonlinear PDEs}
\noindent
 \emph{Nonlinear problems with non unique weak solutions.} 
Following,   \cite {Makr_chall_comp_FCGM18}, consider the nonlinear evolution PDE written in abstract form
\begin{equation*}
u_t(t) + A (u(t)) =0 {}{.}
\end{equation*}
This includes nonlinear evolution PDEs such as Conservation Laws (CL), Hamilton--Jacobi equations (HJ), and equations describing phase separation. For many such problems, one is interested in solutions with limited regularity. These \emph{weak solutions} are well defined mathematically and have been extensively studied. However, they are often not unique, requiring additional selection criteria to identify the \emph{physically relevant} solution, such as entropy solutions for CL, viscosity solutions for HJ equations, or appropriate geometric laws for propagating interfaces. Although complete theories are available only for a limited number of important problems, they provide valuable insight into the delicate design of reliable numerical schemes. A key challenge is to ensure that the numerical method does not converge to a mathematically valid but physically irrelevant weak solution. Since such spurious solutions may appear entirely reasonable computationally, additional criteria are needed to assess the reliability of the resulting approximations.
  
\begin{figure}[th]
\begin{center}
			 \includegraphics[scale=0.40]{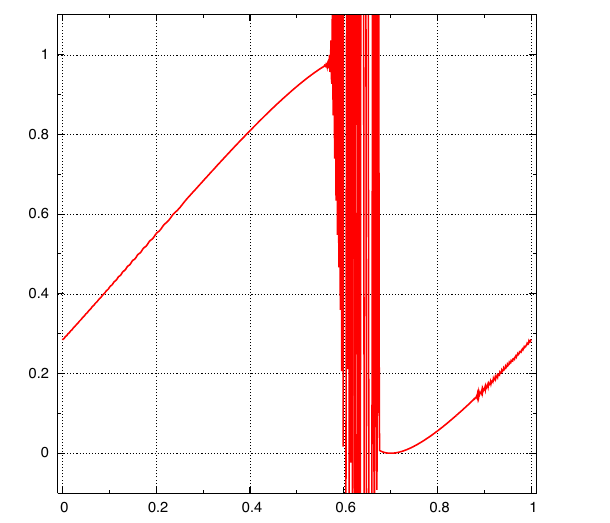}
			 \includegraphics[scale=0.4]{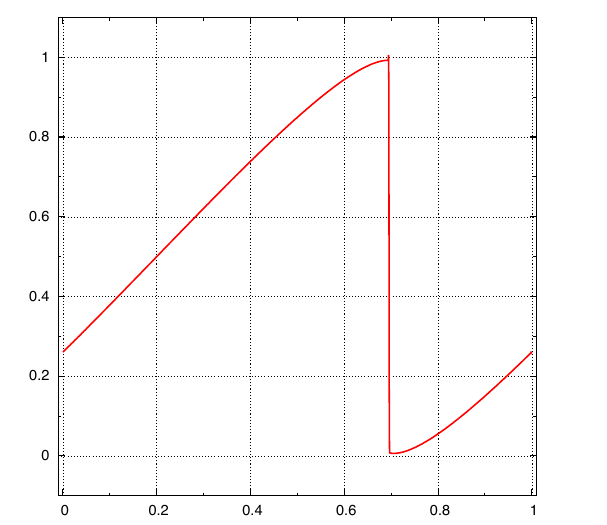}
		 	  		\includegraphics[scale=0.49]{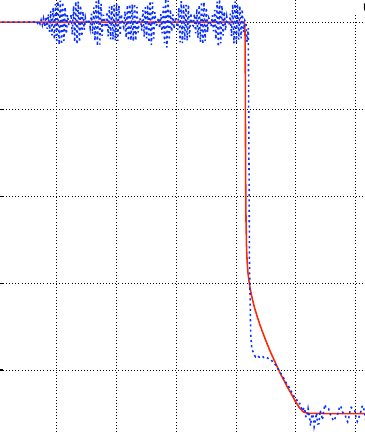}
		\end{center}
	\caption{Central scheme, Upwind scheme, Diffusion and dispersion An oscillatory central scheme (A) and a shock capturing upwind scheme (B) for Burgers equation. Two schemes approximating different weak solutions in the case of non-convex flux function (C){}{.}}\label{fig1}
\end{figure}

  \noindent
\emph{A typical example: Scalar Conservation Laws.} Consider the nonlinear equation 
\begin{equation}\label{1.1}
u_t(x, t) + \text{div} F (u(x, t)) =0, \quad x\in \mathbb{R}^d , t>0.
\end{equation}
We{}{a}k solutions are defined through a variational formulation and admit solutions which are even discontinuous. Among them  uniqueness is guaranteed by imposing an infinite number of additional variational inequalities: The unique \emph{entropy} solution satisfies
\begin{equation*}
\eta (u )_t + \text{div} Q (u ) \leq 0, \quad \text{in }\  \mathcal{D}' . 
\end{equation*}
Here, $\eta$ is   convex  and $\eta, Q$  an entropy entropy-flux pair, \cite{Dafermos_book}. An alternative characterisation of the entropy solution emerges as 
the limit of viscosity approximations (``viscosity solution") 
$
u ^\epsilon \to u 
$
where 
 \begin{equation}\label{1.2}
u_t^\epsilon (x, t) +\text{div} F (u^\epsilon(x, t)) =\epsilon \Delta u^\epsilon(x, t), \quad x\in  \mathbb{R}^d , t>0. 
\end{equation}
 This interpretation of the entropy solution has important consequences for scheme design. Successful approximations typically require some form of \emph{artificial diffusion}, \cite{VonNeumann_Richtmyer_1950}. However, different regularisations may select different limiting solutions, so one must ensure that the scheme does not introduce unwanted perturbations. Thus, even seemingly ``reasonable'' schemes may produce unexpected results. The underlying reason is that each numerical approximation of \eqref{1.1} can be interpreted as a perturbation of the form
 \begin{equation*}
v_t^h ( t) +\text{div} F (v^h( t)) = B_h (h,   v^h( t) ) ,
\end{equation*}
where $ B_h (h,   v^h( t) )$ is a differential operator acting on $ v^h$ {}{which} is not always easy to write explicitly. The character of $B_h$ determines the limiting behaviour: dissipative perturbations may select the entropy solution, whereas dispersive ones may produce persistent oscillations whose weak limits are investigated in  \cite{Goodman_Lax:1988, Lax_Levrmore_1983, Hou_Lax_1991}. 
 Reasonably designed schemes may nevertheless produce remarkably different approximations; see Fig.~\ref{fig1}. Such intriguing computational behaviour dates back to the early development of computational fluid dynamics; see \cite{VonNeumann_coll_works_1963, Hou_Lax_1991}. In particular, Hou and Lax \cite{Hou_Lax_1991} provide a fascinating historical account of J.\ von Neumann's pioneering contributions to scientific computing.

\begin{figure}[th]
 \begin{center}
     \includegraphics[width=3.0cm]{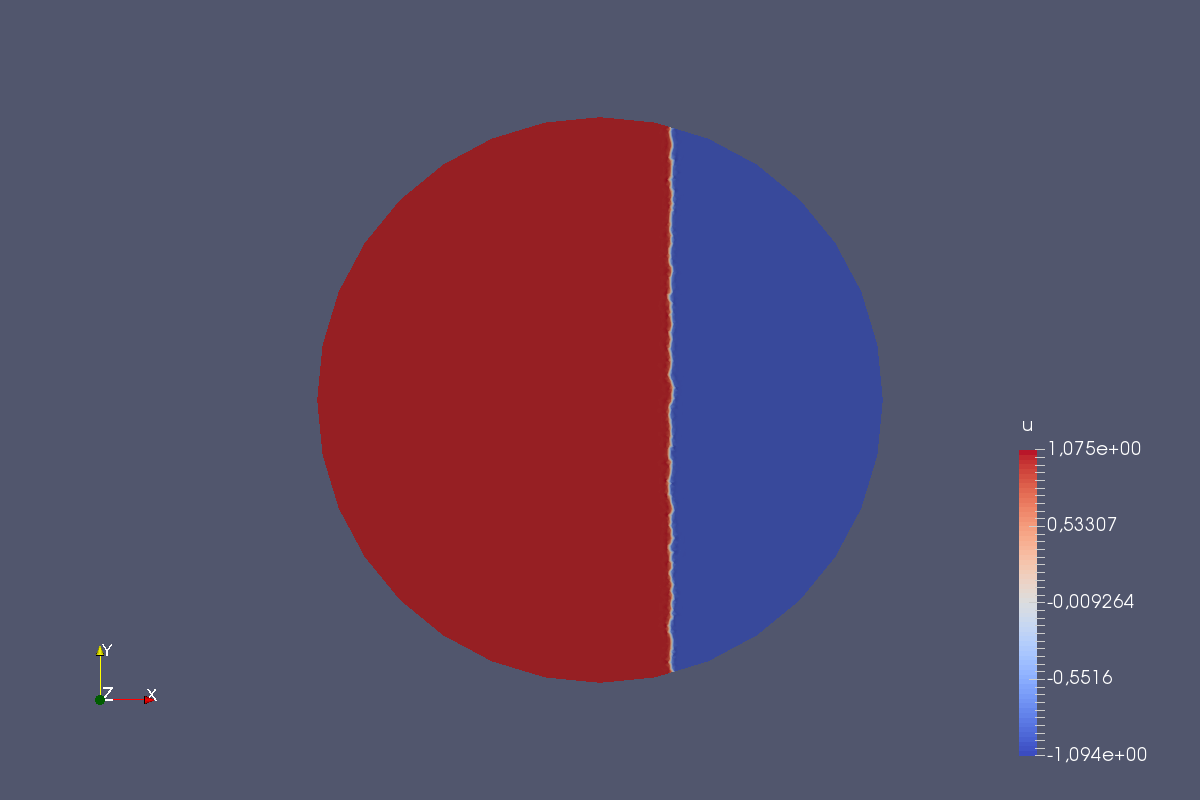}
   \includegraphics[width=3.0cm]{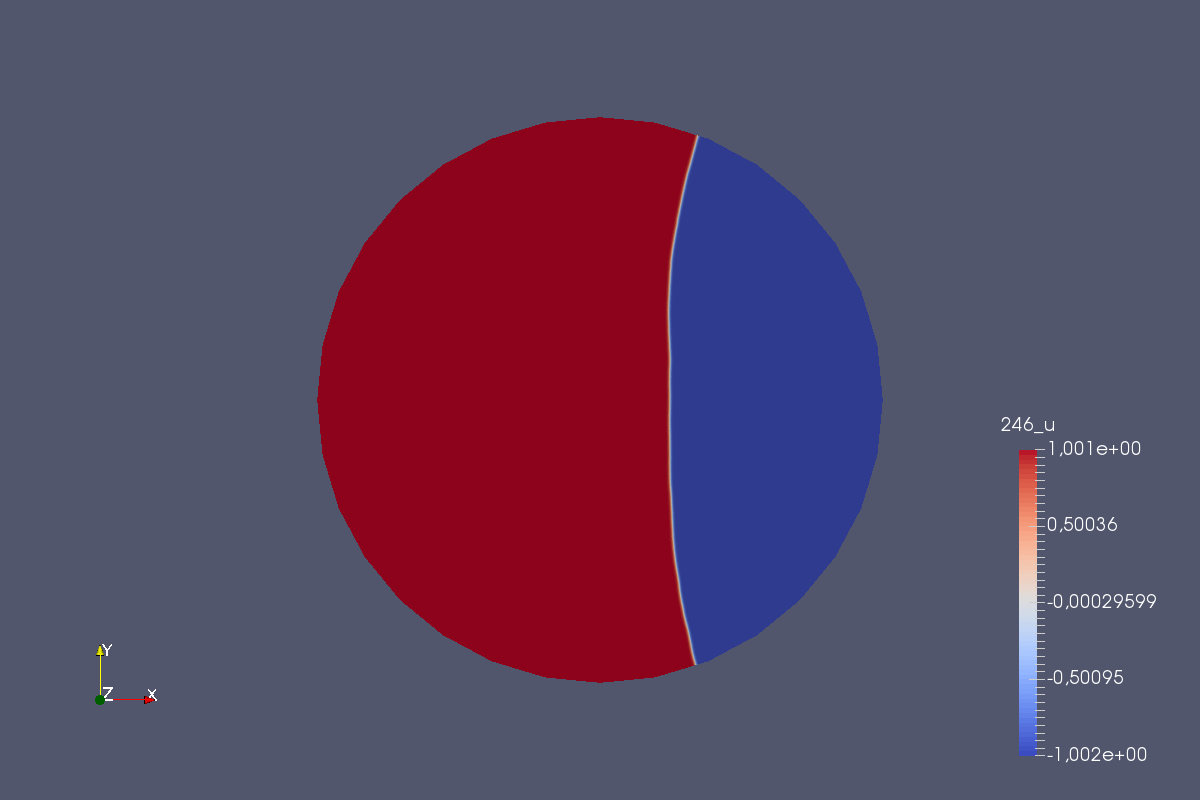}
    \includegraphics[width=3.0cm]{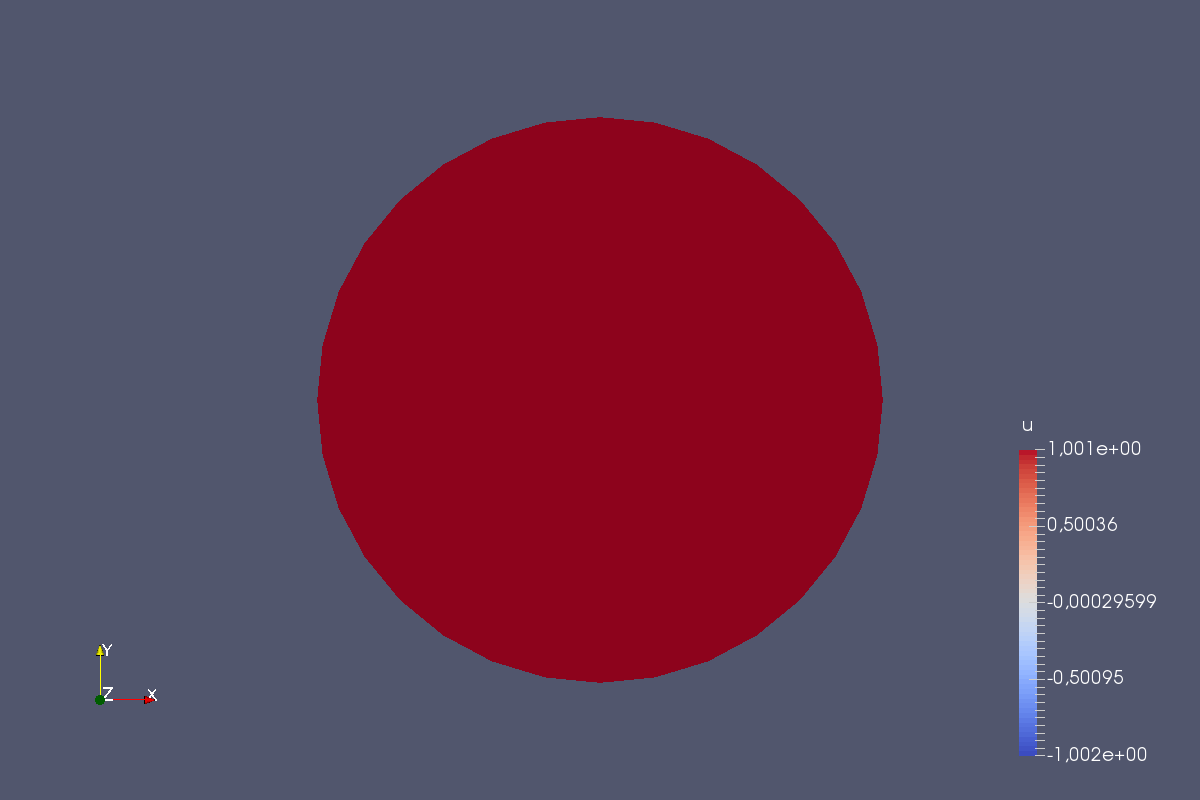}
  \end{center}\
  \caption{Three entirely different approximate minimisers of $\mathcal {E}_{\varepsilon}$ with mesh resolution $h$, $h/2.9$, and $h/3.3$ respectively.}\label{fig_min}
  \end{figure}

\noindent
\emph{Numerical issues on PDE Energy Minimisation.} 
Another important class of \emph{nonlinear problems with non-unique weak solutions} arises in minimisation problems, including phase transitions and other singular phenomena. Although the underlying mathematics is rich, it often provides understanding only in limiting regimes, \cite{Ball_review_2002}. From a computational perspective, 
the numerical analysis of PDEs  energy minimisation is not sufficiently developed.
Consider the nonlinear minimisation problem 
\begin{equation*}
\begin{aligned}
  & \text{find a local minimiser $v$  in $\mathcal{V} $ of : }\\
 & \mathcal {E} (v)   \, {}{,} 
\end{aligned}
\end{equation*}
 The approximate  energy minimisation problem is posed over finite dimensional spaces $V_h,$ and approximating energy functionals $ \mathcal {E}_h\, :$
  \begin{equation*}
\begin{aligned}
  & \text{find a local minimiser $v$  in $\mathcal{V} _h$ of : }\\
 & \mathcal {E}_h (v)   \, . 
\end{aligned}
\end{equation*}
Many questions remain open, since the relationship between $\mathcal{E}_h$ and $\mathcal{E}$ is often unclear. In particular, choosing appropriate approximate energies $\mathcal{E}_h$ can be subtle. Similarly, the relationship between $\mathcal{V}_h$ and $\mathcal{V}$ may lead to unrealistic minimisers, as illustrated by the Lavrentiev phenomenon, see Ball and Knowles \cite{Ball_Knowles_1987}. An indicative example of the sensitivity of discretisation methods is provided by the following simple scalar problem of phase separation with energy functional
\begin{equation*}
\mathcal {E}_{\varepsilon} (u) = \int _\Omega \, \Big [  \, \varepsilon \, |\nabla u |^2 + \, \frac 1 {\varepsilon } (u^2 -1 ) ^2 \Big ]\, dx\, . 
\end{equation*}
Since $\mathcal{E}_\varepsilon(\pm 1)=0$, the constant states $u=1$ and $u=-1$ are global minimisers on any domain. More is true when $\Omega$ is convex: every stable critical point is then constant \cite{Casten_Holland:1978, Matano:1979}, whereas on suitable non-convex domains non-constant local minimisers with interfaces do exist \cite{Kohn_Sternberg_1989}. Nevertheless, the computation of minimisers is far from straightforward: unphysical, but seemingly realistic, interfaces emerge, see Figure~\ref{fig_min}.
 
In view of the above, what is informally described as a ``hallucination'' can often be understood more precisely as a numerical artefact with an identifiable source, and the analogy can be made concrete. A numerical scheme for a nonlinear PDE solves, in effect, a perturbed problem, its modified equation, and the perturbation acts as a hidden selection mechanism among the weak solutions. Similarly, training selects one local minimiser of a non-convex loss among many; as discussed in Section~\ref{SubS:GD}, gradient descent itself approximately follows the gradient flow of a modified loss containing an implicit regulariser. In both cases the computed object is genuine, a weak solution of the equations or a local minimiser of the loss, which explains why it may appear entirely reasonable while being far from the solution sought. 
In both cases, reliability requires information that neither the equations
nor the loss alone provide. For nonlinear conservation laws this is supplied
by the entropy conditions, a selection principle with a long history whose
algorithmic counterparts are themselves delicate. What the corresponding principle should be for the learning problem is
largely open; the stability of the loss functional studied in
Section~\ref{S:Analysis} is, in our view, a necessary ingredient, but
further research is required in this direction.
Rather than treating unexpected features of neural-network approximations as an intrinsic and mysterious property of machine learning, we believe that concepts from numerical analysis, PDEs and the calculus of variations can explain their emergence and provide a principled framework for addressing them.

\subsection{Bibliography and comments}\label{SubS:typical:bib}

Nonlinear approximation, including free-knot splines and adaptive partitions,
is surveyed in \cite{DeVore_ActaNumer}, and the approximation properties of
neural networks, viewed as a form of nonlinear approximation and compared with
classical methods, in \cite{DeVore_Hanin_Petrova:2021}. The topological properties of the sets $V_{\mathcal N}$   are studied in \cite{Petersen_Raslan_Voigtlaender:2021}, and the representation of finite element functions by ReLU networks in \cite{He_Li_Xu_Zheng:2020}. Finite element approximations in which the positions of the nodes are
themselves unknowns appear, often independently, in several areas, including
the study of optimal meshes, moving and adaptive finite elements,
adaptive methods and shape optimisation.

For a review of optimisation methods in machine learning, including
stochastic gradient methods, see \cite{Bottou_Curtis_Nocedal_SIREV:2018}.
Viewing optimisation algorithms as discretisations of gradient flows makes
classical tools of numerical analysis available, such as backward error
analysis, which identifies a modified loss whose gradient flow is followed
more closely by gradient descent \cite{Barrett_Dherin_IGR:2021}, and linear
stability, whose bound $2/\lambda_{\max}$ is observed to govern training in
practice \cite{Cohen_EoS:2021}. The same viewpoint applies to accelerated
methods. Nesterov's method \cite{Nesterov:1983} combines a gradient step with
an extrapolation, or momentum, step,
\[
y_n=\theta_n+\frac{n-1}{n+2}\,(\theta_n-\theta_{n-1}),\qquad
\theta_{n+1}=y_n-\eta\, D_\theta\mathcal{E}_h(y_n),
\]
and, for convex losses with Lipschitz continuous gradient and step size
$\eta\le 1/L$, it reduces the error in the loss at the rate $O(n^{-2})$,
compared with $O(n^{-1})$ for gradient descent; this rate is optimal among
first-order methods. 
It is known that the iterates are appropriate discretisations of  damped second-order equation
\[
\ddot\theta+\frac{3}{t}\,\dot\theta+D_\theta\mathcal{E}_h(\theta)=0,
\]
so that acceleration, too, can be understood through an underlying
continuous dynamical system \cite{Su_Boyd_Candes:2016}. 

The discussion of nonlinear problems with non-unique weak solutions follows \cite{Makr_chall_comp_FCGM18}. For scalar conservation laws, entropy solutions and their uniqueness see \cite{Dafermos_book}; artificial viscosity as a tool in computational science goes back to \cite{VonNeumann_Richtmyer_1950}. The weak limits of dispersive schemes are analysed in \cite{Goodman_Lax:1988, Hou_Lax_1991}, in connection with the small dispersion limit \cite{Lax_Levrmore_1983}, and the selection of nonclassical shocks by numerical schemes in \cite{Lefloch_Rohde_2000}. For energy minimisation, the Lavrentiev phenomenon and its numerical consequences see \cite{Ball_Knowles_1987, Ball_review_2002}; for the absence of stable non-constant equilibria of the Allen--Cahn equation on convex domains see \cite{Casten_Holland:1978, Matano:1979}, and for the existence of non-constant local minimisers on non-convex domains via $\Gamma$-convergence see \cite{Kohn_Sternberg_1989}.
 

\newpage

\section{Key examples and corresponding abstract problems}\label{S:3:Examples}


\subsection{Model Problem 1: Learning a function}\label{SubS:ModP:function-learning}
\smallskip
To motivate the different levels of abstract problems involved in a typical ML algorithm we consider as a model  the learning of a given function $f : \Omega \to \mathbb {R}, $ $\ \Omega \subset \mathbb {R} ^d\, .$ 
We assume that we have access to its point values in a given finite set of points (\emph{interpolation points}).  
Then the \emph{deep neural network interpolant} is the resulting approximation of $f$  which is obtained by 
  a typical machine learning algorithm involving a given DNN architecture and 
  the minimisation step which it is assumed that can be solved exactly. 
 Ideally, one is interested to show convergence of such algorithms as the number of parameters of the DNN and the number of the interpolation points both tend to infinity in an appropriate sense.

\subsubsection*{Loss and minimisers
}

 Considering  $\{y_i\}_i=\{f(x_i)\}_i$ as data, the training loss is    
\begin{equation*}
\label{E_h}
\mathcal{E}_{N}( g )  =     \frac 1 N \sum _{i=1} ^N \, \Big  | g(x_i) - f(x_i) \Big  | ^p  =  \frac 1 N \sum _{i=1} ^N \, \Big  | g(x_i) -y_i\Big  | ^p  \, . 
\end{equation*}
and the learning algorithm is to seek minimisers of this functional, namely,
 \begin{equation}\label{mm_nn:abstract_theta}
		\min  _ {v \in  V _{\mathcal{N}} } \mathcal{E}_N(v)\, .
\end{equation}
We are interested in understanding the behaviour of the minimisers of this problem.
It is interesting first to compare it to standard interpolation in linear finite dimensional spaces. The problem \eqref{mm_nn:abstract_theta} can be seen as \emph{``interpolation" through a minimisation procedure.}
 We are going to see that even this problem is not trivial to analyse. Of course,  the mathematical description of the problem is different depending on whether (i) the point values are given externally (e.g., scattered data)
  or (ii) the point values are chosen by us as it is done in the typical interpolation. To obtain a consistent mathematical theory it is important to emphasise that in the case (i)  the problem is modelled as a probabilistic problem
 while in  the second case (ii) the problem can be modelled as a deterministic problem with a different loss which is consistent with a typical quadrature, see Section~\ref{SubS:Conv:MCInterpolation} for details. 
  
\subsubsection*{One level of abstraction up: NN minimisation without training
}
One may consider the problem \eqref{mm_nn:abstract_theta} as an approximation scheme for the \emph{continuous} minimisation problem 
\begin{equation}\label{mm_nn:abstract_cont}
	\min  _ {v \in  V   } \mathcal{E} (v)  , 
\end{equation}
where $V = L^p (\Omega)$ and
\begin{equation*}
\mathcal{E} (u)=
\int_{\varOmega}    |  u - f|^p    \d x  .
\end{equation*}
The unique minimiser of \eqref{mm_nn:abstract_cont} is obviously $f\, .$
It will be useful to consider as well the $V _{\mathcal{N}}$ minimisation problem \emph{without training:}
\begin{definition}[Semidiscrete idealised problem]\label{abstract_mm_nn} 
Assume that the problem 
\begin{equation*}\label{mm_nn:abstract-d2}
	\min  _ {v \in  V _{\mathcal{N}} } \mathcal{E} (v)   
\end{equation*}
%
has a solution $v^{*} \in V _{\mathcal{N}} .$ We call  $ v^{*} \,  $ a $\mathcal{N}-$DNN minimiser of $\mathcal{E}  $
\end{definition}
As the original  problem, this idealised NN problem is equivalently formulated as a minimisation problem in $\mathbb{R} ^{\text{dim}{\mathcal {N}}}:$
\begin{equation*}\label{mm_nn:abstract_theta-d2}
	\min _ { \theta \in  \Theta } \mathcal{E} (u_\theta).
\end{equation*}
This problem does not lead to a computable algorithm, even if $f$ were known everywhere, since for $v \in V_{\mathcal{N}}$ the quantity $\mathcal{E}(v)$ is not computable. One therefore needs to approximate $\mathcal{E}(v)$ using a quadrature rule or other numerical methods. We will refer to the procedure based on the approximation $\mathcal{E}_h(v) \cong \mathcal{E}(v)$ as \emph{training}.


 \subsubsection*{Quadrature based Deterministic Training}

 We would like to approximate $\mathcal{E}_h(v) \cong \mathcal{E}(v)$ and assume for a moment that  $f$ is known  at a finite set of points 
 $z\in K_h,$ where $K_h$ is a set of points linked to a quadrature rule.  
 We aim   to construct an approximation of $f$ everywhere in its domain.  
We consider    appropriate quadrature for integrals over $\varOmega.$ Such a quadrature requires 
 a set $K_h$ of discrete points $z\in K_h$  and corresponding nonnegative weights $w_z$
 such that 
 \begin{equation}\label{quadr_lf}
\sum _{z\in K_h} \, w_z \, G(z) \cong \int _{\Omega} \, G(x) \, \d x .
\end{equation}
Based on \eqref{quadr_lf}
 we  obtain the ``loss" functional approximating $\mathcal{E}(v)$
\begin{equation*}
\label{E_h-d2}
\mathcal{E}_{Q, h}( g )  =     \sum _{z\in K_h} \, w_z \, \Big  | g(z) - f(z) \Big  | ^p  \, . 
\end{equation*}

\begin{definition} \label{D-dnnint}{\sc [Deep neural network interpolation]}  
Assume that the problem 
\begin{equation*}\label{ieE-minimize_NNQ}
	\min  _ {v \in  V _{\mathcal{N}} } \mathcal{E} _ {Q, h}(v)
\end{equation*}
has a solution $u_h \in V _{\mathcal{N}} .$ We call  $ u_h \,  $ a $\mathcal{N}-$\emph{deep neural network interpolant of $f$ at $K_h  \, .$}
\end{definition}

%
%

%
%
%
%
%
%
%
%
%
%

\subsubsection*{Probabilistic Training} 
 
When we wish to use a data set $\{y_i\}_i = \{f(x_i)\}_i$ to learn $f$, we may not be able to choose the locations of the points $\{x_i\}_i$. In this case, the loss
\begin{equation*}
    \frac{1}{N} \sum_{i=1}^N \Big| g(x_i) - y_i \Big|^p ,
\end{equation*}
can only be regarded as a \emph{Monte Carlo} approximation of $\mathcal{E}(v),$ i.e., the integral is   approximated by
\begin{equation}
\label{quadrature}
\int_{\Omega} G(x)\,{\d\mu(x)} \cong \frac{1}{N} \sum_{i=1}^N G(x_i),
\end{equation}
{where the $x_i$ are random points in $\Omega$, distributed according to a probability measure $\mu$ on $\Omega$. Unless stated otherwise, $\mu$ is the uniform distribution on $\Omega$, so that $\int_\Omega G\,\d\mu=|\Omega|^{-1}\int_\Omega G\,\d x$; the constant factor $|\Omega|$ is irrelevant for the minimisation, and in the probabilistic setting all integrals over $\Omega$ are understood with respect to $\mu$.}


To be more precise, for a sample $\omega \in \bbOmega$ fixed, and a collection $  X_1,$ $X_2,\ldots$ of i.i.d.\ $\Omega$-valued random variables, defined on a probability space $(\bbOmega,\mathcal{F},\mathbb{P})$
\footnote{We consider a collection $Y, X_1,$ $X_2,\ldots$ of i.i.d.\ $\Omega$-valued random variables, defined on a probability space $(\bbOmega,\mathcal{F},\mathbb{P})$, with common law $\mu.$ 
The role of $X_1,X_2,\ldots$  is that they model the random selection of the points $x_1, x_2, \ldots .$ 
For $g\in {L^p(\Omega,\mu)}$ we have 
\begin{equation}\label{repr_one}
	\int _\Omega \, g (x) \, \d \mu (x) = \int _\bbOmega \, g ( X_i (\omega)) \, \d \PP (\omega) = \int _\bbOmega \, g ( Y (\omega)) \, \d \PP (\omega)= \EE [ \, g(Y) \, ]\, .
\end{equation} 
 } 
we consider the discrete energy, 
\begin{equation}\label{prob_E}
\mathcal{E}_{N, \omega}( g )  =  \frac 1N  \sum _{i=1}^N   \Big  | g(X_i(\omega)\, ) - f(X_i (\omega) \, ) \Big  | ^p  \, . 
\end{equation}
This will be the loss approximating $\mathcal{E}\, .$ 
This leads to the problem
\begin{definition} \label{D-dnnint-d2}{\sc [probabilistic deep neural network interpolation]}  
Assume that for each $\omega \in \bbOmega$ the problem 
\begin{equation}\label{prob_E_def}
	\min  _ {v \in  V _{\mathcal{N}} } \mathcal{E} _{N, \omega} (v)
\end{equation}
has a solution $u ^\star(\omega, \cdot)   \in V _{\mathcal{N}} .$  We denote by   $ u_{\bm{X} (\omega)   }  = u ^\star(\omega, \cdot)  \,  $ the  $V_{\mathcal{N}}-$\emph{ probabilistic deep neural network interpolant of $f$ at $ \bm{X}  = (X_1, \dots, X_N  )  \, .$}
\end{definition}
 We assume that the minimiser can be selected so that $\omega\mapsto u^\star(\omega,\cdot)$ is measurable, which is the case, for instance, when the parameters range over a compact set and the loss depends continuously on them.

A key challenge is the following. For a fixed function $g$, the quadrature error
\begin{equation*}
\frac 1 N \sum _{i=1} ^N \, g(X_i) - \int _{\Omega} \, g(x) \, \d \mu(x)
\end{equation*}
can be controlled by standard tools, such as the central limit theorem or large deviations. These tools do not apply directly to
\begin{equation*}
\frac 1 N \sum _{i=1} ^N \, u_{\bm X}(X_i) - \int _{\Omega} \, u_{\bm X}(x) \, \d \mu(x) ,
\end{equation*}
since the interpolant $u_{\bm X}$ depends on the sample points themselves. In statistical learning this difficulty is addressed by uniform laws of large numbers over the class $V_{\mathcal{N}}$, through Rademacher- or VC-type complexity bounds; an alternative approach, based on the convergence of empirical measures, is developed in Sections~\ref{S:Analysis} and~\ref{S:reliability}.

\subsubsection*{Function learning in infinite dimensions: operator learning}

The formulation of Model Problem 1 can be extended  to infinite dimensional domains instead  of
subsets of $\mathbb{R}^d$. Let $X$ and $Y$ be separable Banach spaces of
functions, for instance $X=L^2(D)$ or $X=C(\overline D)$ and $Y=H^1_0(D)$ for a
bounded domain $D\subset\mathbb{R}^d$, and let $\mathcal{G}:X\to Y$ be a given
operator. A typical example is the solution operator of a parametric PDE,
such as $a\mapsto u_a$, where $u_a$ solves $-\operatorname{div}(a\nabla u_a)=f$
in $D$ with homogeneous Dirichlet conditions, in which case $X$ is a space of
coefficient fields. The learning of $\mathcal{G}$, known as \emph{operator
learning}, is function learning with $\Omega$ replaced by $X$: we fix a
probability measure $\mu$ on $X$, which models the inputs of interest (for
example a Gaussian measure, or the law of random coefficient fields), and
consider the continuous functional
\begin{equation*}
\mathcal{E}(G)=\int_X \|G(a)-\mathcal{G}(a)\|_Y^p\,\d\mu(a),
\qquad G\in V:=L^p(X,\mu;Y),
\end{equation*}
whose unique minimiser is $\mathcal{G}$. The role of $V_{\mathcal{N}}$ is
played by a class of \emph{neural operators}, maps $G_\theta:X\to Y$ of a
prescribed architecture, such as deep operator networks or Fourier neural
operators, see e.g.,   \cite{lu2021learning, Mishra_DeepONET_Kran_est, Stuart_NeuralOper2023, Mishra_convNO:2023}. Every such architecture
combines an finite dimensional approximation  of the input $a$ (encoding into finitely many numbers 
such as point values
on a grid, Fourier or other coefficients), a neural network acting on these
coefficients, and a decoding of the output into a function in $Y$. The
semidiscrete problem is $\min_{G\in V_{\mathcal{N}}}\mathcal{E}(G)$, and
training is based on data $(a_i,\mathcal{G}(a_i))$, $i=1,\dots,N$, with
$a_i=X_i(\omega)$ i.i.d.\ samples from $\mu$, through the loss
\begin{equation*}
\mathcal{E}_{N,\omega}(G)=\frac1N\sum_{i=1}^N
\|G(X_i(\omega))-\mathcal{G}(X_i(\omega))\|_Y^p .
\end{equation*}
In practice the data $\mathcal{G}(a_i)$ are computed by a numerical solver and
the norm of $Y$ is evaluated by quadrature on $D$, so that a second level of
discretisation is present; we do not pursue this here.

One feature distinguishes operator learning from the finite-dimensional case.
Since there are no quadrature rules on an infinite-dimensional space, training
is necessarily probabilistic.
When the inputs depend on
finitely many parameters, as for parametric PDEs with $a=a(x;\xi)$,
$\xi\in\mathbb{R}^m$, the problem reduces to function learning on
$\mathbb{R}^m$, and deterministic training becomes possible again.

\subsection{Model Problem 2: Solving elliptic   PDEs  
}\label{SubS:Ex:EllipticPDE}
\smallskip
 
We consider  simple  boundary value problems  of the form,
\begin{equation}\label{EllipticPDE}
\left \{
\begin{alignedat}{3}
&L\, u = f\quad &&\text{in}\,\,&&\varOmega \\[2pt]
& u =0 \quad &&\text{on}\,\,&&\partial\varOmega \\[2pt]
\end{alignedat}
\right .
\end{equation}
 where $ u : \Omega \subset \mathbb{R}^d \rightarrow \mathbb{R} , \: \Omega $ is an open, bounded set with Lipschitz  boundary, 
$f \in L^2 (\Omega)  $ and $ L $ a self-adjoint elliptic operator of the form
$Lu := - \sum_{1 \leq i,j \leq d} \big ( a_{ij} u_{x_i } \big )_{  x_j}  +cu , $
$c>0, $ and where $ a_{ij}$ are smooth enough and 
satisfy standard ellipticity and positivity assumptions.  More general elliptic operators and alternative boundary conditions can be considered as well. There are two main approaches to discretising this problem, motivated by either minimising the Dirichlet energy (the Deep Ritz method) or minimising the  $L^2$ residual (PINN methods). As in the case of function learning, these neural network methods can ultimately be formulated by introducing an appropriate quadrature for the corresponding functionals and seeking minimisers of the resulting loss functionals over suitable neural network spaces.

\subsubsection*{Deep Ritz method  
}

The \emph{(idealised without training/quadrature)} Deep Ritz  method corresponds to the loss
\begin{equation*}\label{en-func-DR} 
\mathcal{E}_D (u)=  \frac12 B(u, u)  - \int_{\Omega}  f\, u  \d x  \, ;
\end{equation*}
where 
$$B(u, v) : = \int_{\Omega}  a(u, v) \d x\, : = \int_{\Omega}  \Big(  \sum_{1 \leq i,j \leq d}  a_{ij} u_{x_i } v _{  x_j}  +cu \, v    \Big) \d x\,  , $$
is the standard bilinear form corresponding to $L.$ {Here $$a(u,v)(x):=\sum_{1\le i,j\le d} a_{ij}(x)\,u_{x_i}(x)\,v_{x_j}(x)+c(x)\,u(x)\,v(x)$$ denotes the pointwise form, so that $B(u,v)=\int_\Omega a(u,v)(x)\,\d x$.}
When $L$ is self-adjoint 
   the   solution of \eqref{EllipticPDE} is the unique minimiser of the problem 
\begin{equation}\label{Dirichlet_Energy}
	\min  _ {v \in  H^1 _0 (\Omega ) } \mathcal{E}_D (v)\, .
\end{equation}
%
%
%

As before, to define the loss functional one may approximate integrals   through Monte-Carlo integration.  
 If $X_i,$ $i=1, \dots N, $ i.i.d.\ random variables, each one corresponding to a random point $x_i$ of $\Omega.$ The integrals are approximated, up to the constant factor $|\Omega|$, which is irrelevant for the minimisation, by 
 \begin{equation*}
\label{quadrature-d2}
\frac 1 N \sum _{i=1} ^N \,  \, F(x_i) \approx \int _{\Omega} \, F(x) \, \d x .
\end{equation*}
and the \emph{Deep Ritz method with probabilistic training} is based on the loss functional
\begin{equation*}
\label{E_h-d3}
 \mathcal{E}_{D, N, \omega }( g )  =  \frac 1N  \sum _{i=1 } ^N    \, {\frac 12\, a(g,g)(X_i(\omega))}    - f(X_i(\omega)\, ) g(X_i(\omega)\, )  \, ,
\end{equation*}
leading to the discrete minimisation problem  (for each instance $\omega \in \bbOmega$),
\begin{equation*}\label{ieE-minimize_NNQ-d2}
	\min  _ {v \in  V _{\mathcal{N}} }  \mathcal{E}_{D, N, \omega }(v)\, . 
\end{equation*}
%
Notice that the requirement $V_{\mathcal{N}} \subset H^1_0(\Omega)$ is subtle and may, in fact, lead to computational complications. It is assumed here for simplicity of the exposition.
 Notice also that, in contrast to the continuous energy, the sampled energy
$\mathcal{E}_{D, N, \omega}$ controls $g$ and $\nabla g$ only at the sample
points, and this control depends on the quadrature, that is, on the sample
points and their number. Nothing in the loss constrains $g$ outside the
sample points, so the coercivity of $\mathcal{E}_D$ on $H^1_0(\Omega)$ is not
automatically inherited by $\mathcal{E}_{D, N, \omega}$ on $V_{\mathcal{N}}$.
Whether, and in which sense, a discrete loss retains the coercivity of the
continuous functional is precisely the question of discrete variational
stability addressed in Section~\ref{S:Analysis} and, for Deep Ritz methods
trained with finite elements, in Section~\ref{SubS:FE-training}.

 Alternatively, one may employ quadrature-collocation training
to approximate the integrals, yielding different loss functionals.   Utilising  a quadrature   
 over a set $K_h$ of discrete points $z\in K_h$  and corresponding   weights $w_z$
 such that 
 $
\sum _{z\in K_h} \, w_z \, F(z) \approx \int _{\Omega} \, F(x) \, \d x .
$
  we may  define
%
 %
\begin{equation}
\label{E_DQh}
\mathcal{E}_{D, Q, h}( g )  = \sum _{z\in K_h} \, w_z \,\Big (\frac 1 2 \,  {a(g,g)(z)} - f(z) g(z)\,\Big ) \, . 
\end{equation}
Then the discrete minimisation problem  becomes
\begin{equation*}\label{ieE-minimize_NNQ-d3}
	\min  _ {v \in  V _{\mathcal{N}} } \mathcal{E}_{D, Q, h} (v)\, .
\end{equation*}


\subsubsection*{PINN methods  
}
  
Although methods based on the discretisation of the Dirichlet energy have a solid physical foundation, they can only be used for elliptic problems with self-adjoint operators $L$. Residual-based methods are   more general and are based on the minimisation of the functional
\begin{equation*}\label{mm_nn:abstract-d3}
	\min  _ {v \in  V _{\mathcal{N}} } \mathcal{E}_R (v)   
\end{equation*}
where
\begin{equation*}\label{Functional_L}
\mathcal{E}_R (v) = \int_{\Omega} | L v  - f |^2 dx \, . 
\end{equation*}
The Monte Carlo approximation of the integral yields loss functionals whose minimisation leads to what are commonly referred to as \emph{Physics-Informed Neural Networks (PINNs)}.
To this end, by employing the same notation used previously
one may define the loss functional 
\begin{equation*}
\label{E_h_PINN}
 \mathcal{E}_{R, N, \omega }( g )  =  \frac 1N  \sum _{i=1 } ^N    \, \Big | {(L g)(X_i(\omega))} -        f(X_i(\omega)\, )   \Big |^2  \, ,
\end{equation*}
and the PINN method is 
\begin{equation*}\label{ieE-minimize_NNQ-d4}
	\min  _ {v \in  V _{\mathcal{N}} } \mathcal{E}_{R, N, \omega } (v)\, \, . 
\end{equation*}
Here we need to select neural network spaces $V_{\mathcal{N}}$ with higher regularity. The direct imposition of boundary conditions then becomes more intricate, especially in Lipschitz domains, where singular solutions may arise.  In the definition of $\mathcal{E}_R$ we have assumed that the functions in $V_{\mathcal{N}}$ satisfy the boundary condition exactly, so that $V_{\mathcal{N}}\subset H^2(\Omega)\cap H^1_0(\Omega)$; otherwise an appropriate discretised  boundary term, such as $\int_{\partial\Omega}|v|^2\,\d S$ or a stronger norm of $v$, has to be added to the loss.

The two approaches illustrate a central theme of these notes: the same
boundary value problem \eqref{EllipticPDE} admits several variational
principles with different mathematical properties. Which of the two leads to
stable neural network discretisations, and under which training, is to be
investigated. Elliptic problems in this setting are quite simple; however,
considering alternative variational principles for more complicated learning
tasks might have a significant influence on the behaviour of the algorithm.

\subsection{Model Problem 3: Linear Parabolic case 
}\label{SubS:Ex:parabolic}
 Let as before $ \Omega \subset \mathbb{R}^d $, open, bounded and set $ \Omega_T = \Omega \times (0,T] $ for some fixed time $ T>0. $
We consider the parabolic problem

\begin{align}\label{ParabolicPDE}
\begin{cases}
u_t + Lu = f ,\;\;\; \; \textrm{in} \;\: \Omega_T, \\
u= 0, \;\;\;\;\;\;\;\;\; \; \textrm{on} \;\: \partial \Omega \times (0,T],  \\
u = u^0, \;\;\;\;\;\;\;\;\, \textrm{on} \;\: \Omega \times \lbrace t=0 \rbrace \, .
\end{cases}
\end{align}
where $L$ is an elliptic operator as in Section \ref{SubS:Ex:EllipticPDE}.  In this section we formulate approximations of \eqref{ParabolicPDE} obtained by minimisation of continuous and time-discrete energy functionals over appropriate sets of neural network functions; their convergence properties are studied in Section~\ref{SubS:PINN-Parabolic}. 

For simplicity in the exposition we shall assume that $\Omega$ is a convex Lipschitz domain; the general case of a non-convex domain can be treated as well. 

\subsubsection*{Exact time integrals and probabilistic training}
So now we define $ \mathcal{G} : H^1(0,T ; L^2 (\Omega)) \cap L^2(0,T ; H^2(\Omega)  \cap  H^1_0(\Omega))\rightarrow \overline{\mathbb{R}} $ as follows
\begin{equation}\label{ParabFunctional}
\begin{gathered}
\mathcal{G}(v) =
\int_0^T \| v_t(t) + L v(t) - f(t) \|^2_{L^2(\Omega)} \d t +  | v(0) - u^0  |^2_{H^1(\Omega)}\, .
\end{gathered}
\end{equation}
Treating boundary conditions is typically complex; since we would like to focus on the general ideas, for simplicity we select neural network spaces that exactly satisfy the boundary conditions, see e.g., \cite{sukumar2022exact}. The space-time versions of the neural spaces 
$
V _{\mathcal{N}, 0} \subset H^1(0,T ; L^2 (\Omega)) \cap L^2(0,T ; H^2(\Omega) \cap H^1_0(\Omega))$ 
are still denoted by $V _{\mathcal{N}}$ and $V _{\mathcal{N}, 0}$ (the $0$ in the latter indicates that the boundary conditions are enforced strongly). 
As in the elliptic case, residual-based methods are then   based on the minimisation of the functional
\begin{equation*}\label{mm_nn:abstract-d4}
	\min  _ {v \in  V _{\mathcal{N}, 0} } \mathcal{G} (v) \, .  
\end{equation*}
The standard PINN method requires a Monte Carlo approximation of the space-time integral:
\begin{equation*}
\label{E_h_PINN-d2}
 \mathcal{G}_{\Omega_T, N, \omega }( g )  =  \frac 1N  \sum _{i=1 } ^N    \, \Big |g_t (\overline X_i(\omega)) +  {(L g)(\overline X_i(\omega))} -        f(\overline  X_i(\omega)\, )   \Big |^2  \, ,
\end{equation*}
and of the spatial integral corresponding to the initial condition
\begin{equation*}
\label{E_h_PINN-d3}
 \mathcal{G}_{\Omega, N, \omega }( g )  =  \frac 1N  \sum _{i=1 } ^N    \, \Big |\nabla g (  X_i(\omega) {, 0}) - \nabla u^0(   X_i(\omega)\, )   \Big |^2  \, ,
\end{equation*}
 where $\overline X_i(\omega)$ correspond to random points of the space-time domain $ \Omega_T = \Omega \times (0,T] $ and $ X_i(\omega)$  correspond to random points of the spatial domain $\Omega.$
 Then,  denoting $ \mathcal{G}_{R, N, \omega } (g) : =  \mathcal{G}_{\Omega_T, N, \omega } (g)+ \mathcal{G}_{\Omega, N, \omega } (g)$ we seek minimisers of the problem
\begin{equation*}\label{ieE-minimize_NNQ-d5}
	\min  _ {v \in  V _{\mathcal{N}, 0} } \mathcal{G}_{R, N, \omega } (v)\, \, . 
\end{equation*}
  

 \subsubsection*{Time discrete Training
}
 
Applying quadrature only to the time integral yields methods based on \emph{time-discrete training.} This framework is particularly well suited both for analysing the stability of the method and for developing alternative algorithmic designs.  
Let $0=t^0< t^1< \cdots <t^N=T$ define a partition of $[0,T]$ and $I_n:=(t^{n-1},t^n],$ $k_n:=t^n-t^{n-1}.$
We shall denote by $v^m(\cdot)$ and $f^m(\cdot)$ the values $v (\cdot,t^m)$ and $f(\cdot ,t^m).$ Then we define 
the discrete in time quadrature by 
 \begin{equation*}
\label{quadrature_time}
\sum _{n=1}^N \,  k_n \, g(t^{n}) \approx \int _{0}^T \, g(t) \, \d t.
\end{equation*}
We then define the time-discrete version of the functional \eqref{ParabFunctional} as  
\begin{equation}\label{Functional_k:intro}
\mathcal{G}_{k, IE} (v) = \sum _{n=1}^N \,  k_n \, \int_{\Omega }\big | \frac{v ^n-  v^{n-1}}{k_n} + L v ^n  - f^n \big |^2 \, \, \d x +\, 
|v {^0}  - u^0 |_{H^1 (\Omega)}^2 
\end{equation}
Then the corresponding minimisation problem over $  V _{\mathcal{N}, 0}$
%
\begin{equation}\label{ieE-minimize_k}
	\min  _ {v \in  V _{\mathcal{N}, 0} } \mathcal{G} _{k, IE}(v)\, ,
\end{equation}
  forms the core algorithmic component of the time-discrete training approach. Fully discrete formulations are then obtained by approximating the spatial integrals with suitable quadrature rules, which may be either probabilistic or deterministic. 
%
%

{From the approximation point of view, an equally natural choice is the discrete functional}
\begin{equation}\label{Functional_k_Ex_intro}
\mathcal{G}_{k, EE} (v) = \sum _{n=1}^N \,  k_n \, \int_{\Omega }\big | \frac{v ^n-  v^{n-1}}{k_n} + L v ^{n-1}  - f^{n-1} \big |^2 \, \, \d x +\, |v {^0}  - u^0 |_{H^1 (\Omega)}^2\, .
\end{equation}
{In Section~\ref{SubS:PINN-Parabolic} we compare its properties with those of the functional $\mathcal{G}_{k, IE}$ and the associated minimisers over $V_{\mathcal{N},0}$.}
This comparison demonstrates that the stability of the algorithms is intrinsically linked to the coercivity of the loss functional.
%
 

\subsection{Generating probability measures}\label{SubS:Ex:measures}

One of the key successes of neural network methods is their remarkable ability to effectively represent and approximate probability distributions. This capability lies at the heart of many modern machine learning applications, where the objective is not merely to approximate a deterministic function but rather to learn the underlying probability distribution from which the observed data are drawn. In many problems of practical interest, the available data consist only of a finite random sample, viewed as specific realisations of an unknown probability distribution. The goal is then to infer or generate this underlying distribution in a way that faithfully captures its essential statistical properties and can be used to produce a new sample with similar characteristics.
To be more precise, suppose that the data is a   random sample $X_1,\dots,X_N$ of an unknown probability measure $\nu_{\, \text{true}}\in\mathcal{P}(\mathbb{R}^k)$, represented by the empirical measure
\begin{equation}\label{emp_meas}
	\nu_{\, \text{data}} =   \frac 1 N \sum _{i=1}^N \   \delta _{X_i } \, .
\end{equation}
Following the pattern of the previous model problems, we distinguish three levels. The \emph{continuous} problem is
\begin{equation*}
	  \min_{ \mu  \in  \mathcal{P}(\mathbb{R}^k) } d_{\mathscr {M}} ({\mu ,\nu _{\, \text{true}} }),
\end{equation*}
where $d_{\mathscr {M}}({\cdot ,\cdot })$ is a specified distance (or divergence) on $\mathcal{P}(\mathbb{R}^k)$; its unique minimiser is $\nu_{\, \text{true}}$ itself. The \emph{semidiscrete} problem restricts the minimisation to a \emph{neural measure space} $V_\Theta\subset\mathcal{P}(\mathbb{R}^k)$, a parametrised family of measures determined by neural networks of a given architecture, defined below:
\begin{equation*}
	  \min_{ \mu  \in  V  _{\Theta } } d_{\mathscr {M}} ({\mu ,\nu _{\, \text{true}} }).
\end{equation*}
Finally, since $\nu_{\, \text{true}}$ is known only through the data, \emph{training} replaces it by $\nu_{\, \text{data}}$, 
in the sense of Section~\ref{SubS:ModP:function-learning}:
\begin{equation}\label{mm_nn:abstract}
	  \min_{ \mu  \in  V  _{\Theta } } d_{\mathscr {M}} ({\mu ,\nu _{\, \text{data}} }),
\end{equation}
and we hope that
\begin{equation}\label{mm_nn:abstract:2}
\mu_{\theta ^*} \cong \nu _{\, \text{true}}, \qquad \text {where } \mu_{\theta ^*}
\in \text{argmin} _ {\mu \in V _{\Theta } } d_{\mathscr {M}}({\mu ,\nu _{\, \text{data}} } ).
\end{equation}
In practice $d_{\mathscr {M}}(\mu_\theta,\nu_{\, \text{data}})$ is itself evaluated through samples of $\mu_\theta$, so that a further empirical approximation is present; we do not pursue this here.

We note that over the
  whole of $\mathcal{P}(\mathbb{R}^k)$ the trained problem is minimised by $\nu_{\, \text{data}}$ itself, an atomic measure which reproduces the data. 
  The ability of the model to generalise therefore rests 
  on the restriction to $V_\Theta$, and the relevant question is, as before, the behaviour of the minimisers as $N\to\infty$. 


We refer to the parameterised families $V_\Theta \subset \mathcal{P}(\mathbb{R}^k)$, induced by suitable classes of   neural networks, as neural measure spaces. As we shall see, although this formulation provides a natural framework for probability learning, it gives rise to a broad spectrum of methods, depending on two fundamental choices: the distance used to compare probability measures and the choice of the neural measure space itself.

In all cases considered here, the measures in $V_\Theta$ are generated from a fixed, known reference measure $\gamma_0$ on some $\mathbb{R}^m$, for instance a standard Gaussian, which is easy to sample. In the simplest case,
\begin{equation*}
V_\Theta=\big\{(T_\theta)_\sharp\gamma_0:\ \theta\in\Theta\big\},
\end{equation*}
where $T_\theta:\mathbb{R}^m\to\mathbb{R}^k$ is a neural network and $(T_\theta)_\sharp\gamma_0$ denotes the push-forward of $\gamma_0$ through $T_\theta$, defined below; new samples are then produced by sampling $\gamma_0$ and applying $T_\theta$. We call such spaces \emph{static} neural measure spaces.  Of particular interest for these notes are the \emph{flow-based} neural
measure spaces, generated by an evolution process on $[0,T]$ started from
$\gamma_0$ and driven by suitable  neural networks: in the deterministic case $T_\theta$
is the time-$T$ flow map of the process, while in the stochastic case the
measure is the law of its terminal value, and samples are produced by
sampling $\gamma_0$ and simulating the process. The presence of an evolution introduces a time discretisation, and with it stability questions of the kind encountered in Section~\ref{SubS:Ex:parabolic}; these spaces are discussed separately below.


\subsubsection*{Wasserstein distance}

A popular choice of distance between probability measures is the Wasserstein distance, \cite{Villani_OT_book_2009}, which leads naturally to the formulation and computation of a wide range of quantities of interest within the framework of optimal transport. Let $p \geq 1 $ and $\mu, \nu$ probability measures on $ {\mathbb{R}^k}$   with finite $p$ moments.  
The Wasserstein distance of order $p$ between  {$\mu$ and $\nu$} is defined by
\begin{equation}\label {W-dist}
	W_p(\mu, \nu) = \Big \{ \min_{\pi \in A (\mu, \nu)} \int _{ {\mathbb{R}^k\times\mathbb{R}^k}} |x-y|^p\, \d \pi(x,y) \Big \} ^{1/p} ,
\end{equation}
where $A(\mu,\nu)$ denotes the set of \emph{transport plans} between $\mu$ and $\nu$, namely the Borel probability measures on $ {\mathbb{R}^k\times\mathbb{R}^k}$ whose first and second marginals are $\mu$ and $\nu$, respectively.
For $\mu $ a probability measure on $ {\mathbb{R}^k}$, and $ \xi  :  {\mathbb{R}^k} \to \R $ a  Borel measurable map, we denote by $\xi _\sharp \mu $  the push-forward of %
$\mu $ through $\xi$, i.e., the probability measure on $\R $  defined by 
$$ \xi _\sharp  \mu (B) :=  \mu(\xi ^{-1} (B)),\qquad \text{Borel } B\subset \R\, . $$ 
In particular, 
\begin{equation}\label {pf-int}
	 \int _{ {\mathbb{R}^k} } |\xi (x)|^p\, \d \mu(x)  =  \int  |s|^p\, \d \xi _\sharp \mu (s)  .
\end{equation}
With the obvious modification in the definition of $W_p$, 
we have, 
\begin{equation}\label {W-dist-pf}
	W_p(\xi _\sharp  \mu, \xi _\sharp \nu) = \Big \{ \min_{\pi \in A (\mu, \nu)} \int _{ {\mathbb{R}^k\times\mathbb{R}^k}} |\xi (x) -\xi (y) |^p\,  \d \pi(x,y) \Big \} ^{1/p} \, .
\end{equation}
Wasserstein distances between empirical measures can be computed efficiently, for instance by the Sinkhorn algorithm, which solves an entropically regularised, and hence approximate, optimal transport problem; see \cite{Peyre:CompOT:2019} for this and other computational approaches.


\subsubsection*{Generative Adversarial Networks}
Generative Adversarial Networks (GANs) provide a flexible variational framework for approximating probability distributions by minimising suitable statistical distances. These distances include Integral Probability Metrics (IPMs) and, under appropriate conditions, $f$-divergences. In the case of Integral Probability Metrics
we have for an appropriate  class of measurable test functions $\mathcal{F}$ { (denoted by $\varphi$, to distinguish them from the data $f$ of the previous model problems)},  
\begin{equation}\label{eq:IPM}
d_{\mathcal{F}}(\mu,\nu)
:=
\sup_{\varphi\in\mathcal{F}}
\left|
\int \varphi\,\d\mu
-
\int \varphi\,\d\nu
\right|.
\end{equation}
Many important distances between probability measures are of the form
\eqref{eq:IPM}. These include for example the total variation distance, the $1$-Wasserstein distance, the Maximum Mean Discrepancy among others. 
%
%

Considering such distances in
\eqref{mm_nn:abstract} we naturally arrive to a saddle point formulation.
%
Indeed, if
$$
\overline V_{\mathcal{N}}
=\Big \{ \varphi_\theta:\theta\in\overline{\Theta}\Big \}
$$
is a neural-network-based parameterised family of functions that provides a sufficiently rich approximation of the class $\mathcal{F}$, then the approximation of \eqref{mm_nn:abstract} takes the form
\begin{equation}\label{mm_nn:abstract:GAN}
	  \min_{ \mu  \in  V  _{\Theta } } \sup_{\varphi\in \overline V_{\mathcal{N}}}
\left|
\int \varphi \,  \d \mu
-
\int \varphi \,  \d \nu _{\, \text{data}}
\right|   
 \end{equation}
which forms the basis of adversarial methods.
{In Generative Adversarial Networks (GANs) \cite{goodfellow2014generative}, the neural measure space is static, $V_\Theta=\{(g_\theta)_\sharp\gamma_0:\theta\in\Theta\}$, generated by a family of \emph{generator} networks $g_\theta\in V_{\mathcal{N}}$; the original GANs correspond to the Jensen--Shannon divergence rather than to an IPM.}

%

{ 
An example of this framework is provided by Wasserstein GANs (WGANs), 
introduced in~\cite{arjovsky2017wassersteingan}, where the underlying IPM is given 
by the Wasserstein-1 distance $W_1$. By the 
Kantorovich--Rubinstein duality theorem~\cite{Villani_OT_book_2009}, the 
Wasserstein-1 distance admits the dual representation
\begin{equation*}
  W_1(\mu,\nu)
  = \sup_{\varphi\in \mathcal{F}_1 }
\int \varphi \,  \d \nu
-
\int \varphi \,  \d \mu
\end{equation*}
where the class of test functions is given by
\begin{equation*}
  \mathcal{F}_1 :=\Big \{\varphi: {K}\to \R:\,\,\varphi\text{ is Lipschitz with }\|\varphi\|_L\leq 1\Big \},
\end{equation*}
for a compact set $ {K\subset\mathbb{R}^k}$ containing the supports of the measures.
Hence, the Wasserstein distance can be interpreted as an IPM associated with the class $\mathcal{F}_1$ of Lipschitz test functions. In the WGAN framework we aim to optimise 
\begin{equation}\label{mm_nn:abstract:WGAN}
	  \min_{ \mu  \in  V  _{\Theta } } \sup_{\varphi\in \overline V_{\mathcal{N}}}
\
\int \varphi \,  \d \nu _{\, \text{data}} -\int \varphi \,  \d \mu
 \end{equation}
 where the discriminator class $\overline V_{\mathcal{N}}$ is selected as a suitable neural-network-based approximation of the class $\mathcal{F}_1$.

\subsubsection*{Flow based Generative Modelling}

Most modern generative models, such as normalising flows, score-based and diffusion models, involve \emph{flows} of measures that transform a simple reference measure into a complex one, or vice versa. In this case the neural measure space is flow-based: for each $\theta$ an evolution on $[0,T]$ generates a curve of measures $\mu_\theta(t)$ with $\mu_\theta(0)=\gamma_0$, and
\begin{equation*}
V_{\mathcal{N},\mathcal{P}} = \{ \mu_\theta(T) : \theta \in \Theta \}\, .
\end{equation*}
In the deterministic case $\mu_\theta(t)=(\Phi_\theta(t))_\sharp\gamma_0$, where $\Phi_\theta(t)$ is the flow map of an ordinary differential equation whose velocity field is a neural network; in the stochastic case $\mu_\theta(t)$ is the law at time $t$ of a stochastic differential equation (Markovian flow), and the flow of measures is described by the corresponding Fokker--Planck equation. The trained problem is again
\begin{equation}\label{abstract_min_genmod}
 	 \min  _ {\mu \in  V_{\mathcal{N},\mathcal{P}} } d_{\mathscr {M}} (\mu, \nu _{\, \text{data}} )\, .
\end{equation}
This interpretation is motivated by the approach in \cite{Katsoulakis_z2023meanfieldlab}, where the main mathematical models suggested are mean field games.

The space $V_{\mathcal{N},\mathcal{P}}$ must be sufficiently expressive, despite the fact that all its elements are obtained from the single reference measure $\gamma_0$. If successful, this approach \emph{links} $\gamma_0$ to
\begin{equation}\label{mm_nn:abstract:2:generative}
\mu_{\theta ^*} \cong \nu _{\, \text{true}}, \qquad \text {where } \mu_{\theta ^*}
\in \text{argmin} _ {\mu \in  V_{\mathcal{N},\mathcal{P}} } d_{\mathscr {M}} ({\mu ,\nu _{\, \text{data}} } ).
\end{equation}
Solving \eqref{abstract_min_genmod} determines both $\mu_{\theta^*}$ and the evolution that produces it, and the latter provides the generative mechanism: in the deterministic case, if $X\sim\gamma_0$ then $\Phi_{\theta^*}(T)(X)\sim\mu_{\theta^*}$, so that arbitrarily many \emph{synthetic data} are obtained by sampling $\gamma_0$ and transporting the samples through the learned flow. In practice the evolution is discretised in time and the distance is replaced by a numerical approximation, which yields a fully discrete problem; the evolution structure makes the time discretisation a central issue. Some important models arising from this framework are outlined below.

\vspace{2mm}

\noindent
{\bf Normalising flows.}
A deterministic flow is generated by the ordinary differential equation
\begin{equation}\label{eq:ODE-flow}
\frac{dX(t)}{dt}=v_\theta(t,X(t)),
\qquad
X_0\sim\gamma_0,
\end{equation}
where the velocity field $v_\theta$ is represented by a neural network, i.e., is an element of an appropriate $V_{\mathcal{N}} $.
The associated probability measures are given by
$$
\mu {_\theta} (t)=( {\Phi_\theta} (t))_\sharp\gamma_0,
$$
where $ {\Phi_\theta (t)}$ denotes the flow map generated by \eqref{eq:ODE-flow}. The
corresponding densities satisfy the continuity equation (for simplicity we use equations for densities; the equations of flows of measures are defined in the sense of distributions, see \cite{AGS:2005}[Chapter 8]) 
\begin{equation}\label{eq:continuity}
\partial_t\rho(t)+\nabla\cdot(\rho(t) v_\theta)=0.
\end{equation}
Flows of this type are known as continuous normalising flows; the classical normalising flows use instead compositions of invertible networks. They are trained by maximum likelihood, that is, with $d_{\mathscr {M}}$ the
Kullback--Leibler divergence: since
$\mathrm{KL}(\nu_{\, \text{true}}\,\|\,\mu_\theta(T))
=-\int\log\rho_\theta(T,x)\,\d\nu_{\, \text{true}}(x)+\text{const}$, with
$\rho_\theta(T,\cdot)$ the density of $\mu_\theta(T)$, replacing
$\nu_{\, \text{true}}$ by $\nu_{\, \text{data}}$ in the integral gives the
trained loss $-\frac1N\sum_i\log\rho_\theta(T,X_i)$.
The required density is available through the change of variables formula $\frac{\d}{\d t}\log\rho(t,X(t))=-\nabla\cdot v_\theta(t,X(t))$ along the trajectories of \eqref{eq:ODE-flow}.

\vspace{2mm}

\noindent
{\bf Diffusion models.}
Diffusion models employ stochastic dynamics. A \emph{fixed} forward stochastic
differential equation
\begin{equation}
dX(t)=b(t,X(t))\,dt+\sigma(t)\,dW(t),\qquad X(0)\sim\nu_{\, \text{true}},
\end{equation}
with no learned parameters, typically an Ornstein--Uhlenbeck process,
transports the data distribution over $[0,T]$ into a distribution close to the
Gaussian $\gamma_0$; the density $\rho$ of the law $\mu(t)$ is governed by the
Fokker--Planck equation
\begin{equation}
\partial_t\rho(t)
=
-\nabla\cdot(b\,\rho(t))
+
\frac12
\nabla\cdot
(\sigma\sigma^{\mathsf T}\nabla\rho(t)).
\end{equation}
The generative model is the time reversal of this process, which is again a
stochastic differential equation, started from $\gamma_0$, whose drift
$-b+\sigma\sigma^{\mathsf T}\nabla\log\rho$ involves the \emph{score}
$\nabla\log\rho(t,\cdot)$ of the forward process, the only unknown quantity.
The neural network enters through an approximation $s_\theta(t,x)$ of the
score: $\mu_\theta(T)$ is the law of the terminal value of the reverse equation
with $s_\theta$ in place of the score, started from a sample of $\gamma_0$,
and the neural measure space $\{\mu_\theta(T):\theta\in\Theta\}$ does not
depend on the data. The distance is again the Kullback--Leibler divergence. In many applications,
however, the density of $\mu_\theta(T)$ is computationally intractable, and
the loss cannot be evaluated as for normalising flows; training then
minimises instead a computable upper bound of
$\mathrm{KL}(\nu_{\, \text{true}}\,\|\,\mu_\theta(T))$, the \emph{score
matching} objective, which measures the distance between $s_\theta$ and the
score along the forward process and is expressed through the explicit
Gaussian transition law of the forward process and the data. 
 
\vspace{2mm}

\noindent
{\bf Wasserstein flows.}
Flows of probability measures play a central role in the geometry of Wasserstein space and are closely related to the generative models discussed above. While the latter typically construct probability measures by transporting a prescribed reference distribution $\gamma_0$ through a learned map, Wasserstein flows formulate the evolution directly on the space of probability measures. The strong connection between these two viewpoints motivates the following brief discussion.

The geometry of the Wasserstein space provides a natural framework for the analysis of flows of probability measures,
\cite{AGS:2005}, treating probability measures as points of a metric space
$$
(\mathcal P_2(\R ^ k),W_2),
$$
equipped with the $2$-Wasserstein distance. Instead of approximating a single
probability measure, one seeks an evolution
$$
(\mu(t))_{t\in[0,T]}
$$
connecting two probability distributions while remaining optimal with respect
to an underlying variational principle.
The basic kinematic constraint is the continuity equation for the densities $\rho (x, t) \geq 0$
\begin{equation}\label{eq:continuity2}
\partial_t\rho +
\nabla\cdot(\rho  v )
=
0,
\end{equation}
where $v=v(x,t) \in \R ^k$ denotes the velocity field transporting the probability mass.

The Benamou--Brenier formulation, \cite{Benamou-Brenier:2000},
\cite{Villani_OT_book_2009}, shows that the squared Wasserstein distance is,
up to the factor $T$, the minimal kinetic energy among all admissible flows
on $[0,T]$:
\begin{equation}\label{eq:BB}
\begin{split}
W_2^2(\mu_0,\mu_T)
=
\inf_{\rho ,v}
\Big \{ T\,
& \int_0^T
\int_{\R ^k}
|v (x,t)|^2
\,\rho (x, t) \,\d x \,\d t
\ :\ \eqref{eq:continuity2}
\text{ holds},\\
&
\mu(0)=\mu_0,
\;
\mu(T)=\mu_T
\Big \} .
\end{split}
\end{equation}
Thus the Wasserstein distance is obtained by solving an optimisation problem
over flows of probability measures rather than over transport maps. In the
context of flow-based neural measure spaces this formula has a natural
interpretation. The kinetic energy
$\int_0^T\int_{\R^k}|v_\theta|^2\,\d\mu_\theta(t)\,\d t$ of a flow generated
by $v_\theta$ is, by \eqref{eq:BB}, bounded from below by
$W_2^2(\gamma_0,\mu_\theta(T))/T$, with equality exactly for the flows that
realise the optimal transport from $\gamma_0$ to $\mu_\theta(T)$. Adding it
to the loss as a regularising term,
\begin{equation}\label{abstract_min_genmod_BB}
 	 \min_{\theta\in\Theta}\Big\{ d_{\mathscr {M}}(\mu_\theta(T),\nu_{\, \text{data}})
 	+\lambda \int_0^T
\int_{\R ^k}
|v_\theta (x,t)|^2
\,\d \mu_\theta(t)(x)\,\d t \Big\}, \qquad \lambda>0,
\end{equation}
therefore favours, among the flows with a good terminal fit, those which are
close to optimal transport from $\gamma_0$, and generators that are
$W_2$-close to the reference measure. Note that the terminal loss alone does
not determine the flow: many velocity fields transport $\gamma_0$ to the same
terminal measure, so that the minimisation over flows has many minimisers,
and the kinetic energy term favours among them those of minimal energy. 
Flows obtained in this way are closely related to Wasserstein
gradient flows, in which the evolution is itself determined by a variational
principle on $(\mathcal{P}_2(\R^k),W_2)$ \cite{AGS:2005}. 
\vspace{2mm}
\noindent
{\bf Wasserstein gradient flows of the loss.}
In the models above the evolution that generates the neural measure space is prescribed in advance. Alternatively, it can be generated by an energy $F$ on $\mathcal{P}_2(\R^k)$
through the \emph{minimizing movement} scheme, introduced by De Giorgi for
gradient flows in general metric spaces and first applied to the space of
probability measures with the Wasserstein metric by Jordan, Kinderlehrer and
Otto \cite{JKO}, see \cite{AGS:2005}:
\begin{equation}\label{eq:JKO}
\mu^{n+1}\in\operatorname*{argmin}_{\mu\in\mathcal{P}_2(\R^k)}\Big\{F(\mu)+\frac{1}{2\tau}W_2^2(\mu,\mu^n)\Big\},
\end{equation}
an analogue of the implicit Euler discretisation, with time step $\tau$, of the Wasserstein
gradient flow of $F$; the $W_2^2$ term plays the role of the proximal term in Section~\ref{SubS:GD}. By the Benamou--Brenier formula, $\frac{1}{2\tau}W_2^2(\mu,\mu^n)$ is the minimal kinetic energy $\frac12\int_0^\tau\!\int_{\R^k}|v|^2\rho\,\d x\,\d t$ among flows from $\mu^n$ to $\mu$ on $[0,\tau]$, so that \eqref{eq:JKO} can be written over flows with a free end, and the velocity field becomes the neural object:
\begin{equation}\begin{split}
\label{eq:JKO-nn}
\theta^n\in\operatorname*{argmin}_{\theta\in\Theta,\ \partial_t\rho_\theta+\nabla\cdot(\rho_\theta v_\theta)=0,\  \rho_\theta(0)=\mu^n }\Big\{ &F\big(\rho_\theta(\tau)\big)+\frac12\int_0^\tau\!\!\int_{\R^k}|v_\theta|^2\rho_\theta\,\d x\,\d t\Big\},
\end{split}\end{equation}
and $\mu^{n+1}=\rho_{\theta^n}(\tau)$. A scheme of this form, with an additional Fisher information regularisation, is considered in \cite{Li_Lu_Wang:2020}. The resulting neural measure space is flow-based in discrete time, and the levels of approximation are the time step $\tau$, the restriction of each step to the velocity fields $v_\theta$, and the samples through which every quantity is evaluated. Several choices of $F$ and of the starting point are of interest. 
For example, if $F=d_{\mathscr {M}}(\cdot,\nu_{\, \text{data}})$ and $\mu^0=\gamma_0$ the flow moves from the reference measure towards   $\nu_{\, \text{data}},$ i.e.,  can be seen as the Wasserstein gradient flow analog of \eqref{mm_nn:abstract}.
For $F=\mathrm{KL}(\cdot\,\|\,\gamma_0)$ and $\mu^0=\nu_{\, \text{data}}$ the gradient flow is the 
Fokker--Planck equation, which is the original motivation of \eqref{eq:JKO}; the scheme is then a time-discrete, variational construction of the forward process of a diffusion model, whose reversal is the generative model, and score-based models can indeed be interpreted as Wasserstein proximal operators of this type \cite{zhang2024wassersteinproximaloperatorsscorebased}. Finally, the gradient flow of $F=d_{\mathscr {M}}(\cdot,\nu_{\, \text{data}})$ converges, as $t\to\infty$, to $\nu_{\, \text{data}}$ itself, so that memorisation appears here as the long-time limit.

In summary, the models above are all instances of the trained problem \eqref{mm_nn:abstract}. 
The analytical framework developed in these notes can, in principle, be
applied to all of these models. As an illustration, we consider in
Sections~\ref{S:Analysis} and~\ref{S:reliability} the simplest case, a
static neural measure space with $d_{\mathscr {M}}=W_p$. Flow-based models
are considerably more involved
and are left for future work.

Beyond the classification, the framework makes visible the successive
levels of approximation involved in a flow-based generative model: the
continuous and the semidiscrete problems, the replacement of
$\nu_{\, \text{true}}$ by the data, the time discretisation of the evolution,
and the approximation of the distance, by a discriminator, by entropic
regularisation or by score matching, together with the sampling of
$\mu_\theta.$   Separating these levels is the counterpart, for generative models,
of the hierarchy discussed in Section~\ref{SubS:GD}, and it allows the
questions of consistency and stability to be posed at each level separately;
this is the perspective adopted in these notes.

 \subsection{Bibliography and comments}\label{SubS:Ex:bib}

The notation and terminology for the function learning problem is from
\cite{loulakis2023newap}. The probabilistic function learning problem is a
central problem in Statistical Learning and ML, see the review
\cite{Bartlett_Montanari_Rakhlin_actaN:2021} and its references. Other works
include \cite{Bartlett_VCdimension:2019, Oneto_local_Rademacher:2018,
Yousefi_local_Rademacher:2018, Bartlett_Benignoverfitting:2020,
Dimakis_ExactlyComputingtheLocalLipschitz:2020, Montanari_neural_kernel:2021},
\cite{Jentzen_full_err:2022}, \cite{Bartlett:gen_2023}. For operator learning, that is, the learning of maps between function
spaces, see e.g., \cite{lu2021learning, Mishra_FNO_est,
Mishra_DeepONET_Kran_est, Stuart_NeuralOper2023,
Mishra_Conv_operNEURIPS2023}.

Numerical solution of PDEs utilising neural networks has received a lot of
attention. Deep Ritz methods \cite{e2017deep} and Physics Informed Neural
Networks (PINNs) \cite{Karniadakis:pinn:orig:2019} have become very
influential. Residual based methods were considered in \cite{Lagaris_1998},
\cite{Berg_2018}, \cite{SSpiliopoulos:2018} and their references. Other
neural network methods for differential equations and related problems
include, for example, \cite{kevr_1992discrete}, \cite{kharazmi2019variational},
\cite{Xu}, \cite{chen2022deep}, \cite{georgoulis2023discrete},
\cite{Grohs:space_time:2023},
\cite{lubich2025regularizeddynamicalparametricapproximation},
\cite{FLLubichN2026}. For a broader perspective of the related
methodologies, see e.g., \cite{karniadakis_kevr_2021physics}. The
formulation of time-discrete training method for parabolic problems follows
\cite{GGM_IMA_2025}.

Representative works considering various ML approaches on generating
probability distributions include \cite{adler2018banach,
goodfellow2014generative, arjovsky2017wassersteingan,
zhang2018wassersteinganmodeltotal, creswell2018generative, jospin2022hands,
Chakraborty_Bartlett:2025, chen2026highaccuracysamplingdiffusionmodels}.
The formulation of Section~\ref{SubS:Ex:measures}, the minimisation of a
distance to the empirical measure over a parametrised family, is in
statistical terms a minimum distance estimation problem; for GANs,
large-sample properties of the resulting estimators are studied in
\cite{Biau_Cadre_Sangnier_Tanielian:2020}. Generating measures by
transporting a reference measure through parametrised maps, the viewpoint
underlying the neural measure spaces, goes back   to
\cite{Tabak_VandenEijnden:2010}, and a mathematical introduction to
generative models from this viewpoint, with emphasis on its relation to
optimal transport, is given in \cite{Ruthotto_Haber:2021}. Continuous
normalising flows originate in \cite{Chen_NeuralODE:2018}, and the
formulation of diffusion models through a forward stochastic differential
equation and its time reversal in \cite{Song_SDE:2021}. For Wasserstein
gradient flows and the minimizing movement scheme see \cite{JKO, AGS:2005};
the dynamic form of the scheme, in which each step is written through the
Benamou--Brenier formula, is used computationally in
\cite{Benamou_Carlier_Laborde:2016, Carrillo_Craig_Wang_Wei:2022,
Li_Lu_Wang:2020}. A broad survey of the role of flows of probability measures in machine
learning is given in \cite{Peyre_transports_survey:2025},
where the models are unified through the Eulerian and Lagrangian
descriptions of an evolution of measures and the non-uniqueness of the
velocity fields that generate it is identified as both an obstacle and a
design opportunity. A unified treatment of flow-based and diffusion-based models, different in
perspective from the one adopted here, in which both arise from a single
stochastic process interpolating between the reference and the data
distribution, is given in
\cite{Albergo_VandenEijnden:2023, Albergo_Boffi_VandenEijnden:2025}.  
The work \cite{Katsoulakis_z2023meanfieldlab} provides a
modelling perspective that interprets various generative models as mean
field games, see also \cite{zhang2024wassersteinproximaloperatorsscorebased}. Spaces based on neural measures were used to approximate randon differential equations in \cite{arampatzis2025neuralmeasureslearningdistributions}.
Unifying viewpoints on generative models exist, in particular through
push-forwards of a reference measure; the organisation by levels of
approximation adopted here appears, to our knowledge, to be new.

 \newpage
 \newpage
 

\section{Stability, a priori, a posteriori estimates in numerical analysis}\label{S:apriori_apost}

\subsection{Error analysis}\label{SubS:aPaP:Error_analysis}

\subsubsection*{Stability of the PDE vs stability of the numerical scheme}
We briefly outline the main principles underlying the derivation of error estimates for numerical methods for differential equations and PDEs. Although the general strategy is relatively transparent, its rigorous implementation can be technically involved and far from straightforward.
Consider the abstract problem $$\mathscr{A} u =0,$$ $u  \in \mathcal{V } $
and its discretisation
$$\mathscr{A}_h\, u _h=0, \quad u _h \in \mathcal{V }_h\, .$$

Assume for simplicity that both are linear. Then
\emph{PDE stability} can be roughly expressed as
$$\mathscr{A}v =g \implies \| v\|_{N_1}\leq C _{pde}\, \| g\| _{N_2}\, ,$$
where $ \| \cdot \|_{N_1}, \, \| \cdot \| _{N_2}\, ,$ are appropriate norms not specified.

On the other hand,
\emph{discrete stability} or stability of the numerical scheme, typically can be  expressed as
$$\mathscr{A}_h\,v_h =g_h \implies \| v_h\|_{N_{h,1}}\leq C _{discr}\, \| g_h\| _{N_{h,2}} \leq C\, ,$$
where $ \| \cdot \|_{N_{h,1}}, \, \| \cdot \| _{N_{h,2}}\, ,$ are appropriate discrete norms, which may be related to
$ \| \cdot \|_{N_1}, \, \| \cdot \| _{N_2}\, .$

\subsubsection*{A priori estimates}
A priori estimates are based on the  stability of the numerical scheme and yield convergence, provided the exact solution is smooth enough. To achieve this one has to define first an appropriate projection operator $ P _h:  \mathcal{V } \to \mathcal{V }_h
$ and to evaluate
$$\mathscr{A}_h\, P_h u  \, ,$$
which is called the \emph{consistency error}, and quantifies the extent to which $u,$ or more precisely $P_h u,$ fails to satisfy the equation defining the numerical scheme. \emph{Consistency error} may be of high-order of accuracy, depending on the design of the scheme. \emph{Consistency error} and \emph{discrete stability} yield \emph{a priori} error bounds:
$$\mathscr{A}_h\,(P_h u -u_h)  =(\mathscr{A}_h\,P_h u  ) \implies \| (P_h u -u_h)\|_{N_{h,1}}\leq C _{discr}\, \| \mathscr{A}_h\,P_h u\| _{N_{h,2}}\, .$$

The principle that consistency and discrete stability imply convergence is the content of the Lax equivalence theorem \cite{Lax_Ricth_1956}. The paper did more than prove an implication: it introduced the very notion of stability of a numerical scheme, as the uniform boundedness of the family of discrete solution operators with respect to the discretisation parameters, and showed that for well-posed linear problems this notion is exactly what separates consistent schemes that converge from those that do not. Earlier, stability had been a property to be tested case by case, for instance by von Neumann's Fourier analysis; after \cite{Lax_Ricth_1956} it became a structural concept, with its own meaning independently of any particular scheme. The concept is now so standard that it is often applied mechanically, as a condition to be verified, and its depth tends to be forgotten: discrete stability is a statement about the discretised problem alone, which is a requirement in any scheme design process. It is precisely in this sense that we shall use it. In Section~\ref{S:Analysis} the principle reappears in variational form.

\subsubsection*{A posteriori estimates}
A posteriori estimates are based on the  stability of the PDE and yield error control which depends solely on the approximate computed solution $u_h\, .$  This time one has to   define first an appropriate reconstruction operator  $ R:  \mathcal{V } _h \to \mathcal{V }\, .
$ Then, the quantity
$$\mathscr{A}  (R u _h )   \, ,$$
which is often called  \emph{residual}, and quantifies the extent to which $R u _h$ fails to be the solution of the PDE, controls the error. \  \emph{Residual} and \emph{PDE stability} yield \emph{a posteriori} error bounds:
$$\mathscr{A} (R u_h  -u)  =(\mathscr{A} (R u_h ) ) \implies \| R u_h  -u \|_{N_{1}}\leq C _{PDE}\, \| \mathscr{A} (R u_h )  \| _{N_{2}}\,
 .$$

\subsection{A linear parabolic model problem}\label{SubS:aPaP:parabolic}

The notion of stability in numerical analysis takes several forms and depends on the numerical method and on the continuous problem being discretised. For the sequel, it
will be important to discuss a very elementary example in order to reveal
the hidden structures of the discretisation that govern stability. To this
end, consider the
\emph{Evolution Problem: }\ Find $u: [0,T] \to D(A)$ satisfying
	\begin{equation*}
\label{1}
\begin{aligned}
&u' (t) + Au(t)= f(t)  , \quad 0<t<T,\\
& u(0)=u^0,
\end{aligned}
\end{equation*}
with $A$ a positive definite, self-adjoint, linear operator
on a Hilbert space $(H, \< \cdot , \cdot \> )$ with domain $D(A)$
dense in $H.$
We discretise in time this equation by considering a partition $0=t^0<t^1<\cdots <t^N=T$ of $[0,T]$ and, as in Section~\ref{SubS:Ex:parabolic}, $I_n:=(t^{n-1},t^n],$ $k_n:=t^n-t^{n-1},$ $n=1,\dots,N.$ A form of the \emph{Backward Euler} discretisation is:
\begin{equation*}
\label{BE}
  \frac 1{k_n} (U^{n}-U^{n-1}) +  A\, U^{n}  = \bar f ^{n}\, , \qquad n=1,\dots,N,
\end{equation*}
where
$
  \bar f^n = \frac1{k_n}\int _{I_n} f(s)\, ds\, .
$
It will be useful to consider the \emph{forward Euler} discretisation as well:
 \begin{equation*}
\label{FE}
  \frac 1{k_n} (U^{n}-U^{n-1}) +  A\, U^{n-1}  =   f (t^{n-1}) \, , \qquad n=1,\dots,N.
\end{equation*}
It is well known that these two schemes have fundamentally different properties and computational behaviour.

\subsubsection*{Reconstructions}
Our aim is to rewrite the scheme as a perturbation of the evolution equation. This will be useful for both the stability and \emph{a posteriori} analysis.  We shall introduce the \emph{time reconstruction}
   $\widehat U$ 
 of a time dependent function  $U$   to be the piecewise linear approximation
of $U$ defined by linearly
interpolating between the nodal values $U^{n-1}$ and $U^{n}$:
\begin{equation}\label{def_hatU}
\widehat U(t):=\ell_0^n(t) U^{n-1}+ \ell_1^n(t) U^{n}, \quad t\in I_n,
\end{equation}
here
$$
\ell_0^n(t) := (t^{n}-t)/k_n\quad \text{and } \ell_1^n(t) := (t-t^{n-1})/k_n.$$
Correspondingly, the piecewise constant interpolant of $U^{j}$  is denoted by $\overline U, $
\begin{equation*}
\label{def_overU}
\overline U(t):= U^{n}, \quad t\in I_n\, .
\end{equation*}
Later we shall need the \emph{lower} piecewise constant interpolant of $U^{j}$
\begin{equation*}
\label{def_overU-d2}
\underline U(t):= U^{n-1}, \quad t\in I_n\, .
\end{equation*}
Then  we \emph{rewrite} the Backward Euler scheme:  \   $\widehat    U, \overline U : [0,T] \to D(A)$
	\begin{equation*}
 \widehat    U  ' (t) + A\,  \overline  U(t) =       \overline
f (t)\, , \qquad  {\overline f(t):=\bar f^{\,n}, \quad t\in I_n\, .}
\end{equation*}
Correspondingly, the \emph{forward Euler} scheme is
 \begin{equation*}
 \widehat    U  ' (t) + A\,  \underline  U(t) =   \underline f (t )
 \, , \qquad  {\underline f(t):=f(t^{n-1}), \quad t\in I_n\, .}
\end{equation*}

\subsubsection{Stability and \emph{a priori} analysis}
To avoid some technicalities we assume in this section that $f=0\, .$
  Then
the   \emph{Backward Euler} method   seeks    $\widehat    U, \overline U : [0,T] \to D(A)$  such that
	\begin{equation*}
 \widehat    U  ' (t) + A\,  \overline  U(t) =        0\, .
\end{equation*}
A first important remark is that by solving
$
  \frac 1{k_n} (U^{n}-U^{n-1}) +  A\, U^{n}  = 0\, ,
$
 we obtain
 \begin{equation*}
\label{BE_step}
   \left [ I + k _n A\right ]U^{n} = U^{n-1} \, .
\end{equation*}
   Thus, at every time step,
$U^{n}$
  is the solution of a stationary (elliptic) problem. This explains why the backward scheme exhibits a rapid increase in regularity as time evolves:
$$\text{ if $U^0 \in H$ then $U^1 \in D(A),$ ...., $U^j \in D (A^j).$}$$
This fundamental property is the discrete analogue of the smoothing effect inherent in the continuous evolution
problem.
%

\subsubsection*{Stability}
Next, we
  multiply by $ \overline  U(t)$ to obtain,
	\begin{equation*}
\int_0 ^T  \< \widehat   U  ' (t) ,  \overline  U(t) \> + \int_0 ^T  \< A\,  \overline  U(t) ,  \overline  U(t)\> =        0\, .
\end{equation*}
 then
\begin{equation} \label{dMR}
	\begin{split}
 &
  \sum _{n=1}^N \, \int_{I_n}  \,
  \<   \widehat U ' ,   \overline U     \> \, \ d t
 =
  \sum _{n=1}^N \, \int_{I_n}  \,   \<    \big [ \frac{U ^n-  U^{n-1}}{k_n}  \big ],  U ^n          \> \, \ d t \\
& =  \frac 1 2  \|      {U ^N }   \| _H^2          + { \frac 1 2 \sum _{n=1}^N \,    \,    \|     {U ^n-  U^{n-1}}  \|   _H^2  }      \,  -  \frac 1 2 \| {U ^0 }   \|^2 _H   \, .
 \end{split}
\end{equation}
Hence
the discrete stability follows:
\begin{equation*} \label{dMR-d2}
	\begin{split}
   \frac 1 2  \|      {U ^N }   \| _H^2          &+  \frac 1 2 \sum _{n=1}^N \,    \,    \|     {U ^n-  U^{n-1}}  \|   _H^2       \\
             &+   \sum _{n=1}^N \,    \, k _n   \| A^{1/2}    {U ^n }  \|   _H^2        \,  \leq  \frac 1 2 \| {U ^0 }   \|^2 _H   \, .
 \end{split}
\end{equation*}
This expression can be written differently by using the globally defined functions $ \widehat   U$ and $ \overline  U $:

	\begin{equation*}
\frac 1 2  \|   \widehat   U  (T)   \| _H^2      + \frac 1 2 \sum _{n=1}^N \int_{I_n} k_n \| \widehat   U  ' (t)  \| ^2 _H  + \int_0 ^T  \|  A ^{1/2}\,  \overline  U(t) \| _H^2\leq  \frac 1 2 \| {U ^0 }   \|^2 _H\, .
\end{equation*}

\subsubsection*{Error estimate}
Let $\widehat u$ and $\overline u$ denote the piecewise linear and piecewise constant interpolants of the exact solution at the nodes, defined as $\widehat U$ and $\overline U$ with $U^n$ replaced by $u(t^n)$; in the terminology of Section~\ref{SubS:aPaP:Error_analysis} they play the role of $P_h u$. Let $\widehat e= \widehat u - \wU, $  $\overline e= \overline u - \overline U.$ Since $\widehat U'+A\overline U=0$, the error satisfies $\widehat e\,'+A\overline e=\widehat u\,'+A\overline u$, and the stability estimate above, applied to $\widehat e,\overline e$ with this right-hand side, yields
\begin{equation*}\begin{split}
		  \|   \widehat  e  (T)   \| _H^2      +    \sum _{n=1}^N \int_{I_n} k_n \| \widehat   e  ' (t)  \| ^2 _H  + &   \int_0 ^T  \|  A ^{1/2}\,  \overline  e(t) \| _H^2\leq   \| {U ^0 - u^0}   \|^2 _H\, \\
		 + &  \int_0 ^T  \|  A ^{-1/2} \Big [ \, \widehat u ' + A\,  \overline    u \Big ]\, \| _H^2 .\\
\end{split}
\end{equation*}
 This is an \emph{a priori} error estimate, which can be used to show convergence and corresponding rates: the right-hand side involves only the data and the consistency error
 $ \rho_{k}  (u) :=  \widehat u ' + A\,  \overline    u , $
and $\rho_{k}(u) \to 0$ as the maximal step $ k:=\max_n k_n\to 0,$ provided $u$ is smooth enough. 

  \subsubsection{Forward Euler for $f=0$  }\  Using the reconstructions $\widehat    U, \underline  U : [0,T] \to D(A)$  the \emph{forward Euler} scheme takes the form
	\begin{equation*}
 \widehat    U  ' (t) + A\,  \underline  U(t) =        0\, .
\end{equation*}
By solving $  \frac 1{k_n} (U^{n}-U^{n-1}) +  A\, U^{n-1}  = 0\, $
we observe that this method has a \emph {de-smoothing effect,} in fact,
 \begin{equation*}
\label{FE_step}
   U^{n} = \left [ I - k _n A\right ] U^{n-1} \, .
\end{equation*}
Thus, it is required that $U^{n-1} \in D(A)$ in order to define $U^{n} \in H$. Hence, for $U^N$ to be well defined, one needs to assume that $U^0 \in D(A^N)$. An alternative calculation, which reveals the well-known instability of the forward Euler scheme, is given by the estimate obtained by multiplication by $\underline{U}(t):$
	\begin{equation*}
\int_0 ^T  \< \widehat   U  ' (t) ,  \underline  U(t) \> + \int_0 ^T  \< A\,  \underline  U(t) ,  \underline  U(t)\> =        0\, .
\end{equation*}
 then
\begin{equation*} \label{dMR-d3}
	\begin{split}
 &
  \sum _{n=1}^N \, \int_{I_n}  \,
  \<   \widehat U ' ,   \underline U     \> \, \ d t
 =
  \sum _{n=1}^N \, \int_{I_n}  \,   \<    \big [ \frac{U ^n-  U^{n-1}}{k_n}  \big ],  U ^{n-1}          \> \, \ d t \\
& =  \frac 1 2  \|      {U ^N }   \| _H^2         \\
&  \quad     { - \frac 1 2 \sum _{n=1}^N \,    \,    \|     {U ^n-  U^{n-1}}  \|   _H^2  }      \,  -  \frac 1 2 \| {U ^0 }   \|^2 _H   \, ,
 \end{split}
\end{equation*}
and, since $U^n-U^{n-1}=-k_nAU^{n-1}$, this leads to the  \emph{conditional} discrete stability
\begin{equation*} \label{dMR-d4}
	\begin{split}
   \frac 1 2  \|      {U ^N }   \| _H^2
             &+   \sum _{n=1}^{N} \,    { \, \left \{  k _n   \| A^{1/2}    {U ^{n-1} }    \|   _H^2  - \frac12\, k_n^2\, \|  A   {U ^{n-1} }  \|   _H^2  \right \}}     \,  \leq  \frac 1 2 \| {U ^0 }   \|^2 _H   \, .
 \end{split}
\end{equation*}
The term in braces is nonnegative on the eigenmodes of $A$ with eigenvalue $\lambda\le 2/k_n$ and negative on the others; for an unbounded operator $A$ the scheme is therefore unstable, and for a spatial discretisation $A_h$ of $A$, whose largest eigenvalue is of order $h^{-2}$, one obtains the very strict mesh condition $k_n \leq C h^2$, which is required for the stability of fully discrete forward Euler schemes for parabolic problems. This is the same mechanism as the condition $\eta_n<2/\lambda_{\max}$ for gradient descent in Section~\ref{SubS:GD}, gradient descent being the forward Euler method for the gradient flow of the loss.

\subsubsection{A posteriori estimates}
  Such estimates are very important for several reasons. In our setting, simple stability of the evolution problem yields \emph{a posteriori} error control, i.e., estimates in which the upper bound depends solely on the computed approximate solution and the data of the problem. To this end, notice that the \emph{Backward Euler} scheme for $f=0$ consists in finding $\widehat{U}, \overline{U} : [0,T] \to D(A)$ such that
 $
 \widehat    U  ' (t) + A\,  \overline  U(t) =        0\, .
$
or
 \begin{equation*}
 \widehat    U  ' (t) + A\,  \widehat   U(t) =  A\, [   \widehat   U(t)- \overline  U(t)]  = : R (\wU,  \overline  U)   \, .
\end{equation*}
Then the error estimate for $u - \wU$ is immediate:
\begin{equation*}\begin{split}
			 \|   u - \wU  \| _H^2       + &   \int_0 ^T  \|  A ^{1/2}\, \big (u - \wU \big ) \| _H^2\leq   \| {U ^0 - u^0}   \|^2 _H\, \\
		& \quad  +    \int_0 ^T  \|  A^{1/2} \, [  \widehat   U(t)- \overline  U(t) ] \, \| _H^2 .\\
\end{split}
\end{equation*}
It is very interesting that a similar bound holds for the
  \emph{Forward Euler} scheme as well. The scheme is now: Seek   $\widehat    U, \underline  U : [0,T] \to D(A)$
 such that $
 \widehat    U  ' (t) + A\,  \underline U(t) =        0\, .
 $
 or
 \begin{equation*}
 \widehat    U  ' (t) + A\,  \widehat   U(t) =  A\, [   \widehat   U(t) -\underline  U(t) ]  = : R (\wU,  \underline  U)
\end{equation*}
which leads to the
 error estimate for $  u - \wU : $
\begin{equation*}\begin{split}
			 \|   u - \wU  \| _H^2       + &   \int_0 ^T  \|  A ^{1/2}\, \big (u - \wU \big ) \| _H^2\leq   \| {U ^0 - u^0}   \|^2 _H\, \\
		& \quad  +    \int_0 ^T  \|  A^{1/2} \, [ \widehat   U(t) -\underline  U(t) ] \, \| _H^2 .\\
\end{split}
\end{equation*}
  Recall that
  $\wU(t) = \ell ^n _0 (t) U^{n-1} + \ell ^n _1 (t) U^{n} ,$
 and thus     $\wU(t) - \overline U(t) = \ell ^n _0 (t) (U^{n-1} - U^{n} ),$
   while $\wU(t) - \underline U(t) = \ell ^n _1 (t) ( U^{n} -U^{n-1}  ),$ for $t\in I_n$. Since $\int_{I_n}|\ell_0^n|^2\,dt=\int_{I_n}|\ell_1^n|^2\,dt=k_n/3$,
   this leads to the \emph{same a posteriori} error estimate for $u - \wU$:
 \begin{equation*}\begin{split}
	\max_{t\in [0, T]}	  \|   u - \wU  \| _H^2       + &   \int_0 ^T  \|  A ^{1/2}\, \big (u - \wU \big ) \| _H^2\leq   \| {U ^0 - u^0}   \|^2 _H\, \\
		& \quad  +    \sum _{n=1}^{N} \,  \frac{k _n}{3} \|  A^{1/2} \, [ U^{n} -U^{n-1}  ] \, \| _H^2 .\\
\end{split}
\end{equation*}
This estimate holds for \emph{both backward and forward Euler} schemes.
We note finally that the two schemes are precisely the time discretisations underlying the time-discrete training functionals $\mathcal{G}_{k,IE}$ and $\mathcal{G}_{k,EE}$ of Section~\ref{SubS:Ex:parabolic}. As we shall see in Section~\ref{SubS:PINN-Parabolic}, the difference in stability exhibited above is inherited by the corresponding losses: only the implicit one is coercive.

\subsection{Relevance for ML algorithms}

Error estimates, whether of \emph{a priori} or \emph{a posteriori} type, provide fundamental tools for assessing the quality and performance of numerical algorithms. \emph{A priori} estimates rely on the stability of the discrete scheme, with bounds depending only on properties of the exact solution. It is important to emphasise that discrete stability is relevant even in the pre-asymptotic regime. Indeed, it is well known that when restrictions on the mesh size or time step, such as those imposed by stability or CFL conditions, are violated, the resulting methods can exhibit poor numerical performance.

\emph{A posteriori} estimates, on the other hand, can be particularly valuable even when a complete convergence theory is unavailable. In the context of machine learning, upper bounds for the generalisation error may be computable from quantities available during the training process. If these bounds remain consistently small, they provide a posteriori evidence that the resulting approximation is accurate, without requiring a complete \emph{a priori} convergence theory for the algorithm.


A central approach in both statistical learning and ML for PDEs is to derive estimates of hybrid type, combining \emph{a posteriori} computable quantities with an \emph{a priori} assumption — typically that the discrete minimisers have bounded derivatives, or that all members of the neural-network space satisfy a similar uniform bound. 
Such estimates can be useful, and are sometimes the only available means of bounding the generalisation error. However, the analysis of Section~\ref{SubS:aPaP:parabolic} reveals that this approach deserves careful scrutiny. Since the boundedness assumption effectively replaces discrete stability as a hypothesis, hybrid estimates do not, in general, yield information about whether the underlying algorithm is stable. {An alternative, which avoids uniform bounds over the network class, is developed in Section~\ref{S:reliability}.}

There is, moreover, a less familiar point revealed by the same analysis. \emph{A posteriori} estimates can remain valid even for unstable schemes. For example, the estimate for the Forward Euler method remains valid regardless of whether the scheme is stable. Moreover, if the corresponding discrete time derivative is uniformly bounded, {in the sense that $\sum_n k_n\|A^{1/2}(U^n-U^{n-1})/k_n\|_H^2\le C$, the right-hand side of the estimate is of order $k^2$, and} the estimate itself implies convergence. In view of the above, this observation may have important implications for the analysis of ML algorithms, since imposing such bounds a priori on the neural-network space is itself a way of enforcing stability, rather than deriving it from the scheme.


\subsection{Bibliography and comments}\label{SubS:aPaP:bib}

A modern accessible account with interesting historical references on stability, consistency, and convergence of numerical methods can be found in \cite{Arnold2015}.
For the error analysis of time discretisations of parabolic problems in this spirit see \cite{Thomee_book_2006}. 
A posteriori error control in finite elements was pioneered by I. Babu\v{s}ka, see e.g.,  \cite{
Babuska_Rheinboldt_SINUM:1978}, and is now a central tool for error control
and adaptivity in many areas of numerical PDEs. The first proofs of optimal convergence rates for adaptive methods were
obtained in the wavelet setting by Cohen, Dahmen and DeVore
\cite{Cohen_Dahmen_DeVore:2001}. For the theory of adaptive finite element methods we refer to  \cite{Bonito_Canuto_Nochetto_Veeser:2024} and its references.

The time reconstruction $\widehat U$ of the discrete solution has been proven useful in various instances, see \cite{baiocchi_brezzi_par_1983}, \cite{nochettoSV2000posteriori}, \cite{Dem_Stuart_Tzavaras_2001}, and \cite{MakrN2006posteriori} for higher-order versions; the corresponding reconstruction in space, the elliptic reconstruction, is introduced in \cite{MN:siam2003}, see also \cite{LM:mc2006}, and \cite{AMN2, AMN3} for the extension to higher-order time discretisations; for a posteriori error estimation for elliptic problems  see \cite{verfuerthBook}.

The derivation of a posteriori estimates for the implicit Euler scheme follows the ideas of the seminal work \cite{nochettoSV2000posteriori}, where nonlinear evolution problems, including gradient flows, are considered. The resulting bounds provide a fundamental ingredient in the analysis of minimizing movement schemes within the framework developed in \cite{AGS:2005}.
 \newpage

\section{Analysis: Stability and Convergence}\label{S:Analysis}
 
 In this section, we focus on the \emph{stability and convergence} of machine learning algorithms. We begin with a somewhat abstract discussion, introducing the main concepts in an informal manner without aiming for complete mathematical precision. We then illustrate how this analytical framework can be applied to several of the examples presented in the previous section. 

 \subsection{Analysis at an abstract level}\label{S:Abstract-liminf-limsup}\label{S:Analysis-abstract}

Our aim is to investigate to what extent a machine learning algorithm can converge to the exact solution of the underlying problem. Such convergence is, of course, only possible if the corresponding neural network spaces are sufficiently expressive. In mathematical terms, this requires that, for every function $u$ belonging to an appropriate function space and for every prescribed tolerance $\varepsilon>0$, there exists a neural network function whose approximation error is smaller than $\varepsilon$ in the norm relevant to the problem under consideration;  Such approximation results are now well established and are commonly referred to as \emph{Universal Approximation Theorems} (UATs).

However, this property alone is not sufficient to capture the notion of approximability that underlies the convergence of an   algorithm. It is also necessary that the approximation quality improves systematically as the complexity of the neural network increases. More precisely, one should be able to construct a sequence of neural network spaces, with an increasing number of parameters, together with corresponding approximations whose error tends to zero as the expressive power of the networks grows. This approximation property is one of the fundamental ingredients in the mathematical analysis of machine learning algorithms and provides the basis upon which convergence results can be established.

%

Following the classical paradigm of numerical analysis, we need to introduce an appropriate notion of \emph{stability} for the algorithm. Stability properties of machine learning algorithms have received comparatively little attention, and identifying a suitable notion of stability in this context is far from straightforward. One first needs to determine the appropriate functional framework within which the convergence analysis can be carried out. Our choice is the De Giorgi core liminf--limsup framework, a fundamental ingredient of $\Gamma$-convergence, a key analytical tool in the calculus of variations {~\cite{DeGiorgi_Franzoni:1975, DalMaso_Gamma:1993}}. 

As we shall see, this approach allows us to establish convergence results without requiring structural assumptions to hold uniformly for all elements of the neural network spaces. In this way, it provides a flexible framework for identifying stable algorithmic procedures while avoiding restrictive assumptions on the underlying network architectures.
%
%
%
%
%
%
%
%
 
Next, we describe the main ingredients needed for our abstract setting. 
 
 
 \subsubsection*{The three core minimisation problems}
  
  The starting point is a generic ML algorithm of  the form
 \begin{equation}\label{ML-FD} 
	\min  _ {v \in  V _{\mathcal{N}} } \mathcal{E}_h (v),   
\end{equation}
where $  V _{\mathcal{N}}$ is the NN space and $\mathcal{E}_h $ is the loss functional. {We use the term liminf--limsup framework rather than
$\Gamma$-convergence on purpose: the standard $\Gamma$-convergence theorems
are not invoked here, because the topology in which the limits are taken
changes with the discretisation, the energies change with the training, and the results on
approximate local minimisers use only the liminf part of the argument.}
 
    Typically there is a corresponding \emph{``continuous" minimisation  problem}, posed on a space $V:$ 
 \begin{equation} \label{ML-Cont} 
	\min  _ {v \in  V  } \mathcal{E} (v) \, .  
\end{equation}
Ideally, minimisers of \eqref{ML-FD} should provide approximations to the minimisers of \eqref{ML-Cont}, which we assume to exist. 
  It will be useful to consider the \emph{(idealised) semi-discrete problem} 
 \begin{equation} \label{ML-SD} 
	\min  _ {v \in  V _{\mathcal{N}} } \mathcal{E}(v)\, .   
\end{equation}
The minimisation problem \eqref{ML-SD} is called idealised, since in practice the loss $\mathcal{E}(v)$ cannot typically be computed exactly over $V_{\mathcal{N}}$, and a further approximation step is required, leading to \eqref{ML-FD}. Notice that, in general, the two sources of discretisation, namely the replacement of $V$ by $V_{\mathcal{N}}$ and the approximation of $\mathcal{E}(v)$ by $\mathcal{E}_h(v)$, are independent and need not be related.

 Next, we aim to justify the following \emph{claim}: if the continuous problem is ``well behaved'' and $\mathcal{E}_h$ is a ``stable and accurate'' approximation of $\mathcal{E}$, then the corresponding machine learning algorithm produces accurate predictions.


\subsubsection*{Expressiveness  of    $V _{\mathcal{N}} $
}

As mentioned above, a fundamental question in the mathematical analysis of machine learning algorithms is the following:

\begin{center}
\emph{How well can elements of the neural network space $V_{\mathcal{N}}$ approximate elements of the underlying function space $V$?}
\end{center}

Throughout this discussion we regard the neural network space as a finite-dimensional approximation of the continuous space,
$$
V_{\mathcal{N}} \subset V.
$$
The approximation properties of the family ${V_{\mathcal{N}}}$ play a central role in determining whether a machine learning algorithm can converge to the exact solution of the underlying problem.

Over the past twenty-five years, significant progress has been made in understanding the approximation capabilities of neural networks. A large body of results, commonly referred to as \emph{Universal Approximation Theorems}, establishes that neural networks are capable of approximating arbitrary functions in a variety of function spaces to any prescribed accuracy, provided that the network architecture is sufficiently expressive. More recently, quantitative approximation results have shown that suitably designed neural network spaces exhibit approximation properties comparable to those of the most advanced finite element spaces, particularly for functions with low-dimensional structure, anisotropic features, or singularities.

Despite this remarkable progress, many important questions remain open. Most available approximation results are of an existential nature: they guarantee the existence of a neural network achieving a prescribed accuracy but provide only limited information about how such a network can be constructed or trained efficiently. Moreover, in contrast to the classical theory of finite element methods and other approximation schemes, where sharp error estimates are available in terms of the regularity of the target function, the geometry of the domain, and the discretisation parameters, comparable quantitative estimates for neural network approximations are still available only in relatively restricted settings.

%

 \subsubsection*{Consistency and stability of the semidiscrete problem}
}
%
%
Next, we assume that
$ 
V_\ell
$ 
is a sequence of increasingly rich neural network spaces that provide progressively more accurate approximations of the underlying function space $V$. We denote by $u_\ell$ a minimiser of the corresponding semidiscrete problem
$$
\min_{v\in V_\ell}\mathcal{E}(v).
$$

The first fundamental question is whether
$$
u_\ell\to u,
$$
where $u$ is a minimiser of the continuous problem
$$
\min_{v\in V}\mathcal{E}(v).
$$

This question is of central importance. Although, as we shall see later, it does not by itself provide a complete convergence theory, an affirmative answer is a necessary prerequisite for the reliable performance of the learning algorithm. Indeed, if the semidiscrete minimisers fail to converge even when the exact functional $\mathcal{E}$ is minimised over increasingly expressive neural network spaces, then no subsequent training or numerical approximation can be expected to recover the correct solution. The question remains relevant even though, in practice, the optimisation
procedure rarely returns a global minimiser. Indeed, the behaviour of the
global minimisers of the exact functional isolates the roles of the energy
and of the neural network spaces from that of the optimisation, and is the
natural benchmark against which the computed approximations are to be
judged.

Our main claim is that, for a broad class of machine learning models, the answer to this question depends essentially only on the stability properties of the \emph{continuous} functional $\mathcal{E}$, and not on the nonconvex nature of the neural network parameterisation. More precisely, whenever $\mathcal{E}$ satisfies suitable convexity or variational stability assumptions, one can establish the convergence
$$
u_\ell\to u.
$$

A particularly striking aspect of this result is that \emph{the intrinsic nonconvexity of the semidiscrete optimisation problem does not obstruct the convergence analysis.} Instead, it is effectively bypassed through variational arguments at the level of the functionals, showing that the convergence mechanism is governed primarily by the structure of the continuous problem rather than by the nonlinear parameterisation of the neural network spaces.\\

 \subsubsection*{Stability / the role of training / quadrature 
}
The final computations are performed using a computable loss functional
\begin{equation*} 
	  \mathcal{E}_h (v)   \approx   \mathcal{E}  (v) 
\end{equation*}
 We now denote by $u_\ell$ the corresponding minimiser of
 \begin{equation*} 
	\min  _ {v \in  V _{\ell} } \mathcal{E}_{h(\ell)}  (v) \, .  
\end{equation*}

The resulting question is whether we can show that
$$
u_\ell\to u,
$$
where $u$ is a minimiser of the continuous problem
$$
\min_{v\in V}\mathcal{E}(v).
$$
The answer is affirmative, provided that $\mathcal{E}_h(v)$ is a ``stable and accurate'' approximation of $\mathcal{E}(v)$. Clearly, $\mathcal{E}_h(v)$ must approximate $\mathcal{E}(v)$ on an appropriate class of elements of $V$ containing the exact solution. The precise meaning of stability, however, requires a more careful formulation taking into account the fact that the variational properties of $\mathcal{E}_h$ might change as $h$ varies: 

\noindent
\emph{\textbf{Discrete Variational Stability.}
A family of discrete energies $\mathcal{E}_h$ is called $V$-stable if there exist discrete seminorms  $\| \cdot\| _{V_h}$  such that:}
 \begin{enumerate}
  	  \item   $\mathcal E_h (v) \geq \alpha \| v\|^{\beta} _{V_h}  {- C}  \quad (\text{{discrete coercivity}}) $
  \item   \emph{uniformly bounded sequences in  $\| v_h\| _{V_h}$ have  convergent subsequences in a weaker topology $H$ of $V$} 
 \item \emph{the lower semicontinuity property holds: If $v_h\rightarrow v$ in $H$ as $h\to 0$ then} 
$$
\mathcal{E}(v)\leq
\liminf_{h\rightarrow 0}
\mathcal{E}_{h}(v_h).
$$
\end{enumerate}
In other words, we require a suitable form of discrete compactness that
guarantees the existence of convergent subsequences together with the
lower semicontinuity needed to identify the limiting minimiser. The
relative weight of the three requirements varies with the application: in
some cases coercivity is not needed, compactness being available for other
reasons, while in others lower semicontinuity follows immediately from
coercivity and convexity.

As in the semidiscrete problem, the discrete optimisation problems may remain highly nonconvex due to the neural network parameterisation. Nevertheless, the essential variational properties, such as lower semicontinuity and compactness, are inherited from the underlying continuous energies and remain the key ingredients in the convergence analysis.

This viewpoint contrasts with approaches that focus primarily on the nonconvexity of the discrete optimisation problem while overlooking the variational structure and stability properties of the underlying continuous functional. The analysis suggests that, for convergence purposes, the properties of the parent continuous energy may be more fundamental than the nonconvex character of the discrete neural network optimisation problem.

\subsubsection*{A Numerical Analysis perspective on the liminf–limsup framework}

{As explained at the end of this subsection, a uniform bound for the energies of the discrete minimisers follows from the limsup inequality introduced below.}
 Therefore,   if the method is stable, meaning it satisfies the \emph{Discrete Variational Stability} property discussed above, we can conclude that there exists an  $   \tilde u\in V$ 
such that $   u_\ell     \rightarrow   \tilde  u$ in 
$H $ up to a subsequence not relabelled.

Next, to fix ideas, assume that the continuous  problem \eqref{ML-Cont} has a unique solution $u$ and  we   prove that $ \tilde  u$ is a global minimiser of $ \mathcal{E}  , $ i.e., $ \tilde  u=u,$ and the convergence  of the entire sequence to $u.$

To this end, \emph{assume} that given an element $w\in V$ we can construct a sequence $w_\ell \in V_\ell$ such that 
\begin{align}\label {limsup-NA}
   \limsup_{\ell  \rightarrow \infty}  \mathcal{E}  _{\ell}[ w_\ell] \le  \mathcal{E} [  w ],  
\end{align}
where $ \mathcal{E}  _{\ell} :=\mathcal{E}_{h(\ell)}\, $ and $h(\ell) \to 0$ as 
$\ell \to \infty\, .$ In fact, in most applications one can show $ \lim_{\ell  \rightarrow \infty}  \mathcal{E}  _{\ell}[ w_\ell] =  \mathcal{E} [  w ].$ Then 
by the fact that $ u_\ell $ are   minimisers of the discrete problems we finally obtain
\begin{align}\label{abstract-linflsup}
 \mathcal{E} [   \tilde  u ] \le  \liminf_{\ell \rightarrow \infty}  \mathcal{E} _{\ell}[ u_\ell ]  \le  \limsup_{\ell  \rightarrow \infty}  \mathcal{E} _{\ell}[ u_\ell] 
\le  \limsup_{\ell  \rightarrow \infty}  \mathcal{E}  _{\ell}[ w_\ell] \le  \mathcal{E} [  w ], 
\end{align}
for all $ w \in V$. 
Therefore $  \tilde u $  is an absolute minimiser of $ \mathcal{E}  , $ and $  \tilde u =u.$  Since all weakly convergent subsequences of $ \{ u_{\ell} \, \}$ can converge only to $u,$ the entire sequence converges to $u.$  \\

\noindent 
\textbf{Stability and approximability.} The above argument can be interpreted from a Numerical Analysis viewpoint as another application of the principle 
\begin{center}
{\emph{stability $+$ consistency/approximability $\implies $ convergence.}}
\end{center}
	 In fact, stability in the sense of \emph{Discrete Variational Stability} implies the existence of $\tilde u$ and the liminf inequality stated above. As we shall see later in this section, although this may not be apparent at first sight, this is indeed an appropriate notion of stability for the minimisation problems under consideration. It captures, uniformly with respect to the discretisation parameters, the boundedness of the corresponding solution operators {, which is precisely the notion of stability of Lax and Richtmyer discussed in Section~\ref{SubS:aPaP:Error_analysis}}. Moreover, whenever a connection with classical notions of stability can be established, it provides a clear link between the two frameworks. In particular, methods that fail to satisfy this stability property exhibit poor pre-asymptotic behaviour and, in numerical experiments,  appear to diverge.

	The limsup inequality provides a clear analogue of the classical notions of consistency and approximability in numerical analysis. It applies to all \emph{sufficiently smooth} elements of $V$ and, in particular, to the exact solution $u$. More precisely, it is sufficient to construct, for every such $w\in V$, a sequence $w_\ell\in V_\ell$ such that
$$
\lim_{\ell  \rightarrow \infty}  \mathcal{E}  _{\ell}[ w_\ell] =  \mathcal{E} [  w ].
$$
The sequence $w_\ell$ can be viewed as a suitable discrete approximation or recovery sequence for $w$. Thus, the limsup condition expresses the requirement that the discrete spaces and discrete energies are sufficiently rich and accurate respectively, to reproduce, in the limit, the energy of the continuous problem.

The analogy with classical numerical methods is particularly transparent. For example, in finite difference methods, if $\mathcal{E}$ represents a differential operator and $\mathcal{E}_\ell$ its discrete approximation, consistency requires that the discrete operator applied to suitable discrete representations $w_\ell$ of a smooth function $w$ converges to the corresponding continuous quantity. Typically, $w_\ell$ is obtained by evaluating $w$ at the grid points. In finite element methods, the corresponding sequence $w_\ell$ is usually constructed by means of an appropriate interpolation or projection operator onto the finite element space.

The limsup condition therefore plays the role of an abstract \emph{consistency} requirement that is independent of the particular mechanism used to construct the discrete approximation. In the present setting, this is particularly useful since the spaces $V_\ell$ are neural network spaces and the appropriate approximation or recovery sequence need not be given by a conventional interpolation or projection operator. What is essential is the existence of a sequence of neural network functions whose discrete energies converge to the energy of the target function.

\subsubsection*{What the framework delivers}

Although the discussion above is deliberately informal, with the precise meaning of each ingredient supplied by the application at hand, the framework it describes has fixed consequences. First, the uniform bound on the energies of the discrete minimisers, which starts the argument, is not an additional assumption: it follows from the limsup inequality applied to a single recovery sequence, for instance that of the exact solution. Second, the two parts of the framework can be used separately. The liminf inequality alone, applied to any sequence with uniformly bounded losses, shows that every accumulation point has energy bounded by the computed losses; this is an a posteriori statement, and it applies to the approximate local minimisers that training actually produces. Combined with the limsup inequality, it shows that accumulation points are minimisers of $\mathcal{E}$, and when only one discretisation parameter tends to infinity, for instance the number of samples with the network space fixed, their error is bounded by the best approximation error of that space. Both statements are made precise, for function learning, in Section~\ref{SubS:Conv:MCInterpolation}, and, for the generation of probability measures, in Section~\ref{SubS:Conv:measures}.

Implicit in this is a partial answer to the question of why machine learning algorithms work. If, for a given application, the requirements of the framework are met, that is, the loss retains the coercivity and lower semicontinuity of the continuous functional and the network spaces are expressive enough to provide recovery sequences, then reliable predictions are to be expected, and small computed losses (or, more generally, losses approaching the minimal value of the exact energy) verify this fact. Notice that the good behaviour is established at the level of the function space, where the energies live, and not at the level of the parameter space, where the optimisation takes place: the nonconvexity in parameter space, discussed in Section~\ref{S:typical}, is not relevant here. The answer is partial because this approach does not provide information about whether the optimisation procedure will find minimisers with small loss; however, if the procedure does, the framework provides guarantees for the
reliability of the resulting approximations.

Finally, nothing in the approach requires $V$ to be a space of functions or $V_{\mathcal{N}}$ to be a linear subspace of it. A set equipped with a weak notion of convergence, loss functionals with the coercivity and lower semicontinuity properties above, and recovery sequences are the main ingredients. This is what allows the same approach to be followed for spaces of probability measures, in the generative models of Section~\ref{SubS:Ex:measures}, and for operator learning, where $\Omega$ is replaced by a space of functions.


\subsection{A simple application: function learning with deterministic interpolation}
 
Recall the function learning problem  \eqref{mm_nn:abstract_theta} and the corresponding  the \emph{continuous} minimisation problem 
\begin{equation}\label{mm_nn:abstract_cont_2}
	\min  _ {v \in  V   } \mathcal{E} (v)  , 
\end{equation}
where $V = L^p (\Omega)$ and
\begin{equation*}
\mathcal{E} (u)=
\int_{\varOmega}    |  u - f|^p    \d x  .
\end{equation*}
The functional $\mathcal{E}$ is  convex and the unique minimiser of \eqref{mm_nn:abstract_cont_2} is obviously $f\, .$
The semidiscrete problem is 
%
\begin{equation}\label{mm_nn:abstract_semd}
	\min  _ {v \in  V _{\mathcal{N}} } \mathcal{E} (v)   
\end{equation}
%
and we denote by $v^{*} \in V _{\mathcal{N}}  $  a $\mathcal{N}-$DNN minimiser of $\mathcal{E} \, . $ 
Before considering further discretisations we shall examine whether minimisers of the \emph{semidiscrete problem} \eqref{mm_nn:abstract_semd} converge to $f\, .$ First we would like to pose this question, on a sequence of neural network spaces $V_\ell$ which have been chosen such that for each $\ell$ there exists a $w_\ell\in V_\ell$ such that as $\ell \to \infty, $
\begin{equation}\label{mm_nn:Vell-appr-semd}
	w_\ell \to f \quad \text {in } L^p (\Omega) \, .
\end{equation}
In this case, the convergence is immediate in view of the direct estimate
\begin{equation}\label{mm_nn:Vell-appr-semd_2}
	\int_{\varOmega}    |  u_\ell - f|^p    \d x =  \mathcal{E} [ u_\ell] 
\le  \mathcal{E}   [ w_\ell] \, .
	\end{equation}
As we shall see, it is not always possible to extend this estimate when training takes place, i.e., when $\mathcal{E}$ is replaced by a computable approximation $\mathcal{E}_h.$ 

To see why the approach described in Section \ref{S:Abstract-liminf-limsup} is readily applicable, notice that the minimisers have uniformly bounded energies, since, 
\begin{equation}\label{mm_nn:Vell-appr-semd_3}
	   \mathcal{E} [ u_\ell] 
\le  \mathcal{E}   [ w_\ell] \leq C \, ,
	\end{equation}
in view of \eqref {mm_nn:Vell-appr-semd}. Then there exists a subsequence, still denoted by $\{ u_\ell \}$, and an element $\tilde u $ such that     $  u_\ell     \rightharpoonup    \tilde  u $ in $L^p(\Omega) .$
Furthermore, since $\mathcal {E} $ is convex, 
  weak lower semicontinuity   implies
 $$\liminf\limits_{\ell \to \infty}
\int_\Omega |  ( u_\ell  -f)  |^p \ge \int_\Omega |\tilde u-f|^p.$$
These arguments confirm that in this case the  \emph{Discrete Variational Stability} holds with $H$ being the topology of weak $L^p$-convergence. The inequalities in  \eqref{abstract-linflsup} imply 
 $\mathcal{E} [   \tilde  u ] \le    \mathcal{E} [  f ]=0$, 
and thus  $  \tilde u =f.$  Since all weakly convergent subsequences of $ \{ u_{\ell} \, \}$ can converge only to $f,$ the entire sequence converges to  {$f$}.  {For global minimisers the direct estimate \eqref{mm_nn:Vell-appr-semd_2} already gives strong convergence; the weak argument is nevertheless the one that survives training and applies to approximate local minimisers, as we shall see.}

 \subsubsection*{Deterministic Training}

 As discussed in Section \ref{SubS:ModP:function-learning} we approximate $\mathcal{E}_h(v) \cong \mathcal{E}(v)$ by assuming that   $f$ is known  at a finite set of points 
 $z\in K_h,$ where $K_h$ is a set of points linked to a quadrature rule.  
 Then the loss functional  %
\begin{equation*}
\label{E_h-d4}
\mathcal{E}_{Q, h}( g )  =     \sum _{z\in K_h} \, w_z \, \Big  | g(z) - f(z) \Big  | ^p  \, , 
\end{equation*}
leads to the minimisation problem  
\begin{equation}\label{mm_nn:EQh}
	\min  _ {v \in  V _{\mathcal{N}} } \mathcal{E} _ {Q, h}(v)\, .
\end{equation}
 Problem \eqref{mm_nn:abstract_semd} is the semidiscrete analogue of \eqref{mm_nn:EQh}. 
As we shall see, a crucial observation is that 
 there exist appropriate interpolation   operators $ I_ {{\mathbb{V}}_h}  : V_S \to L^p (\Omega),$ where here 
 $V_S\subset V$ consists of   functions with appropriate smoothness, such that
 \begin{equation*}
\label{E_h_int}
\mathcal{E}_{Q, h}( g )  =  \int _\Omega  \,I_ {{\mathbb{V}}_h} \Big ( \big  |  g -     f   \big  |^p \Big )(x) \d x \, . 
\end{equation*}
An alternative loss functional which still uses interpolation at the   points   $z\in K_h,$ is 
\begin{equation}\label{tildeEQh}
\tilde {\mathcal{E}} _{Q, h}( g )  =  \int _\Omega  \, \Big  | I_ {{\mathbb{V}}_h} \Big (g(x) -     f(x) \Big ) \Big  | ^p \d x \, . 
\end{equation}
Still, $\tilde {\mathcal{E}}_{Q, h} \cong \mathcal{E}(v),$ and the new loss, allows more general choices for the operators $ I_ {{\mathbb{V}}_h}\, .$ The corresponding minimisation problem is 
\begin{equation}\label{mm_nn:tildeEQh}
	\min  _ {v \in  V _{\mathcal{N}} } \tilde {\mathcal{E}} _ {Q, h}(v)\, .
\end{equation}
The analysis of the problem \eqref{mm_nn:EQh} is postponed to the Section \ref{SubS:Anal:functionQh}. 
The main steps of the convergence for the problem \eqref{mm_nn:tildeEQh} are as follows: 
 First we would like to pose this question, on a sequence of neural network spaces $V_\ell$ which have been chosen such that for each $\ell$ there exists a $w_\ell\in V_\ell$ such that as $\ell \to \infty, $ and $h\to 0, $
\begin{equation}\label{mm_nn:tildeVell-appr-semd}
	I_ {{\mathbb{V}}_h} w_\ell \to f \quad \text {in } L^p (\Omega) \, .
\end{equation}
We also require that as   $h\to 0, $
\begin{equation}\label{mm_nn:tildeVell-appr-semd-2}
	I_ {{\mathbb{V}}_h} f \to f \quad \text {in } L^p (\Omega) \, .
\end{equation}
For convenience we link the discretisation parameters by choosing $h=h(\ell)$ where $h(\ell) \to 0$ as $\ell \to \infty .  $ We denote by $u_\ell$ a minimiser of \eqref{mm_nn:tildeEQh} written as 
\begin{equation}\label{mm_nn:tildeEQh-ell}
	\min  _ {v \in  V _{\ell} } \tilde {\mathcal{E}} _ {Q, h(\ell)}(v)\, .
\end{equation}
The approach described in Section \ref{S:Abstract-liminf-limsup} can be  applied in this setting. 
First we observe that  
\begin{equation}\label{mm_nn:tildeVell-appr-semd_3}
	   \tilde {\mathcal{E}} _ {Q, h(\ell)} [ u_\ell] 
\le    \tilde {\mathcal{E}} _ {Q, h(\ell)}  [ w_\ell] \leq C \, ,
	\end{equation}
in view of \eqref {mm_nn:tildeVell-appr-semd} and \eqref{mm_nn:tildeVell-appr-semd-2}. Then, taking into account \eqref{mm_nn:tildeVell-appr-semd-2},  there exists a subsequence, still denoted by $\{ u_\ell \}$ and an element $\tilde u \in V$ such that     $ I_ {{\mathbb{V}}_{h(\ell)}} ( u_\ell  -f)   \rightharpoonup    (\tilde u-f)$ in $L^p(\Omega) .$
By  
  weak lower semicontinuity   we get 
 $$\liminf\limits_{\ell \to \infty}
\int_\Omega | I_ {{\mathbb{V}}_{h(\ell)}} ( u_\ell  -f)  |^p \ge \int_\Omega |\tilde u-f|^p.$$
These arguments confirm that in this case the  \emph{Discrete Variational Stability} holds with $H$ being the topology of weak $L^p$-convergence and choosing the discrete seminorms as 
$ \| v \|_{V_h} = \|  I_ {{\mathbb{V}}_{h }}v\| _{L^p}. $ Then, in view of \eqref {mm_nn:tildeVell-appr-semd} and \eqref{mm_nn:tildeVell-appr-semd-2},  the inequalities in  \eqref{abstract-linflsup} imply 
 $\mathcal{E} [   \tilde  u ] \le    \mathcal{E} [  f ]=0$, 
and thus  $  \tilde u =f.$  Therefore the sequence  $ \{ u_{\ell} \, \}$ satisfies, $ I_ {{\mathbb{V}}_{h(\ell)}}  u_{\ell} \to f \, .  $

\subsection{Function learning through MC Interpolation}\label{SubS:Conv:MCInterpolation}
 \def\E{{\mathcal{E}}}

We consider now the typical function learning problem. Our aim is to investigate its convergence properties from various perspectives, and without any structural hypotheses on the nonlinear parametric spaces other than their expressiveness/approximability. 
As introduced earlier, the problem of MC-Deep neural network interpolation  is to optimise
\begin{equation}\label{ieE-minimize_NNQ_prob-2}
	\min  _ {v \in  V _{\mathcal{N}} } \E _{ {N,\omega}}(v)\, . 
\end{equation}
The minimisers are denoted by  $u_X \in V _{\mathcal{N}}$ and  $X=(X_1,X_2, \dots, X_N)$ are i.i.d.\ random variables, each one corresponding to a random point $x_i$ of $\Omega.$ 
A key mathematical challenge is that  
  although for a fixed function $g$ the error
  \begin{equation*}
\frac 1 N \sum _{i=1} ^N \,  \, g(X_i(\omega)\,) - \int _{\Omega} \, g(x) \, \d x ,
\end{equation*}
can be controlled using standard tools, these are not applicable to control
\begin{equation*}
\frac 1 N \sum _{i=1} ^N \,  \, u_X(X_i(\omega)\,) - \int _{\Omega} \, u_X(x) \, \d x ,
\end{equation*}
since $u_X$ depends on the random variables $X_i$ through \eqref{ieE-minimize_NNQ_prob-2}, and thus they change as we increase the resolution of the Monte-Carlo quadrature.

%
%
%

We are   interested to study the behaviour of the generalisation error:
$$ \|  f- u_{\bm{X}     }\| _p\, .$$
It will be useful to consider the application of quadrature as an empirical measure approximation of the continuous integral with respect to the measure $\mu\, .$ 
Indeed one may rewrite the loss functional for $ \omega\in \bbOmega ,$
\begin{equation*}\label{prob_E-d2}
\mathcal{E}_{N, \omega}( g )  =  \frac 1N  \sum _{i}   \Big  | g(X_i(\omega)\, ) - f(X_i (\omega) \, ) \Big  | ^p  \, 
\end{equation*}
using  the empirical measures on $\Omega,$
\begin{equation*}\label{emp_meas-d2}
	\mu _{N, \bm{X} (\omega)} =   \frac 1 N \sum _{i} \   \delta _{X_i(\omega)} \, ,
\end{equation*}
as
\begin{equation*}\label{energy_emp_meas}
	\mathcal{E}_{N, \omega}( g )  = \int _\Omega |g (x)-f(x)|^p\, \d \mu _{N, \bm{X} (\omega)} (x) 
\end{equation*}
then
\begin{equation*}\label{prob_E_error-d4}
 {E}_{MC} (u_{\bm{X} (\omega)    }) =  \int_{\Omega}    |  u_{\bm{X}  (\omega)   } - f|^p    \d x -  \int _\Omega | u_{\bm{X}  (\omega)   }-f|^p\, \d \mu _{N, \bm{X} (\omega)}   \, . 
\end{equation*}

%

\subsubsection{Weak convergence of measures and functions}\label{SubS:weakconv}
Interpreting quadrature through measure approximation naturally provides tools to study the limiting behaviour  of discrete minimisers. We shall use therefore standard notions of weak convergence of measures. 
For consistency we follow the notation and terminology of \cite{AGS:2005}:
	\begin{itemize}
		\item A sequence of measures $\nu_n $ on $D$ is said that \emph{weakly converges} to $\nu$
		
		if for all bounded and continuous functions, $g : D \to \R$ we have 
 \begin{equation*}
 \begin{split} 
\lim _{n\to \infty} \int _D \,  g (x) \, \d \nu_n (x) = \int _D \,  g (x) \, \d \nu (x) \, . 	
\end{split}
\end{equation*}

\item A sequence of functions $v_n \in L^1 (D, \nu_n)$  is said that \emph{weakly converges to 
$v  \in L^1 (D, \nu)$ with respect to $(\nu_n, \nu), $} \cite[Definition 5.4.3]{AGS:2005} if

(a) $\nu_n $ weakly converges to $\nu$

(b) and

\begin{equation*}\label{w-conv-u_n0}\begin{split} 
\lim _{n\to \infty} \int _D \, \zeta (x) v_n(x) \, \d \nu_n (x) = \int _D \,  \zeta (x) v(x)  \, \d \nu (x) \, , \qquad\text{for all }\zeta \in C_0^\infty (D)\, .
\end{split}
\end{equation*}

	\end{itemize}
We also say that $v_n$ strongly converges to $v$ in $L_p,\ p > 1$, with respect to $(\nu_n, \nu), $ if in addition, 
\begin{equation*} 
\limsup _{n\to \infty}\|v_n\|_{L_p(D,\nu_n)} \leq \|v\|_{L_p(D,\nu)} .
\end{equation*} 

A standard useful starting point is that  the empirical measures weakly converge to $\mu$ in the sense that, \cite{V-58}, 
\begin{equation*}\label{perklln}
 \mu _{N, \bm{X}_N (\omega  )} \to \mu, \qquad \text{ for all }\omega \in \bbOmega _1,
 \end{equation*}
where
 $X_1,X_2,\ldots$ is a sequence of i.i.d. random variables with common law $\mu$ and  $\bbOmega _1\subset \bbOmega $ with $\P ( \bbOmega \backslash  \bbOmega _1)=0.$ 
This follows 
since there exists a fixed countable subset $\mathcal {G}$ of the set of bounded continuous functions on $D$, and for $g\in \mathcal {G},$ there holds  
\begin{equation*} \begin{split}
\int_\Omega g(x)\, \d \mu _{N, \bm{X}_N (\omega  )}(x)= \frac 1N  \sum _{i=1}^N    g (X_i)    \,   \longrightarrow  \int_\Omega g(x)\, \d \mu(x)\, , \qquad \mathbb{P}-\text{a.s.},	
\end{split}
\end{equation*}
by the strong law of large numbers. Since any countable union of sets of measure zero has measure zero, the assertion follows.

 \subsubsection{Stability via $\liminf$ inequality}
	
	The following result forms the backbone of the stability analysis of the probabilistic learning algorithm. It relies on the fact that the empirical measures converge weakly to $\mu$ and yields several intermediate results concerning not only global minimisers, but also approximate local minimisers, provided that their losses remain uniformly bounded.

	{For a fixed DNN space $V_{\mathcal{N}}$ the problem under consideration is}
\begin{equation}\label{prob_E_def_theorem0}
	\min  _ {v \in V _{\mathcal{N}} } \E _{ {N } , \omega} (v)\, ,
\end{equation}
{and we recall that no minimality property is used in the following lemma.}
\begin{lemma}[\sc the $\liminf$ inequality]
\label{lim-inf-lemma}
Let $p>1$, $f\in L^p(\Omega,\mu)$, and $V _{\mathcal{N}}$ a fixed DNN space. For
$\omega\in \bbOmega$, let $\big(u_{\bm{X}_{N} (\omega)}\big)_N$ be a sequence of
elements of $V _{\mathcal{N}}$, where $ \bm{X} _{N} = (X_1, \dots, X_{N}  )$, and
assume that there exists a subset $\bbOmega _2\subset \bbOmega$ with
$\P ( \bbOmega \backslash  \bbOmega _2)=0$ such that
 \begin{equation}\label{w-conv-u_n_bound0}\begin{split}
   \E _{N, \omega}  \big( u_{\bm {X}_N (\omega)} \big) \leq M (\omega)<+\infty\, ,
   \qquad \text{ for all  }\ \omega \in \bbOmega _2.
\end{split}
 \end{equation}
Then, for $\mathbb{P}$-a.e.\ $\omega\in\bbOmega$, the sequence
$(u_{\bm{X}_{N} (\omega)   } )_N$ has weak limits along subsequences with respect
to $ (\mu _{N, \bm{X}_{N} (\omega  )}, \mu) $, and every such weak limit
$ u _{\omega} ^  {\mathcal{N}}$ belongs to $L^p(\Omega,\mu)$ and satisfies
 \begin{equation*}\label{w-conv-u_n_0}\begin{split}
\E \big(u_\omega^  {\mathcal{N}}\big) \leq \liminf\limits_{N\rightarrow \infty}  \E _{{N}, \omega}  \big( u_{\bm {X}_{N} (\omega)} \big)\, .
\end{split}
 \end{equation*}
\end{lemma}

\begin{remark}
No minimality property of $u_{\bm{X}_{N} (\omega)}$ is assumed: the lemma applies
to global minimisers of \eqref{prob_E_def_theorem0} as well as to approximate
local minimisers, as long as their losses are uniformly bounded in the sense
of \eqref{w-conv-u_n_bound0}. Likewise, the DNN structure of $V _{\mathcal{N}}$
plays no role in the proof, which only uses that the $u_{\bm{X}_{N} (\omega)}$
are measurable functions with bounded losses; the conclusion therefore holds
verbatim when the spaces vary with $N$, i.e.\ for sequences
$u_{\bm{X}_{N} (\omega)}\in V_{\mathcal{N}(N)}$. This is the form in which the
lemma is used in Theorem~\ref{learningMC_conv}.
\end{remark}

The proof is an adaptation to our present setting of
\cite[Theorem 5.4.4]{AGS:2005}. The main steps are the following, where we
set $v_N(\omega):=u_{\bm {X}_N (\omega)}-f$ and write $\mu_N$ for
$\mu _{N, \bm{X}_N (\omega)}$.

\begin{itemize}
\item By the boundedness assumption we have
 \begin{equation*}\label{energy_emp_meas1}
 \int _\Omega |v_N(\omega) (x)|^p\, \d \mu _{N} (x)=\E _{N, \omega}  \big( u_{\bm {X}_N (\omega)} \big) \leq M (\omega ),  \quad \omega \in \bbOmega _2.
 \end{equation*}

\item The sequence of measures $\big([{v_N (\omega)}]_\sharp \mu _{N}\big)_N $
on $\R$ has $p$-th moments bounded by $M(\omega)$, by the previous step; it
is therefore tight, by Markov's inequality, and relatively compact with
respect to weak convergence, by Prokhorov's theorem.

\item Define the plans  $\bm{\gamma} _N  = (\bm{i} \times  {v} _N )_\sharp \mu _{N}$
on $\Omega\times\R$, where $\bm{i}: \Omega \to \Omega $ is the identity map. The
first marginal of $\bm{\gamma} _N $ is $\pi ^1_\sharp \bm{\gamma} _N = \mu _{N} $
and the second is $\pi ^2 _\sharp \bm{\gamma} _N =  ({v} _N )_\sharp  \mu _{N}$.

\item For $\omega\in\bbOmega_1$ the first marginals $\mu_N$ converge weakly to
$\mu$, and for $\omega\in\bbOmega_2$ the second marginals are tight; since a
sequence of plans whose marginals are tight is itself tight, $(\bm{\gamma} _N)_N$
is relatively compact for $\omega\in\bbOmega_1\cap\bbOmega_2$, a set of full
measure. Let $\bm{\gamma}$ be the limit of a subsequence, not relabelled.

\item The limit points $\bm\gamma$ need not be concentrated on the graph of
a function: if the $v_N$ oscillate, the disintegration
$\bm\gamma=\int_\Omega\bm\gamma_x\,\d\mu(x)$ of $\bm\gamma$ with respect to its
first marginal, which is $\mu$, has fibres $\bm\gamma_x$ that are not Dirac
masses. The weak limit of $v_N$ is the \emph{barycentre} of the fibres,
 \begin{equation*}
 \overline {\bm{\gamma}}_1  (x)  := \int _{\R }\, y\,  \d \bm{\gamma}_{x} (y)\, ,
 \end{equation*}
which is well defined for $\mu$-a.e.\ $x$ since
$\int_{\Omega\times\R}|y|\,\d\bm\gamma<\infty$ by the moment bound.

\item Indeed, for all $\zeta \in C_0^\infty (\Omega)$, the definition of the
push-forward gives
$$\int _\Omega   \zeta (x) \,v_{N} (x) \,\d \mu _{N} (x) = \int _{\Omega \times \R}   \zeta (x)\,  y \,   \d \bm{\gamma} _{N} (x, y).$$
The integrand $\zeta(x)\,y$ is continuous but unbounded in $y$, so weak
convergence alone does not allow passage to the limit; however, since
$\int|y|^p\,\d\bm\gamma_N\le M(\omega)$ with $p>1$, the function $y$ is uniformly
integrable with respect to $(\bm\gamma_N)_N$, and
\cite[Lemma 5.1.7]{AGS:2005} yields
 \begin{equation*}\label{w-conv-u_n_00}\begin{split}
\lim _{N\to \infty} \int _\Omega  \, \zeta (x) \,v_{N} (x) \,\d \mu _{ {N}} (x)
= \int _{\Omega \times \R}   \zeta (x)\,  y \,   \d \bm{\gamma}   (x, y)
= \int _{\Omega}   \zeta (x)\,     \overline {\bm{\gamma}}_ 1(x)  \,  \d  \mu  (x)  \, ,
\end{split}
 \end{equation*}
the last equality by the disintegration. This is the weak convergence of
$v_N$ to $\overline {\bm{\gamma}}_ 1$ with respect to $(\mu_N,\mu)$. This is the
only step where $p>1$ is used.

\item To pass from $v_N$ to $u_{\bm{X}_N(\omega)}=v_N+f$, note that
$\int_\Omega \zeta f\,\d\mu_N\to\int_\Omega\zeta f\,\d\mu$ for $\mathbb{P}$-a.e.\
$\omega$, by the strong law of large numbers applied to $\zeta f\in L^1(\Omega,\mu)$;
as in the proof of the weak convergence $\mu_N\to\mu$, the exceptional set can
be chosen independently of $\zeta$ by taking $\zeta$ in a countable dense
subset of $C_0(\Omega)$ and extending by density. Setting
$ u_\omega^{\mathcal{N}} (x) :=    \overline {\bm{\gamma}}_ 1 (x) +f (x)$, we obtain the
weak convergence of $u_{\bm{X}_N(\omega)}$ to $u_\omega^{\mathcal{N}}$ with respect
to $(\mu_N,\mu)$ along the subsequence.

\item Finally, the lower semicontinuity of $\bm\gamma\mapsto\int|y|^p\,\d\bm\gamma$
with respect to weak convergence (the integrand is continuous and
nonnegative), the disintegration of $\bm\gamma$, and Jensen's inequality for
the convex function $|\cdot|^p$ give
 \begin{equation*}\label{w-conv-u_n_4}\begin{split}
 \liminf\limits_{N\rightarrow \infty}   \E _{{N}, \omega}  \big( u_{\bm {X}_{N} (\omega)} \big) &=
  \liminf\limits_{N\rightarrow \infty}  \int _{\Omega \times \R}   \, | y | ^p\,   \d \bm{\gamma} _{N  } (x, y) \geq
   \int _{\Omega\times \R}    | y | ^p \,   \d \bm{\gamma}   (x, y)  \\
   &=\int_\Omega \Big( \int_{\R}\, |y|^p \,\d\bm{\gamma}_x(y)\Big)\,\d \mu(x)
   \geq \int _{\Omega}       |  \overline {\bm{\gamma}}_ 1(x)  |^p \, \d  \mu  (x) = \E \big(u_\omega^{\mathcal{N}}\big) \, ,
\end{split}
 \end{equation*}
which is the asserted inequality; in particular $u_\omega^{\mathcal{N}}-f\in L^p(\Omega,\mu)$,
hence $u_\omega^{\mathcal{N}}\in L^p(\Omega,\mu)$.
\end{itemize}
 
\subsubsection{Approximability and robustness as $N\to \infty$} 

We turn our attention to the $N\to \infty$ limit. A consequence of the liminf property is that these limits are remarkably robust. This applies to both global minimisers and approximate local minimisers. 

\subsubsection*{Global minimisers} 	
	  
In the case of global minimisers the expressiveness of 	 $V _{\mathcal{N}}$ through 
\begin{equation*}\label{infCeaV_N}
  \inf_{\varphi   \in  V _{\mathcal{N}} } \|\varphi-f\|_p^p\, ,
   \end{equation*}	 
controls the generalisation error even at the limit as $N\to \infty.$ In fact, 
let  $V _{\mathcal{N}}$ be a fixed DNN space. For $\omega\in \bbOmega$ consider a sequence $ \big(u_{\bm{X}_{N} (\omega) }\big)_N    ,$ of minimisers of 
\eqref{prob_E_def_theorem0}. Then, 
 \begin{equation*}\label{limsupbound}
\limsup_{N\to\infty}\E _{ {N } , \omega} (u_{\bm{X}_{N } (\omega)   }  )\leq  \inf_{\varphi   \in  V _{\mathcal{N}} } \|\varphi-f\|_p^p\, ,\qquad\text{for}\ \mathbb{P}-\text{a.e. } \omega\in\bbOmega.
 \end{equation*}
Indeed this is a consequence of the minimising property of $u_{\bm{X}_{N} (\omega)},$ since obviously, 
 \begin{equation*}\label{prob_E_theorem_2_10}
	\E _{ {N } , \omega} (u_{\bm{X}_{N } (\omega)   }  ) \leq \E _{ {N } , \omega} (v)\, , \quad \text {for all } \ v \in  V _{\mathcal{N}}  \, . 
 \end{equation*}
Then, for a sequence $(v_n)_n$ in $V _{\mathcal{N}}$ that realises the infimum, i.e.,
\[
\lim_{n\to\infty} \|v_n-f\|_p=\inf_{\varphi   \in  V _{\mathcal{N}} } \|\varphi-f\|_p.
\]
%
we can finally prove that 
 \begin{equation*}\label{w-prob_E_theorem_3}
\limsup_{N\rightarrow \infty}   \E _{N, \omega}  ( u_{\bm {X}_{N} (\omega)} )   
 \leq \inf_n\|v_n-f\|_p^p =\inf_{\varphi\in V_\mathcal{N}} \|\varphi-f\|_p^p  \, .
 \end{equation*}

In the next result, we keep $V_{\mathcal{N}}$ fixed and let $N\to\infty$. Remarkably, the minimisers need not converge to a single approximation of $f,$ but nevertheless, all possible subsequential limits have a generalisation error controlled by the best approximation error within $V_{\mathcal{N}}$. Hence, \emph{all such limits provide good approximations of $f$, provided that $V_{\mathcal{N}}$ is sufficiently expressive}.
 
	\begin{theorem}[\sc limit behaviour of minimisers as $N\to \infty$] 
Let $p>1$ and $V _{\mathcal{N}}$  a fixed DNN space. For $\omega\in \bbOmega$, consider a sequence $ \big(u_{\bm{X}_{N} (\omega)   } \big)_N$  of  absolute minimisers of 
\[
	\min  _ {v \in V _{\mathcal{N}} } \E _{ {N } , \omega} (v),
\]
where $ \bm{X} _{N} = (X_1, \dots, X_{N}  )  \, .$
Then, for $\mathbb{P}$-a.e.\ $\omega\in\bbOmega$, the sequence $(u_{\bm{X}_{N} (\omega)   } )_N$  has weak limits along subsequences with respect to $ (\mu _{N, \bm{X}_{N} (\omega  )}, \mu) $, and for any such weak limit, $ u _{\omega} ^  {\mathcal{N}}$, we have
 \begin{equation*}\label {est_final}
  \|u _{\omega} ^  {\mathcal{N} } -f\|_p   \leq  \inf_{\varphi   \in  V _{\mathcal{N}}}  \|\varphi-f\|_p
 \, .
 \end{equation*}
\end{theorem}

The proof relies on the liminf-limsup framework discussed in Section \ref{S:Abstract-liminf-limsup}. To invoke the stability through liminf we observe that  the uniform bound required   follows from the limsup property of global minimisers. 
	 Therefore,  there exists a    $u_\omega  ^  {\mathcal{N} } \in L^p (D, \mu)$ such that $u_{\bm {X}_{N }(\omega)} $ converges weakly {up to a subsequence} to $u_\omega ^  {\mathcal{N} }$ with respect to $ (\mu _{N, \bm{X}_{N} (\omega  )}, \mu) $ and  
 \begin{equation*}\label{w-prob_E_theorem_2}\begin{split} 
\E (u_\omega ^  {\mathcal{N} })  \leq \liminf\limits_{N \rightarrow \infty}  \E _{N, \omega}  [ u_{\bm {X}_{N} (\omega)} ]\, .
\end{split}
 \end{equation*}		
		Then, 
 \begin{equation*}\label{w-conv-u_n_final}
\E \big(u_\omega^  {\mathcal{N}}\big) \leq \liminf\limits_{N\rightarrow \infty}  \E _{{N}, \omega}  \big( u_{\bm {X}_{N} (\omega)} \big)\,\leq \limsup_{N\rightarrow \infty}   \E _{N, \omega}  ( u_{\bm {X}_{N} (\omega)} )   
 \leq \inf_{\varphi\in V_\mathcal{N}} \|\varphi-f\|_p^p  \, .
 \end{equation*}	
 
A natural convergence result follows as the neural network spaces become increasingly expressive, in the sense that they satisfy the natural assumption \eqref{w_ell_2}. Here, $V_\ell$ denotes a suitably chosen sequence of neural network spaces $V_{\mathcal{N}}$.
	
	 \begin{theorem} \label{learningMC_conv}{\sc [convergence of discrete  absolute minimisers]}
For  $f \in L^p (\Omega, \mu) $ 
  we assume that there exists a $w_\ell \in V_\ell$ such that 
  \begin{equation}\label{w_ell_2}\begin{split} 
 \|w_{\ell}-f\|_p \leq \  \beta _\ell \, (f), \qquad 
 \text{and } \ \beta _\ell \, (f) \to 0,  \  \ \ell\to \infty\, .	\end{split}
\end{equation}
Let $\omega\in \bbOmega$ and consider a sequence of    absolute minimisers of 
\begin{equation*}\label{prob_E_def_theorem}
	\min  _ {v \in  V _{\ell} } \E _{ {N(\ell)} , \omega} (v)
\end{equation*}
  denoted by   $ u_{\bm{X}_{N(\ell)} (\omega)   }  ,$ where $ \bm{X} _{N(\ell)} = (X_1, \dots, X_{N(\ell)}  )  \, .$
Then, a.s.\ with respect to $\omega $
the sequence 
$(u_{\bm {X}_{N(\ell)} (\omega)} )_\ell$ converges weakly   to $f$ {, provided $N(\ell)\to\infty$ sufficiently fast, in the sense made precise below} with respect to $ (\mu _{\ell, \bm{X}_{N(\ell)} (\omega  )}, \mu) .$ 
\end{theorem}

The proof follows the  g{framework} of Section~\ref{S:Abstract-liminf-limsup}: the liminf inequality is provided by Lemma~\ref{lim-inf-lemma}, while the limsup inequality requires $\mathcal{E}_{N(\ell),\omega}(w_\ell)\to 0$ almost surely. Since the functions $w_\ell$ change with $\ell$, the strong law of large numbers for a fixed function is not sufficient, and $N(\ell)$ has to grow sufficiently fast; for instance, by Chebyshev's inequality and the Borel--Cantelli lemma, it suffices that $\sum_\ell N(\ell)^{-1}\big\| |w_\ell-f|^p\big\|^2_{L^2(\Omega,\mu)}<\infty$. We refer to \cite{loulakis2023newap} for the precise statement.

 \subsubsection*{Approximate local minimisers} 	
 
 Next, we study the behaviour of \emph{approximate local minimisers} as the number of samples $N$ tends to infinity. Approximate local minimisers may arise naturally from approximate methods for solving the underlying minimisation problem. Typically, in such situations, the value of the corresponding loss is known or can be controlled, but it is not clear a priori how the generalisation error behaves.
{By an \emph{approximate local minimiser} we mean any element of $V_{\mathcal{N}}$ produced by the optimisation procedure; no minimality property is assumed, and the results below use only the values of the loss.} We shall therefore assume certain properties of the loss at these points, without requiring that they correspond to \emph{global minimisers}. The resulting estimates have an \emph{a posteriori} character and provide a clear mathematical justification for the empirical observation that
\begin{center}
\emph{``a consistently small loss, combined with stability, implies a small generalisation error.''}
\end{center}
Here the extra assumption of \emph{stability} appears to be necessary and it is linked to the lim-inf stability Lemma \ref{lim-inf-lemma}.
 
The following result provides a \emph{conditional} estimate, i.e., it is based on in principle \emph{verifiable} conditions on the computed losses.   Its proof is a simple adaptation of the previous analysis, upon noticing that Lemma  \ref{lim-inf-lemma}   still holds, without assuming that the original sequence consists of global minimisers.

\begin{theorem}[\sc conditional convergence of local minimisers]
Let $\omega\in \bbOmega$ and $V _{\mathcal{N}}$  a fixed DNN space. Consider a sequence of  approximate  local minimisers of 
\begin{equation*}\label{prob_E_def_theorem-d2}
	\min  _ {v \in V _{\mathcal{N}} } \E _{ {N } , \omega} (v)
\end{equation*}
denoted by   $ u_{\bm {X}_{N} (\omega), \text{loc}}  ,$ where $ \bm{X} _{N} = (X_1, \dots, X_{N}  )  \, .$
Assume that there exists $M(\omega) $ such that 
\begin{equation*}\label{w-conv-u_n_bound0_lm}\begin{split} 
   \E _{N, \omega}  \big( u_{\bm {X}_{N} (\omega), \text{loc}} \big) \leq M (\omega)<+\infty\, ,
   \qquad \ \text{for }\ \mathbb{P}-\text{a.e. } \omega\in\bbOmega \, .
\end{split}
\end{equation*}	
 Then, for $\mathbb{P}$-a.e.\ $\omega\in\bbOmega$, the sequence $(u_{\bm {X}_{N} (\omega), \text{loc}} )_N$   has weak limits along subsequences with respect to $ (\mu _{N, \bm{X}_{N} (\omega  )}, \mu) $, and for any such weak limit, $u_{\omega, \text{loc}}^\mathcal{N},$  the following estimate holds
\begin{equation*}\label {est_final_lm}
  \|u_{\omega, \text{loc}}^\mathcal{N} -f\|_p {^p}   \leq  \liminf\limits_{N\rightarrow \infty}  \E _{N, \omega}  \big( u_{\bm {X}_{N}(\omega), \text{loc}} \big)\, , \quad  \text{for }\ \mathbb{P}-\text{a.e. } \omega\in\bbOmega
 \, .
\end{equation*}
 Furthermore, if the loss corresponding to $u_{\bm {X}_{N}(\omega), \text{loc}} $ has a limit, 
 \begin{equation*}\label {est_final_lm_2}
 \lim \limits_{N\rightarrow \infty}  \E _{N, \omega}  \big( u_{\bm {X}_{N}(\omega), \text{loc}} \big)= \beta _{\omega, \text{loc}}^\mathcal{N} \, ,
\end{equation*}
then
 \begin{equation*}\label {est_final_lm_3}
  \|u_{\omega, \text{loc}}^\mathcal{N} -f\|_p   \leq  \big(\beta _{\omega, \text{loc}}^\mathcal{N}\big)^{1/p} \quad  \text{for }\ \mathbb{P}-\text{a.e. } \omega\in\bbOmega
 \, .
\end{equation*}
\end{theorem}

The above result provides conditions, based solely on the loss values at these approximate local minima, under which the error with respect to the target function remains small. One such condition is the consistent observation of small loss values as the number of samples increases. This provides yet another indication of the robustness of the algorithm.

\newpage
\subsubsection{Analysis for quadrature based deterministic training}\label{SubS:Anal:functionQh}
 
 
 The ideas developed in the previous section for probabilistic training can also be applied to quadrature-based deterministic training. The key observation is that a deterministic quadrature rule can likewise be associated with an empirical measure. Consequently, the main steps of the preceding analysis can be adapted to this setting as well.
 
 Recall that we assume   that  $f$ is known  at a finite set of points 
 $z\in K_h,$ where $K_h$ is a set of points linked to a quadrature rule.  
 The corresponding  quadrature is based on this 
  set $K_h$ of discrete points $z\in K_h$  and corresponding nonnegative weights $w_z$ 
 \begin{equation*}\label{quadr_lf:S=}
\sum _{z\in K_h} \, w_z \, G(z) \cong \int _{\Omega} \, G(x) \, \d x .
\end{equation*}
The corresponding loss is 
\begin{equation*}
\label{E_h-d5}
\mathcal{E}_{Q, h}( g )  =     \sum _{z\in K_h} \, w_z \, \Big  | g(z) - f(z) \Big  | ^p  \, . 
\end{equation*}


Consider
the empirical measures on $\Omega,$
\begin{equation*}\label{emp_meas-d3}
	\mu  _ {Q, h}  =    \sum _{z\in K_h} \, w_z \, \   \delta _{z} \, ,
\end{equation*}
then
\begin{equation*}\label{energy_emp_meas-d2}
	\mathcal{E}_ {Q, h} ( g )  = \int _\Omega |g (x)-f(x)|^p\, \d \mu _ {Q, h} (x) \, .
\end{equation*}

%

{We use the notions of weak convergence of measures and of functions with respect to $(\nu_n,\nu)$ introduced in Section~\ref{SubS:weakconv}.}

{The role of the almost sure weak convergence of the random empirical measures is now played by an assumption on the quadrature rule, namely that it is convergent for continuous functions, so that the deterministic empirical measures weakly converge to $\mu$:}

\begin{equation*}\label{perklln-d2}
 \mu _ {Q, h} \to \mu, \qquad \text{ as }h \to 0.
 \end{equation*}
{That is, for any bounded and continuous function $g$ on $\Omega$ there holds}

\begin{equation*} \begin{split}
\int_\Omega g(x)\, \d \mu  _ {Q, h} (x)= \sum _{z\in K_h} \, w_z \,   g(z)    \,   \longrightarrow  \int_\Omega g(x)\, \d \mu(x)\, , \qquad \text{ as }h \to 0. \end{split}
\end{equation*}

\noindent
 \emph{The $\liminf$ inequality.\/}
	A variation of the proof of Lemma \ref{lim-inf-lemma} yields,

	\begin{lemma}[\sc the $\liminf$ inequality]
\label{lim-inf-lemmaQh}
Let $p>1$ and $V _{\mathcal{N}}$  a fixed DNN space.  Consider a sequence of  absolute minimisers of 
 \begin{equation*}\label{prob_E_def_theorem0-d2}
	\min  _ {v \in V _{\mathcal{N}} } \E _ {Q, h}  (v)
 \end{equation*}
denoted by   $ u_ {Q, h}      .$  
Assume that  
 \begin{equation}\label{w-conv-u_n_bound0Qh}\begin{split} 
   \E _ {Q, h} \big(  u_ {Q, h} \big) \leq M  <+\infty\, ,
   \qquad \text{ for all  } h>0.
\end{split}
 \end{equation}	
Then,   the sequence $(  u_ {Q, h})_h$  has weak limits along subsequences with respect to $ (\mu   _ {Q, h}, \mu) $, and for any such weak limit, $ u  ^  {\mathcal{N}}$, we have
 \begin{equation*}\label{w-conv-u_n_0-d2}\begin{split} 
\E \big(u ^  {\mathcal{N}}\big) \leq \liminf\limits_{h \rightarrow 0}  \E   _ {Q, h}  \big(  u_ {Q, h} \big)\, .
\end{split}
 \end{equation*}	
\end{lemma}

 \noindent
 \emph{Estimates as $h\to 0.$} 
As before we consider  the $h\to 0$ limit. Still in the deterministic case, a consequence of the liminf property is that these limits, for both global minimisers and approximate local minimisers,  are  robust.

We have the following estimate:   
 \begin{equation}\label{limsupboundQh}
\limsup_{h \rightarrow 0}\E _ {Q, h}  \big(  u_ {Q, h} \big) \leq  \inf_{\varphi   \in  V _{\mathcal{N}} } \|\varphi-f\|_p^p\, .
 \end{equation}
As before,  the minimisers $u _ {Q, h} ,$ satisfy 
 \begin{equation*}\label{prob_E_theorem_2_10Qh}
	\E _ {Q, h} (u _ {Q, h}  ) \leq \E _ {Q, h} (v)\, , \quad \text {for all } \ v \in  V _{\mathcal{N}}  \, . 
 \end{equation*}
Then, consider a sequence $(v_n)_n$ in $V _{\mathcal{N}}$ such that,
\[
\lim_{n\to\infty} \|v_n-f\|_p=\inf_{\varphi   \in  V _{\mathcal{N}} } \|\varphi-f\|_p.
\]
Notice that the sequence $(v_n)_n$ depends solely on $V _{\mathcal{N}}$ and it is 
fixed as $h\to 0\, .$ Then, under some regularity hypotheses on $f$ and $V _{\mathcal{N}},$ for each $n$ we have  
%
 \begin{equation*}\label{llnvn}
\lim_{h\to 0} \E_ {Q, h} (v_n) =\|v_n-f\|_p^p.
 \end{equation*}
Therefore, 
 \begin{equation*}\label{w-prob_E_theorem_3Qh}
\limsup_{h \rightarrow 0}\E _ {Q, h}  \big(  u_ {Q, h} \big)
\leq \limsup_{h \rightarrow 0}\E _ {Q, h}  \big(  v_n \big)
 =  \|v_n-f\|_p^p   \, ,
 \end{equation*}	
and	the estimate \eqref{limsupboundQh} follows.
%


 
	\begin{theorem}[\sc limit behaviour of minimisers as $h\to 0$] 
Let $p>1$ and $V _{\mathcal{N}}$  a fixed DNN space. Consider a sequence $ \big(u_ {Q, h} \big)_h$  of  absolute minimisers of 
\[
	\min  _ {v \in V _{\mathcal{N}} } \E  _ {Q, h} (v)\, .
\]
Then the sequence $ \big(u_ {Q, h} \big)_h$ has weak limits along subsequences with respect to $ (\mu _ {Q, h}, \mu) $ as $h\to 0$ and for any such weak limit, $ u  ^  {\mathcal{N}}$, we have
 \begin{equation*}\label {est_final-d2}
  \|u  ^  {\mathcal{N} } -f\|_p   \leq  \inf_{\varphi   \in  V _{\mathcal{N}}}  \|\varphi-f\|_p
 \, .
 \end{equation*}
\end{theorem}

As before, the proof relies on the liminf-limsup framework discussed in Section \ref{S:Abstract-liminf-limsup}, and the uniform bound 
required for the liminf part is  a consequence of  the limsup property of global minimisers. A similar result to Theorem  \ref{learningMC_conv} holds in this case as well. \\

\noindent
 \emph{Approximate local minimisers.\/} 	
 Similar  \emph{conditional} estimates for approximate local minimisers still hold yielding a posteriori control, as $h\to 0,$ of the generalisation error by the
 computed losses. In fact if we  
consider a sequence of  approximate  local minimisers of our problem denoted by   $ u_{Q, h, \text{loc}}  ,$  and assuming that 
  there exists $M  $ such that 
\begin{equation*}\label{w-conv-u_n_bound0_lm-d2}\begin{split} 
   \E _{Q, h}  \big( u_{Q, h, \text{loc}}  \big) \leq M  <+\infty\, ,
   \qquad \ \text{for }\ h \to 0 .
\end{split}
\end{equation*}	
 Then,  the sequence $ u_{Q, h, \text{loc}}  $   has weak limits along subsequences with respect to $ (\mu _{Q, h}, \mu) $, and for any such weak limit, $u_{ \text{loc}}^\mathcal{N},$  the following estimate holds
\begin{equation*}\label {est_final_lm-d2}
  \|u_{  \text{loc}}^\mathcal{N} -f\|_p {^p}   \leq  \liminf\limits_{h\rightarrow 0}  \E _{Q, h}  \big(u_{Q, h, \text{loc}}  \big)\, . 
\end{equation*}
 Furthermore, if the loss corresponding to $u_{Q, h, \text{loc}}$ has a limit, 
 then this limit controls the generalisation error of the limiting points $u_{  \text{loc}}^\mathcal{N}.$ 

\noindent
\subsubsection{Operator learning}
The analysis of this section extends to the operator learning problem of
Section~\ref{SubS:ModP:function-learning}, with appropriate amendments. Let
$X$ be a separable metric space of inputs with a probability measure
$\mu\in\mathcal{P}_p(X)$, let $Y$ be a separable Hilbert space, and let
$\mathcal{G}:X\to Y$ be measurable with
$\int_X\|\mathcal{G}(a)\|_Y^p\,\d\mu<\infty$. The energies are
$$\mathcal{E}(G)=\int_X\|G(a)-\mathcal{G}(a)\|_Y^p\,\d\mu(a)$$
and, for i.i.d.\ samples $a_i=X_i(\omega)$ of $\mu$,
$$\mathcal{E}_{N,\omega}(G)=\frac1N\sum_{i=1}^N\|G(X_i(\omega))-\mathcal{G}(X_i(\omega))\|_Y^p,$$
minimised over a class $V_{\mathcal{N}}$ of neural operators. To verify the
liminf inequality, we first note that the empirical measures
$\mu_{N,\bm{X}_N(\omega)}$ converge weakly to $\mu$ for $\mathbb{P}$-a.e.\
$\omega$, by Varadarajan's theorem \cite{V-58}, which holds on any separable
metric space. Furthermore, an analogue of Lemma~\ref{lim-inf-lemma} still
holds, but with a series of appropriate modifications. We shall briefly
discuss them next, without aiming at precise statements in these notes. We
start as in the proof of Lemma~\ref{lim-inf-lemma} by considering the
measures $(\bm i\times v_N)_\sharp\mu_{N,\bm{X}_N}$ on $X\times Y$, with
$v_N=G_{\theta_N}-\mathcal{G}$ taking values in $Y$. A key issue is that a
uniform bound on $\int\|v_N\|_Y^p\,\d\mu_{N,\bm{X}_N}$ no longer yields
tightness in the norm topology of $Y$, whose bounded sets are not compact.
This is precisely the situation treated systematically in \cite[Chapter 5]{AGS:2005} for
measures on separable Hilbert spaces. A new norm $\| \cdot \| _\varpi$ is introduced in  $Y$ which induces    a weaker topology with the property  that  bounded sets of $Y$  are relatively compact  and metrisable.  Then narrow convergence of the plans on $X\times Y$ with respect to this
topology is used, together with cylindrical test functions
$\zeta(a)\,\<y,e_j\>_Y$ built from an orthonormal basis $(e_j)$ of $Y$ in
place of the test functions $\zeta(x)\,y$ of the scalar case, and the
barycentre $\overline{\bm\gamma}_1(a)=\int_Y y\,\d\bm\gamma_a(y)$ is well
defined.
The liminf inequality then follows from the convexity and weak lower
semicontinuity of $\|\cdot\|_Y^p$ through Jensen's inequality, exactly as in
the last step of the proof of Lemma~\ref{lim-inf-lemma}, see \cite[Theorem 5.4.4]{AGS:2005} for details. With these
modifications, the estimates for global minimisers and the conditional
estimates for approximate local minimisers hold for operator learning as
stated, with the norms of $L^p(X,\mu;Y)$ in place of those of
$L^p(\Omega,\mu)$. We do not carry out the details here; the point we wish to
make is that the abstract approach discussed in this section may be applied
to operator learning as well, and that the infinite-dimensional character of
$Y$ requires a careful analysis based on \cite[Chapter 5]{AGS:2005}. 
\subsection{Generating measures: the static Wasserstein case}\label{SubS:Conv:measures}

We now follow the approach of Section~\ref{S:Abstract-liminf-limsup} for the
simplest of the generative models of Section~\ref{SubS:Ex:measures}: a
static neural measure space, and $d_\mathscr{M}(\cdot, \cdot )=W_p (\cdot, \cdot ),$
the Wasserstein distance. Let $p\ge 1$, $V=\mathcal{P}_p(\mathbb{R}^k)$ with the
metric $W_p$, and consider the neural measure space
\begin{equation*}
V_\Theta=\big\{(T_\theta)_\sharp\gamma_0:\ \theta\in\Theta\big\}\subset V,
\end{equation*}
where $\gamma_0\in\mathcal{P}_p(\mathbb{R}^m)$ is the reference measure and
the networks $T_\theta:\mathbb{R}^m\to\mathbb{R}^k$, $T_\theta \in V_{\mathcal{N}}$,
are Lipschitz continuous, so that $(T_\theta)_\sharp\gamma_0$ has finite
$p$-th moment. The continuous problem is $\min_{\mu\in V}W_p(\mu,\nu_{\, \text{true}})$,
whose minimiser is $\nu_{\, \text{true}}$. The \emph{semidiscrete} problem uses the same loss, but the minimisation is
restricted to $V_\Theta$. Further steps of discretisation are needed to arrive at the fully discrete
loss used in practice, whose analysis is left for future work. Here we
discuss three intermediate discretisations, each adding one level of
approximation to the semidiscrete problem, in order to show how such a step
by step treatment clarifies the behaviour of the method at each level and
where the ingredients of the framework of
Section~\ref{S:Abstract-liminf-limsup} enter. First consider the
problem with data, in which $\nu_{\, \text{data}}$ replaces
$\nu_{\, \text{true}}$ in the minimisation:
\begin{equation}\label{mm_nn:abstract:Wp}
\min_{\mu\in V_\Theta} W_p(\mu,\nu_{\, \text{data}}),\qquad
\nu_{\, \text{data}}=\frac1N\sum_{i=1}^N\delta_{X_i},\quad X_i\sim\nu_{\, \text{true}}\ \text{i.i.d.},
\end{equation}
in which the distance is assumed to be computed exactly; its approximation
is discussed in the second and third steps. 

\subsubsection*{The analogue of exact integration}
Problem \eqref{mm_nn:abstract:Wp} is the counterpart, for measures, of the
function learning problem in which $f$ is known only through its point values
but the integrals are evaluated exactly, as in \eqref{tildeEQh}: the data
$\nu_{\, \text{data}}$ play the role of the point values, and the exactly
computed distance that of the integral. As in that case, the estimates are
elementary, because $W_p$ is a metric on $V$ and $\nu_{\, \text{data}}$ is a
fixed point of $V$. Indeed, for any minimiser $\mu_N$ of
\eqref{mm_nn:abstract:Wp}, the triangle inequality implies
\begin{equation}\label{Wp:Cea}
W_p(\mu_N,\nu_{\, \text{true}})\le \inf_{\mu\in V_\Theta}W_p(\mu,\nu_{\, \text{true}})
+2\,W_p(\nu_{\, \text{data}},\nu_{\, \text{true}})\, .
\end{equation}
Estimate \eqref{Wp:Cea} is a C\'ea-type estimate.    Since $W_p(\nu_{\, \text{data}},\nu_{\, \text{true}})\to 0$ almost surely as
   $N\to\infty$ (the empirical measures converge weakly and their $p$-th
   moments converge, by the strong law of large numbers, see
   \cite[Proposition 7.1.5]{AGS:2005}), if the spaces $V_{\Theta_\ell}$ are
   expressive enough, in the sense that for every $\mu\in V$ there exist
   $\mu_\ell\in V_{\Theta_\ell}$ with $W_p(\mu_\ell,\mu)\to 0$, then
   \eqref{Wp:Cea} shows that minimisers of \eqref{mm_nn:abstract:Wp} over
   $V_{\Theta_\ell}$, whenever they exist, converge almost surely to
   $\nu_{\, \text{true}}$ in $W_p$, for any $N(\ell)\to\infty$.

In addition, one may recover an a posteriori type estimate by observing
that for an approximate local minimiser $\mu^*_\theta\in V_\Theta$ (and in
fact for \emph{any} $\mu_\theta\in V_\Theta$) we have trivially
\begin{equation}\label{Wp:apost}
W_p(\mu^*_\theta,\nu_{\, \text{true}})\le W_p(\mu^*_\theta,\nu_{\, \text{data}})
+W_p(\nu_{\, \text{data}},\nu_{\, \text{true}})\, .
\end{equation}
The first term on the right-hand side is the computed loss, and the second
converges to zero almost surely as $N\to\infty$. Hence, if a sequence
$(\mu^*_N)_N$ of approximate local minimisers has consistently small losses,
its elements are eventually close to $\nu_{\, \text{true}}$, and any
accumulation point $\bar\mu$ of the sequence, with respect to weak
convergence, satisfies, by the lower semicontinuity of $W_p$,
\begin{equation}\label{Wp:apost:lim}
W_p(\bar\mu,\nu_{\, \text{true}})\le\liminf_{N\to\infty}W_p(\mu^*_N,\nu_{\, \text{data}})\, ,
\end{equation}
which is the counterpart of the conditional estimate for approximate local
minimisers of Section~\ref{SubS:Conv:MCInterpolation}. The liminf--limsup framework itself is needed only for the existence and
the identification of accumulation points, and its ingredients are 
verified as follows: a bound $W_p(\mu,\nu_{\, \text{data}})\le C$ bounds the $p$-th
moment of $\mu$, hence gives tightness \cite[Remark 5.1.5]{AGS:2005} and
weak compactness \cite[Theorem 5.1.3]{AGS:2005}, and $W_p$ is lower
semicontinuous with respect to weak convergence
\cite[Proposition 7.1.3]{AGS:2005}, which is the liminf inequality.

\subsubsection*{Training: sampling the generator measure}
Training begins, in this setting, when the distance is no longer computed
exactly. The first instance is the replacement of $\mu_\theta$ by the
empirical measure of $M$ generated samples,
$\widehat\mu_{\theta,M}=(T_\theta)_\sharp\widehat\gamma_M$, where
$\widehat\gamma_M$ is the empirical measure of $M$ i.i.d.\ samples of
$\gamma_0$. The loss $W_p(\widehat\mu_{\theta,M},\nu_{\, \text{data}})$ is a
distance between two atomic measures and is computable. The triangle
inequality now gives, for any $\theta$,
\begin{equation}\label{Wp:apost:sampled}
W_p(\mu_\theta,\nu_{\, \text{true}})\le W_p(\widehat\mu_{\theta,M},\nu_{\, \text{data}})
+W_p(\widehat\mu_{\theta,M},\mu_\theta)+W_p(\nu_{\, \text{data}},\nu_{\, \text{true}}),
\end{equation}
and the middle term involves the computed $\theta$, which depends on the
generator samples: this is the key difficulty of Section~\ref{SubS:ModP:function-learning} in its generative form. 
  As we shall see in Section~\ref{S:reliability}, this term can be bounded as
   \begin{equation}\label{Wp:pushforward}
   W_p(\widehat\mu_{\theta,M},\mu_\theta)\le \mathrm{Lip}(T_\theta)\,W_p(\widehat\gamma_M,\gamma_0),
   \end{equation}
   where $\mathrm{Lip}(T_\theta)$ is the Lipschitz constant of the computed
   network; the corresponding rates are also given there.

\subsubsection*{Training: approximating the distance through entropic regularisation}
The second instance of training is the approximation of the distance itself. We shall not pursue further discretisations in these notes but we shall briefly examine an important case, both in practice and for our purposes, namely the entropic
regularisation. 
The regularised problem defined next is strictly convex and its
solution is computed by Sinkhorn's matrix scaling algorithm
\cite{Cuturi:2013}, which is what makes $\mathcal{T}_\delta$ an appropriate model for fully discrete losses 
 in  Wasserstein based generative modelling \cite{Genevay_Peyre_Cuturi:2018}. 
For $\delta>0$ let
\begin{equation}\label{T_delta}
\mathcal{T}_\delta(\mu,\nu)=\min_{\pi\in A(\mu,\nu)}\Big\{\int_{\R^k\times\R^k}|x-y|^2\,\d\pi(x,y)
+\delta\,\mathrm{KL}\big(\pi\,\|\,\mu\otimes\nu\big)\Big\},
\end{equation}
where, for probability measures $\pi,\sigma$ on $\R^k\times\R^k$, the
relative entropy is $\mathrm{KL}(\pi\,\|\,\sigma)=\int\log\frac{\d\pi}{\d\sigma}\,\d\pi$
if $\pi$ is absolutely continuous with respect to $\sigma$ and
$\mathrm{KL}(\pi\,\|\,\sigma)=+\infty$ otherwise, so that
$\mathcal{T}_0=W_2^2$. Note that
$\mathcal{T}_\delta$ is not a metric: $\mathcal{T}_\delta(\mu,\mu)>0$ and the
triangle inequality fails, so that the direct arguments of the exact case without further training 
are not available. The trained model problem becomes
$$\min_{\mu\in V_\Theta}\mathcal{T}_\delta(\mu,\nu_{\, \text{data}}),$$
 and we would like to study 
  what happens as $\delta\to 0$. This is a question of discrete
variational stability for the family of energies
$\mathcal{T}_\delta(\cdot,\nu_{\, \text{data}})$, and the   key analytical ingredients   are
available in the literature.
\begin{itemize}
\item \emph{Coercivity and lower semicontinuity.} Since the entropy term is nonnegative
and the admissible plans are the same, $\mathcal{T}_\delta\ge W_2^2$ for every
$\delta>0$. A bound $\mathcal{T}_\delta(\mu,\nu_{\, \text{data}})\le C$
therefore bounds $W_2(\mu,\nu_{\, \text{data}})$, hence the second moment of
$\mu$, uniformly in $\delta$, and yields tightness.
Furthermore, $\mathcal{T}_\delta$ is lower
semicontinuous with respect to weak convergence, as the cost is lower
semicontinuous and the relative entropy is jointly weakly lower
semicontinuous. Combined with $\mathcal{T}_\delta\ge W_2^2$ and the lower
semicontinuity of $W_2$, this gives the liminf inequality: if
$\mu_\delta\rightharpoonup\mu$ as $\delta\to0$ then
$\liminf_{\delta\to0}\mathcal{T}_\delta(\mu_\delta,\nu)\ge W_2^2(\mu,\nu)$. 
\item \emph{Limsup.} For fixed $\mu,\nu$ with finite second moments,
$\mathcal{T}_\delta(\mu,\nu)\to W_2^2(\mu,\nu)$ as $\delta\to0$, and for
compactly supported measures $\mathcal{T}_\delta-W_2^2=O(\delta\log(1/\delta))$.
   For the $\delta$-limit the recovery sequence is thus trivial, the constant
   sequence $\mu$; for the limit over increasingly expressive spaces
   $V_{\Theta_\ell}$ can be  provided by the approximation property of the
   neural measure spaces.

 \end{itemize}
   The main results required to establish the above claims rigorously can
   be traced back to the $\Gamma$-convergence of entropic optimal transport
   to the Monge--Kantorovich problem, due to L\'eonard \cite{Leonard:2012}
   in the setting of the Schr\"odinger problem, and established in the
   present form in \cite{Carlier_Duval_Peyre_Schmitzer:2017}. 

Finally, note that a posteriori type guarantees for accumulation points of
minimisers and approximate local minimisers can be derived from
\begin{equation}\label{Wp:apost:entropic}
W_2(\mu_\theta,\nu_{\, \text{true}})\le \mathcal{T}_\delta(\mu_\theta,\nu_{\, \text{data}})^{1/2}
+W_2(\nu_{\, \text{data}},\nu_{\, \text{true}})\, ,
\end{equation}
for any computed $\mu_\theta\in V_\Theta.$ This is just a straightforward
consequence of $\mathcal{T}_\delta\ge W_2^2$, and implies that the computed
regularised loss controls the error, for every $\delta$. The convergence of
the minimisers over $V_\Theta$ as $\delta\to0$ and $N\to\infty$ involves the
same complications as in the exact case, now for minimisers depending on
both parameters, and is not pursued here. We also note that
$\mathcal{T}_\delta(\mu,\mu)>0$, so that the minimiser of
$\mathcal{T}_\delta(\cdot,\nu_{\, \text{data}})$ is not $\nu_{\, \text{data}}$
but an entropically smoothed version of it: the regularisation can be
viewed as a selection principle, compare to Section~\ref{SubS:hallucinations}.

 \subsection{PINNs for parabolic PDEs
}\label{SubS:PINN-Parabolic}
\smallskip

  Recall the notation used in the parabolic case:  $ \Omega \subset \mathbb{R}^d $, is an open and  bounded and   $ \Omega_T = \Omega \times (0,T] $ denotes the space-time domain for some fixed time $ T>0. $
This section is devoted to problems of the form 
\begin{align}\label{ParabolicPDE:an}
\begin{cases}
u_t + Lu = f ,\;\;\; \; \textrm{in} \;\: \Omega_T, \\
u= 0, \;\;\;\;\;\;\;\;\; \; \textrm{on} \;\: \partial \Omega \times (0,T],  \\
u = u^0, \;\;\;\;\;\;\;\;\, \textrm{on} \;\: \Omega \times \lbrace t=0 \rbrace \, .
\end{cases}
\end{align}
where $L$ is an elliptic operator defined in Section \ref{SubS:Ex:EllipticPDE}. In  this section we discuss  convergence properties of approximations of \eqref{ParabolicPDE:an} obtained by minimisation   time-discrete energy functionals over appropriate sets of  neural network functions. 
For simplicity in the exposition we shall assume that $\Omega$ is a convex Lipschitz domain; the general case of a non-convex domain can be treated as well. 

Consider the space-time energy functional $ \mathcal{G} : H^1(0,T ; L^2 (\Omega)) \cap L^2(0,T ; H^2(\Omega)  \cap  H^1_0(\Omega))\rightarrow \overline{\mathbb{R}} $ as follows
\begin{equation}\label{ParabFunctional:an}
\begin{gathered}
\mathcal{G}(v) =
\int_0^T \| v_t(t) + L v(t) - f(t) \|^2_{L^2(\Omega)} \d t +  | v(0) - u^0  |^2_{H^1(\Omega)}\, .
\end{gathered}
\end{equation}
Following the notation and the assumptions on the problem formulation  of Section \ref{SubS:Ex:parabolic}, we consider the neural network approximating spaces 
$
V _{\mathcal{N}, 0} \subset H^1(0,T ; L^2 (\Omega)) \cap L^2(0,T ; H^2(\Omega) \cap H^1_0(\Omega))$ 
are still denoted by $V _{\mathcal{N}}$ and $V _{\mathcal{N}, 0}$ (the $0$ in the latter indicates that the boundary conditions are enforced strongly). 
Residual-based methods are     based on the minimisation of the functional
\begin{equation*}\label{mm_nn:abstract-d5}
	\min  _ {v \in  V _{\mathcal{N}, 0} } \mathcal{G} (v) \, .  
\end{equation*}
In this section we consider only methods  based on \emph{time-discrete training.} 
In particular we consider 
\begin{equation}\label{Functional_k:an}
\mathcal{G}_{k, IE} (v) = \sum _{n=1}^N \,  k_n \, \int_{\Omega }\big | \frac{v ^n-  v^{n-1}}{k_n} + L v ^n  - f^n \big |^2 \, \, \d x +\, 
|v g{^0}  - u^0 |_{H^1 (\Omega)}^2 
\end{equation}
and  the corresponding minimisation problem over $  V _{\mathcal{N}, 0}$
%
\begin{equation}\label{ieE-minimize_k:an}
	\min  _ {v \in  V _{\mathcal{N}, 0} } \mathcal{G} _{k, IE}(v)\, ,
\end{equation}
 %
%
We shall compare the above method to the corresponding explicit-Euler-type discretisation 
\begin{equation}\label{Functional_k_Ex_intro:an}
\mathcal{G}_{k, EE} (v) = \sum _{n=1}^N \,  k_n \, \int_{\Omega }\big | \frac{v ^n-  v^{n-1}}{k_n} + L v ^{n-1}  - f^{n-1} \big |^2 \, \, \d x +\, |v g{^0}  - u^0 |_{H^1 (\Omega)}^2\, .
\end{equation}
{which, from the point of view of approximation, is equally natural and equally consistent as $\mathcal{G}_{k, IE}$; notice that, as a minimisation problem over space-time neural network functions, it is just as implicit as $\mathcal{G}_{k, IE}$.}
As mentioned earlier, the comparison of these two functionals and their minimisers demonstrates that the stability of the algorithms is intrinsically linked to the coercivity of the loss functional. It also highlights the potentially decisive role of the training process in ensuring coercivity and, consequently, the stability of the algorithm. 

\subsubsection*{Discrete variational stability}

{Before studying the time-discrete case let us mention briefly that the variational stability of the semidiscrete problem, in which $\mathcal{G}$ is minimised over the spaces $V_{\ell,0}$ of Section~\ref{SubS:Ex:parabolic}, hinges on the maximal regularity of the parabolic problem, Theorem 5, p.~382 in \cite{Evans}, and on the Aubin--Lions theorem (an analogue of the Rellich--Kondrachov theorem in the parabolic case), see for example \cite{SZ}.}
\begin{equation}\label{G Equicoercivity}
\begin{gathered}
\textrm{If} \;\; \mathcal{G}(u) \leq C \;\; \textrm{for some} \;\: C >0 \;,\; \textrm{we have} \\
\|u\|_{L^2(0,T; H^2(\Omega))} + \|u_t\|_{L^2(0,T; L^2(\Omega))} \leq C_1 \, .
\end{gathered}
\end{equation}
{This is the maximal regularity estimate for the parabolic problem; the initial term of $\mathcal{G}$ is measured in the $H^1$ seminorm precisely because this estimate requires $u(0)\in H^1_0(\Omega)$.}
 This coercivity bound implies the   liminf inequality:  For a sequence,   $ u_\ell $, such that $ \mathcal{G}_\ell (u_\ell) \leq C $ uniformly in $ \ell $,  \eqref{G Equicoercivity} implies that $(u' :=u_t): $ $ ||u_\ell ||_{L^2(0,T; H^2(\Omega))} + ||u'_\ell||_{L^2(0,T; L^2(\Omega))} $ are uniformly bounded.
%
This leads to $ u_\ell' +  L u_\ell - f \rightharpoonup u' + L u - f \,, $ and the convexity of $ \int_{\Omega_T} | u_\ell' + L u_\ell -f |^2 $ implies weak lower semicontinuity, that is
\begin{equation}\label{proofGlowersemicontinuityeq2}
\int_{\Omega_T} | u' + L u -f |^2  \leq \liminf_{\ell \rightarrow \infty}  \int_{\Omega_T} | u_\ell' + L u_\ell -f |^2
\end{equation}
and therefore we conclude that $ \mathcal{G}(u) \leq \liminf_{\ell \rightarrow \infty} \mathcal{G}_\ell (u_\ell). \\ $

We turn now our attention to the functional ${G}_{k, IE}.$
 Next we utilise   the \emph{time reconstruction} introduced earlier in Section \ref{S:apriori_apost} 
   $\widehat U$ 
    to be the piecewise linear  interpolant    
of $U$ defined as: 
\begin{equation*}\label{def_hatU-d2}
\widehat U(t):=\ell_0^n(t) U^{n-1}+ \ell_1^n(t) U^{n}, \quad t\in I_n,
\end{equation*}
with
$
\ell_0^n(t) := (t^n-t)/k_n$  and  $\ell_1^n(t) := (t-t^{n-1})/k_n.$  The piecewise constant interpolant of $U^{j}$  is denoted by $\overline U, $  
\begin{equation*}\label{def_overU-d3}
\overline U(t):= U^{n}, \quad t\in I_n\, .\end{equation*}
Now, the discrete energy  $ \mathcal{G}_{IE, k}   $ can be written as follows
\begin{equation*}\label{ParabFunctional_k}
\begin{split}
	\mathcal{G}_{IE, k} (U) = &\| \widehat U_t + L \overline U - \overline f \|^2_{L^2(0,T; L^2(\Omega))}  +\,     {|  U ^0  - u^0 |_{H^1(\Omega)}^2}\\
=& \int_0^T \|  \widehat U_t + L \overline U - \overline f \|^2_{L^2(\Omega)}\, \d t  +\,     {|  U^0  - u^0 |_{H^1(\Omega)}^2}\, .
\end{split}
\end{equation*}
 We show next that the functional $ \mathcal{G}_{IE, k}$   is equicoercive with respect to $ \widehat U , \overline U $. That is,
\begin{equation}\label{GEquicoercivity_TD}
\begin{gathered}
\textrm{If} \;\; \mathcal{G}_k(U) \leq C \;\; \textrm{for some} \;\: C >0 \;,\; \textrm{we have} \\
\|\overline U\|_{L^2(0,T; H^2(\Omega))} + \|\widehat U'\|_{L^2(0,T; L^2(\Omega))} \leq C_1
\end{gathered}
\end{equation}
To see why \eqref{GEquicoercivity_TD} is true we first observe that $\mathcal{G}_k(U) \leq C $ and the fact that $ f \in L^2(\Omega_T) $, lead  to the uniform bound
\begin{equation}\label{Proofof G EquicoercEq4}
\| \widehat U_t + L \overline U \|_{L^2(0,T; L^2(\Omega))} \leq C_1\, . 
\end{equation}
The coercivity bound \eqref{GEquicoercivity_TD} hinges on \emph {discrete maximal regularity} bound, i.e., 
 that $\|\overline U\|_{L^2(0,T; H^2(\Omega))} + \|\widehat U'\|_{L^2(0,T; L^2(\Omega))} $ can be controlled by $\| \widehat U_t + L \overline U \|_{L^2(0,T; L^2(\Omega))} \, .$
This follows by 
\begin{equation}\label{dMR_TD}
	\begin{split}
  \|\widehat U_t + L \overline U     \|^2_{L^2(0,T; L^2(\Omega))} 
& =  \|\widehat U_t     \|^2_{L^2(0,T; L^2(\Omega))}  
+  \| L \overline U     \|^2_{L^2(0,T; L^2(\Omega))} \\
&+2 \sum _{n=1}^N \, \int_{I_n}  \, 
  \<   \widehat U_t , L \overline U     \> \, \ d t 
\\
& =  \|\widehat U_t     \|^2_{L^2(0,T; L^2(\Omega))}  
+  \| L \overline U     \|^2_{L^2(0,T; L^2(\Omega))} \\
& \qquad +2 \sum _{n=1}^N \, \int_{I_n}  \,   \<    \big [ \frac{U ^n-  U^{n-1}}{k_n}  \big ], L   U ^n          \> \, \ d t \\
& =  \|\widehat U_t     \|^2_{L^2(0,T; L^2(\Omega))}  
+  \| L \overline U     \|^2_{L^2(0,T; L^2(\Omega))} \\
& \qquad +2 \sum _{n=1}^N \,    \,    \<    \big [  {U ^n-  U^{n-1}}   \big ], L   U ^n          \> \,   \\
& =  \|\widehat U_t     \|^2_{L^2(0,T; L^2(\Omega))}  
+  \| L \overline U     \|^2_{L^2(0,T; L^2(\Omega))}  
+ \<   L\,  {U ^N }   ,  U ^N\>        \\
& \qquad + \sum _{n=1}^N \,    \,   \<  L\,   \big [  {U ^n-  U^{n-1}}   \big ],   {U ^n-  U^{n-1}  }       \> \,  -  \<   L\,  {U ^0 }   ,  U ^0\> \, .
 \end{split}
\end{equation}
{All terms on the right-hand side of \eqref{dMR_TD} are nonnegative except the last one, and $\<LU^0,U^0\>\le C|U^0|^2_{H^1(\Omega)}$ is controlled by the initial term of the energy; this is where the $H^1$ seminorm in the initial condition is used. Hence $\|\widehat U_t\|^2+\|L\overline U\|^2\le \|\widehat U_t+L\overline U\|^2+C$, and the $H^2$ bound on $\overline U$ follows from the bound on $L\overline U$ by elliptic regularity, since $\Omega$ is convex. This implies \eqref{GEquicoercivity_TD}.}
To verify the 
 {$\liminf$ inequality}, 
  we assume there is a sequence, still denoted by $ U_\ell $, such that $ \mathcal{G}_{IE, \ell} (U_\ell) \leq C $ uniformly in $ \ell .$
From the discrete stability estimate, the uniform bound $ \mathcal{G}_{IE, \ell} (U_\ell) \leq C $ implies that $\|\overline U_\ell\|_{L^2(0,T; H^2(\Omega))} + \|\widehat U_\ell'\|_{L^2(0,T; L^2(\Omega))} \leq C_1\, ,$
 are uniformly bounded. By the {weak} relative compactness in $L^2(0,T; L^2(\Omega))$ we have (up to a subsequence not re-labeled)  the existence of $u_{(1)}$ and $u_{(2)}$ such that 
\begin{equation*}\label{proofGlowersemicontinuTDiscr}
\begin{gathered}
L \overline  U_\ell \rightharpoonup L u_{(1)}\;\: \textrm{and} \;\: \widehat U_\ell' \rightharpoonup u_{(2)} '\;\: \textrm{weakly in} \;\: L^2(0,T; L^2(\Omega))\, .
\end{gathered}
\end{equation*}
Notice that, for any space-time test function $\varphi\in C^\infty _0 $ one can prove
(we have set $\tilde \varphi ^n :=  \frac 1 {k_n}\int_{I_n } \varphi     \, \ d t $)
\begin{equation*}\label{IE_mr_estimate}
	\begin{split}
-	&\int_0^T  \<  \widehat U _\ell  , \varphi ' \>   \d t   = 
\int_0^T  \<  \widehat U_\ell ' , \varphi \>   \d t \\
& \ =	 - \int_0^T  \<  \overline U _\ell  , \varphi ' \>   \d t  
+ \sum _{n=1}^N \,    \<   \big [  U _\ell ^n -      U_\ell  ^{n-1}\big ], \big [ \tilde \varphi ^n  - \varphi ^{n-1} \big ]\> \, . \\
	\end{split}
\end{equation*}
By the uniform bound, 
$$\|\widehat U_\ell'\|^2_{L^2(0,T; L^2(\Omega))}  
=  \sum _{n=1}^N \, \frac 1 {k_n}\|  U_\ell  ^n -      U_\ell ^{n-1}\| ^2 _{L^2(\Omega)} \leq C_1 ^2 \, ,$$
{and the fact that the last sum is bounded by $C_1\big(\sum_n k_n\|\tilde\varphi^n-\varphi^{n-1}\|^2\big)^{1/2}=O(k_\ell)$, with $k_\ell:=\max_n k_n\to 0$,} we can conclude   that  $u_{(1)} = u_{(2)} = u$ and thus,  
$  \widehat U_\ell' + L \overline  U_\ell -  \overline  f \rightharpoonup u' + L u - f .$
Hence, 
\begin{equation*}\label{proofGlowersemicontinuityeq2-d2}
\int_{\Omega_T} | u' + L u -f |^2  \leq \liminf_{\ell \rightarrow \infty}  \int_{\Omega_T} |  \widehat U_\ell' + L \overline  U_\ell -  \overline  f|^2\, .
\end{equation*}


 \subsubsection*{$\limsup$ inequality and final convergence result}\label{TD:limsup}
 Let $ w \in H^1(0,T; L^2(\Omega))\cap L^2(0,T; H^2(\Omega)  \cap  H^1_0(\Omega))    $ The aim is to show  the existence of a recovery sequence $(w_\ell)$ such that $w_\ell \to w$ and 
   $$ \mathcal{G}(w) = \lim_{\ell \rightarrow \infty} \mathcal{G}_{IE, \ell} (w_\ell) .$$
   This is indeed possible by exploiting properties of the neural network spaces and using arguments based on smoothing and detailed estimates. The technical details of the proof can be found in \cite{GGM_IMA_2025}.
In the spirit of the abstract framework, we state explicitly what is assumed of the neural network spaces: for every $w \in H^1(0,T; L^2(\Omega))\cap L^2(0,T; H^2(\Omega)\cap H^1_0(\Omega))$ there exist $w_\ell\in V_{\ell,0}$ such that $w_\ell\to w$ in $H^1(0,T; L^2(\Omega))\cap L^2(0,T; H^2(\Omega))$, and $k_\ell\to 0$. Under this assumption the limsup inequality follows from the continuity of $\mathcal{G}$ and the consistency of the time quadrature, $\mathcal{G}_{IE,\ell}(w_\ell)-\mathcal{G}(w_\ell)\to 0$ for such $w_\ell$; the construction of the recovery sequence $w_\ell$
is   technical and it is addressed in detail in  \cite{GGM_IMA_2025}.

Then we can finally prove, 
that  if the implicit time-discrete loss is minimised over increasingly expressive network spaces with time steps tending to zero, the reconstructed minimisers converge to the solution of the parabolic problem, in the energy norm of the time-discrete functional:

\begin{theorem}[Convergence]
Let $u$ be the exact solution of the parabolic PDE and let $ (u_\ell) , $ $u_\ell \in V_{\ell}, $ be a sequence of minimisers of $ \mathcal{G}_{IE, \ell}  $, i.e.
\begin{equation*}\label{minofTDEnergy}
\mathcal{G}_{IE, \ell} (u_\ell) = \inf_{v_\ell \in V_\ell} \mathcal{G}_{IE, \ell} (v_\ell)\, .
\end{equation*}
Then, 
\begin{equation*}\label{TDConvOfDiscrMinE}
\hat u_\ell \rightarrow u,  \;\; \; \textrm{in} \;\: L^2(0,T; H^1(\Omega))  ,\;\ \  
\end{equation*}
where $\hat u_\ell $ is the piecewise linear reconstruction of ${u_\ell}$.
%
%
\end{theorem}

The proof relies on applying the abstract approach discussed in Section \ref{S:Analysis-abstract} in the present setting. For completeness we repeat the key steps: Let  $u\in L^2(0,T; H^2(\Omega))\cap H^1(0,T; L^2(\Omega)) $ be the   solution of \eqref{ParabolicPDE}. 
Consider the sequence of minimisers $(u_\ell)\, .$ Obviously, 
$$ \mathcal{G}_{IE, \ell}  (u_\ell) \leq \mathcal{G}_{IE, \ell} (v_\ell), \qquad \text{for all } v_\ell \in V_\ell\, . $$
In particular,  
$$ \mathcal{G}_{IE, \ell}  (u_\ell) \leq \mathcal{G}_{IE, \ell}  (\tilde u_\ell),  $$
where $\tilde u_\ell $ is the recovery sequence corresponding to $w=u.$  
Hence, we conclude that the sequence $ \mathcal{G}_{IE, \ell}  (u_\ell)$ is   uniformly  bounded. The \emph{discrete variational stability}  implies   
  that
\begin{align}
\|\overline u_\ell\|_{L^2(0,T; H^2(\Omega))} +
\|\widehat u_\ell'\|_{L^2(0,T; L^2(\Omega))} \leq C  \, .
\end{align} 
 We may assume by employing a density argument on the initial value that  $u^0\in H^2(\Omega),$ and thus  $\|\widehat u_\ell \|_{L^2(0,T; H^2(\Omega))} \leq C$ as well.    The  Aubin-Lions theorem ensures that there exists $   \tilde u\in L^2(0,T; H^1(\Omega)) $ 
such that $ \widehat u _\ell   \rightarrow   \tilde u$ in 
$L^2(0,T; H^1(\Omega))  $ up to a subsequence not re-labeled.
Furthermore the previous analysis shows that $L\tilde u \in L^2(0,T; L^2(\Omega))\, .  $
To prove that $   \tilde u $ is the minimiser of $ \mathcal{G} ,$ and hence $   \tilde u= u,$ we apply the framework of Section \ref{S:Analysis-abstract}: Let   $ w \in H^1(0,T; L^2(\Omega))\cap L^2(0,T; H^2(\Omega)) . $ We did show the existence of a recovery sequence $(w_\ell)$ such that $w_\ell \to w$ and 
   $$ \mathcal{G}(w) = \lim_{\ell \rightarrow \infty} \mathcal{G}_{IE, \ell}  (w_\ell) .$$
   Therefore,  the $\liminf$ inequality 
and the fact that $ u_\ell $ are   minimisers of the discrete problems imply that
\begin{align}
 \mathcal{G} (\tilde u)  
 \le  \liminf_{\ell \rightarrow \infty}   \mathcal{G} _{IE, \ell}  (u_\ell)  
 \le  \limsup_{\ell \rightarrow \infty} 
 \mathcal{G}_{IE, \ell}  (u_\ell) 
 \le  \limsup_{\ell \rightarrow \infty} 
 \mathcal{G} _{IE, \ell}  (w_\ell)  =  
 \mathcal{G}  (w) , 
\end{align}
for all $ w \in H^1(0,T; L^2(\Omega))\cap L^2(0,T; H^2(\Omega)) . $  Therefore $  \tilde u $ is   the  minimiser of $  \mathcal{G} ,$ hence 
$   \tilde u= u$ and the entire sequence satisfies 
$$
\hat u_\ell \rightarrow u,  \;\; \; \textrm{in} \;\: L^2(0,T; H^1(\Omega))  ,\;\ \  
$$
which shows the desired result. 

 \subsubsection{Explicit time discrete functional}
It will be interesting to consider a seemingly similar (from the point of view of quadrature and approximation) discrete functional:
\begin{equation}\label{Functional_k_Ex}
\mathcal{G}_{k, EE} (v) = \sum _{n=1}^N \,  k_n \, \int_{\Omega }\big | \frac{v ^n-  v^{n-1}}{k_n} + L v ^{n-1}  - f^{n-1} \big |^2 \, \, \d x +\,     {|  v^0  - u^0 |_{H^1(\Omega)}^2},
\end{equation}
and compare its properties with those of the functional $\mathcal{G}_{k, IE} (v) $ and the corresponding $V_{\mathcal{N}}$ minimisers. The functional \eqref{Functional_k_Ex} is associated with \emph{explicit} Euler discretisation in time, as opposed to the \emph{implicit} Euler discretisation underlying $\mathcal{G}_{k, IE} (v).$ Clearly, in the discrete minimisation framework, both energies are fully implicit, since the evaluation of the minimisers requires the solution of global space-time problems. It is therefore particularly interesting that the two energies lead to fundamentally different stability properties.
Let us first note that a discrete coercivity estimate such as \eqref{GEquicoercivity_TD} does not appear to be possible. Indeed, an argument similar to  {\eqref{dMR_TD}} can be used, but with the crucial difference that the second-to-last term in this relation is negative rather than positive. This is directly related to the well-known instability of the forward Euler method. For finite difference forward Euler schemes, stability typically requires a strong CFL condition of the form $k\leq C h^2 $, where $h $ is the spatial discretisation parameter. A similar phenomenon appears to occur here. Although we do not derive stability bounds when spatial training is performed, our numerical experiments show that the stability of the explicit training method deteriorates as the number of spatial training points increases while $k$ is kept fixed. {These computations, presented in \cite{GGM_IMA_2025}, provide convincing evidence that }
coercivity bounds similar to  \eqref{GEquicoercivity_TD} are necessary for  stable behaviour of the approximations, for details see \cite{GGM_IMA_2025}.

\subsubsection*{The elliptic case}
The two variational principles for the elliptic problem of
Section~\ref{SubS:Ex:EllipticPDE} are not analysed in this section. The
Deep Ritz method with training based on finite element quadrature is
formulated and analysed in Section~\ref{SubS:FE-training}. For the Deep Ritz
method with Monte Carlo training, the sampled energy controls $g$ and
$\nabla g$ only at the sample points, and a coercivity bound of the type
required by discrete variational stability is not readily available; the
same holds for the residual-based (PINN) formulation with Monte Carlo
training. The analysis of these methods through discrete variational
stability is rather involved and is not presented here. For the PINN loss
we derive instead a posteriori estimates, which are the subject of
Section~\ref{S:reliability}.

\subsection{Bibliography and comments}\label{SubS:Conv:bib}

The framework developed in this section is based on the works
\cite{loulakis2023newap} and \cite{GGM_IMA_2025}. The liminf--limsup
framework is the core of De Giorgi's $\Gamma$-convergence
\cite{DeGiorgi_Franzoni:1975, DalMaso_Gamma:1993, braides2002gamma}, used
here, for the reasons explained in Section~\ref{S:Abstract-liminf-limsup},
without invoking the $\Gamma$-convergence theorems themselves; the weak
convergence of functions with respect to varying measures, which is the
topology adapted to probabilistic training, follows
\cite[Chapter 5]{AGS:2005}. \\

Several important contributions over the last years have been devoted to
the approximation capability of DNN spaces, see
\cite{guhring2019errorboundsapproximationsdeep, Dahmen_Grohs_DeVore:specialissueDNN:2022,Schwab_DNN_constr_approx:2022, Xu_appr_cnn:2022, Xu_approx:2022, Schwab_DNN_highD_analystic:2023, Mishra:appr:rough:2022, Grohs_Petersen_Review:2023}
and their references. New results are needed, among other topics, to
identify the neural network architectures required to achieve given
approximation bounds, rather than merely establishing their existence.
Nevertheless, the remarkable advances in the approximation theory of neural
networks constitute a major milestone in our understanding of their
behaviour.\\

Regarding rates in generalisation estimates one can distinguish  three different notions in this context. \emph{Approximation rates} concern the spaces $V_{\mathcal{N}}$ alone, that is, how well the class approximates a target of given regularity; these are now well developed as just mentioned. \emph{Rates for the statistical error} concern the distance between the trained and the exact functional, and are the subject of Section~\ref{S:reliability}. \emph{Rates for the discretisation} concern the error of the minimiser of the actual trained loss, in a continuous norm, and require, as in classical numerical analysis, a stability estimate and a consistency estimate with rates. Rates of this third kind are, to our knowledge, not available for neural network methods beyond the Monte Carlo function learning problem on compact sets, where they are obtained in statistical learning theory for empirical risk minimisers over uniformly bounded classes whose size is tied to the number of samples; even for deterministic quadrature in a continuous norm, where the discrete norm controls nothing between the quadrature points and no inverse estimate is available for network spaces, and for time-discrete schemes, they are open. What is usually presented as a rate for the method is an approximation rate combined with a uniform assumption on the class. This is why the framework of these notes concentrates on stability first: rates for the discretisation cannot be obtained without it.\\

Previous works on the analysis of methods based on residual minimisation
over neural network spaces for PDEs include \cite{SSpiliopoulos:2018},
\cite{Karniadakis:pinn:conv:2019}, \cite{shin2020error}, \cite{hong2022priori},
and \cite{Mishra_gen_err_pinn:2023}. In \cite{SSpiliopoulos:2018},
convergence was established for sufficiently smooth classical solutions of
a class of nonlinear parabolic PDEs, without explicitly considering the
training of the functional. Convergence results under assumptions on the
discrete minimisers or on the neural network space, in the setting of Monte
Carlo training, were obtained in \cite{Karniadakis:pinn:conv:2019},
\cite{shin2020error}, and \cite{hong2022priori}. In addition,
\cite{shin2020error} makes use of continuous stability properties of certain
linear operators in the analysis. In \cite{muller2020deep}
$\Gamma$-convergence was used in the analysis of deep Ritz methods without
training.
Stability based on inf-sup frameworks for Variational PINNs (the residuals
are evaluated in a weak-variational sense), were considered in
\cite{B_Canuto_P_Vpinn_quadrature_2022}, where the role of quadrature was
proven crucial in the analysis of the method.

The results of \cite{Mishra_gen_err_pinn:2023} are based on estimates whose
bounds depend on the discrete minimisers and their derivatives. Convergence
can therefore be obtained under suitable uniform bounds on these functions
in appropriate Sobolev norms. The results of \cite{hong2022priori}, which
consider deterministic training, are related in this respect, as they apply
to neural network spaces for which high-order derivatives are uniformly
bounded in appropriate norms by construction.

When training and quadrature are incorporated into the design and analysis
of the method, one possible approach is to treat the resulting training or
quadrature error as an additional contribution to C\'ea-type estimates. Such
estimates naturally lead to assumptions controlling higher-order
derivatives of the discrete minimisers. Our approach takes a somewhat
different perspective: rather than treating the training procedure as an
additional source of error, we view the entire scheme, including training,
as a unified discrete variational problem and study its stability directly.
This viewpoint is closer to the standard a priori analysis of numerical
methods for differential equations, where the stability of the fully
discrete scheme is analysed as an intrinsic property of the method.

Error estimates for operator learning, for deep operator networks and
Fourier neural operators, are derived in \cite{Mishra_DeepONET_Kran_est,
Mishra_FNO_est}, see also \cite{Stuart_NeuralOper2023} for approximability of neural operators; they
decompose the error into encoding, approximation and reconstruction parts,
and control the generalisation part through covering numbers of the operator
class, which requires uniform Lipschitz bounds on the networks. In  Section~\ref{SubS:Conv:MCInterpolation} we indicate how the generalisation
part can be treated instead by the liminf inequality, without such bounds, by employing the
setting  
of measures on the infinite-dimensional
output space, as developed in \cite[Chapter 5]{AGS:2005}.

The minimisation of a distance to the empirical measure over a parametrised
family of measures is, in statistical terms, minimum distance estimation,
and large-sample results for GANs in this spirit are given in
\cite{Biau_Cadre_Sangnier_Tanielian:2020}; the approximation property of the
neural measure spaces assumed in Section~\ref{SubS:Conv:measures}, that
every measure in $\mathcal{P}_p(\R^k)$ is a $W_p$-limit of push-forwards of
$\gamma_0$ through neural networks, is a universal approximation property
for distributions, see e.g.\ \cite{Lu_Lu:2020} for push-forwards through
gradients of ReLU networks in Wasserstein and other integral probability
metrics, and \cite{Perekrestenko_Muller_Bolcskei:2020} for the generation
of high-dimensional distributions from a one-dimensional source. For
entropic optimal transport see \cite{Cuturi:2013} for the Sinkhorn
algorithm, \cite{Genevay_Peyre_Cuturi:2018} for its use in generative
modelling, and \cite{Leonard:2012, Carlier_Duval_Peyre_Schmitzer:2017} for
the convergence of the regularised cost to the Wasserstein distance as the
regularisation parameter tends to zero. The bias $\mathcal{T}_\delta(\mu,\mu)>0$
is removed by the Sinkhorn divergence
$S_\delta(\mu,\nu)=\mathcal{T}_\delta(\mu,\nu)-\tfrac12\mathcal{T}_\delta(\mu,\mu)-\tfrac12\mathcal{T}_\delta(\nu,\nu)$,
which is nonnegative, vanishes only for $\mu=\nu$, metrises weak convergence,
and interpolates between $W_2^2$ and a maximum mean discrepancy
\cite{Feydy_etal:2019}. The sample complexity of the regularised loss is of
order $N^{-1/2}$, with a constant depending on $\delta$ and on the data
dimension, which blows up as $\delta\to0$ \cite{Genevay_etal_sample:2019};
in the a posteriori estimate \eqref{Wp:apost:entropic} this trades the
dimension-dependent rate of the empirical term for the regularisation error
$O(\delta\log(1/\delta))$, and the choice of $\delta$ in terms of $N$ is a
balance between stability and consistency of the kind familiar from
numerical analysis. The approximation of the distance by an integral
probability metric over a class of discriminator networks, as in Wasserstein
GANs \cite{arjovsky2017wassersteingan}, is different in nature: the loss is
then no longer a metric, and the coercivity and lower semicontinuity of the
discriminator-based losses, uniformly in the discriminator class, would have
to be established.

\newpage
    
 \section{Reliability through estimates}\label{S:reliability}
Error estimates, whether of \emph{a priori} or \emph{a posteriori} character, provide a useful means of assessing the reliability of numerical algorithms. As we have seen, numerical computations may lead to unreliable predictions for several reasons: the algorithm may lack stability, or the underlying problem itself may be sensitive to relatively small perturbations. Error estimates provide a way to quantify and, in favourable circumstances, control these effects.

\emph{A priori} estimates serve a twofold purpose. First, when the underlying continuous problem is well posed and the numerical method satisfies the appropriate stability and consistency properties, they provide a theoretical guarantee that, at least asymptotically, the numerical solution approaches the correct solution. Second, they have an important practical role in the design and assessment of numerical experiments. When computations are performed with different values of the discretisation parameters, we expect the observed errors to behave consistently with the theoretical estimates. In particular, the numerical results should exhibit the predicted convergence rates as the discretisation is refined.

\emph{A posteriori} estimates are particularly valuable when a complete convergence theory is unavailable, as is often the case for machine learning algorithms. As discussed in Section~\ref{S:apriori_apost}, upper bounds for the generalisation error that are computable from quantities available during training can serve as a direct indicator of reliability, without requiring a full \emph{a priori} analysis. This is especially consequential in scientific applications, where a small training loss or good empirical performance alone does not guarantee an accurate or physically meaningful prediction. In such settings, \emph{a posteriori} estimates provide an important additional layer of validation, allowing the reliability of the computed solution to be assessed directly from the output of the algorithm. 

   
%
%
%
%
%
%

   \subsection{Estimates for MC Interpolation}\label{SubS:Rel:MCInterpolation}
 \def\E{{\mathcal{E}}}
	 
Recall the problem of MC-Deep neural network interpolation: 
Let $X=(X_1,X_2, \dots, X_N)$ be i.i.d.\ random variables with common law $\mu$, each one corresponding to a random point $x_i$ of $\Omega.$\\
 We seek minimisers of  
\begin{equation*}\label{ieE-minimize_NNQ_prob-d7}
	\min  _ {v \in  V _{\mathcal{N}} } \E _{N,\omega}(v), 
\end{equation*}
which we denote by $u_X \in V _{\mathcal{N}} .$ 
A central problem in foundational machine learning is to estimate the generalisation error, i.e., the error $\|u_{\bm {X}(\omega)}\!-\!f\|_p,$ since it is a  measure for how well  the algorithm behaves beyond the training points. The following bound is trivial: 
Let  $f\in L^p(\mu),\ \omega\in \bbOmega, $ and $\varphi   \in  V _{\mathcal{N}},$ since $ \E _{N,\omega}(u_{\bm {X}(\omega)})\leq  \E _{N,\omega}(\varphi),$
\begin{equation*}
 \begin{split}
  \|u_{\bm {X}(\omega)}\!-\!f\|^p_p \leq &   \frac 1N  \sum _{i}   \Big  | \varphi  (X_i(\omega)\, ) - f(X_i (\omega) \, ) \Big  | ^p   \\
  & + 
 \int_{\Omega}    |  u_{\bm{X}  (\omega)   } - f|^p    \d \mu(x) -  \frac 1N  \sum _{i}   \Big  | u_{\bm{X} (\omega)   } (X_i(\omega)\, ) - f(X_i (\omega) \, ) \Big  | ^p  \\
= &    \|\varphi -\!f\|^p_p  \\
   & + \Big \{
 \int_{\Omega}    |  u_{\bm{X}  (\omega)   } - f|^p    \d \mu(x) -  \frac 1N  \sum _{i}   \Big  | u_{\bm{X} (\omega)   } (X_i(\omega)\, ) - f(X_i (\omega) \, ) \Big  | ^p  \,  \Big \} \\
& +\Big \{
 \frac 1N  \sum _{i}   \Big  | \varphi  (X_i(\omega)\, ) - f(X_i (\omega) \, ) \Big  | ^p  -  \int_{\Omega}    | \varphi - f|^p \d\mu \Big \}
 \end{split}
\end{equation*}
In order to  study the behaviour of the generalisation error
we need to control
\begin{itemize}
	\item  the error due to the minimisation over a DNN discrete space (approximation error)
	\item  the error due to the approximation of the continuous functional which is associated to $f$ by the sample averages: 
\begin{equation*}\label{prob_E_error-d5}
 {E}_{MC} (u_{\bm{X} (\omega)    }) =  \int_{\Omega}    |  u_{\bm{X}  (\omega)   } - f|^p    \d \mu(x) -  \frac 1N  \sum _{i}   \Big  | u_{\bm{X} (\omega)   } (X_i(\omega)\, ) - f(X_i (\omega) \, ) \Big  | ^p  \, . 
\end{equation*}
\end{itemize}
The key technical issue here is that since $u_{\bm{X} (\omega)   }$ is a solution of a minimisation problem involving $( X_i ) ,$ it obviously depends on  $( X_i ) ,$ and thus the function $u_{\bm{X} (\omega)    }$ changes with $N.$ Thus, the law of large numbers and other standard  tools estimating ${E}_{MC}  ( g) $ for fixed $g$ are not applicable \cite{Kallenberg_book}.

%
%
%
%
%
%
%
%
%

 \subsubsection*{Rademacher complexity}
	In statistical learning, one may   overcome this problem  by uniform estimates within a class, i.e. by estimating \cite{Bartlett_Montanari_Rakhlin_actaN:2021}
\begin{equation*} 
\big | {E}_{MC}(u_{\bm{X}     }) \big | \leq \sup _ {g \in  V _{\mathcal{N}} }   \, \big | {E}_{MC}(g ) \big |\, ,   
\end{equation*}
then
\begin{equation*}
 \begin{split}
  \|u_{\bm {X}(\omega)}\!-\!f\|^p_p \leq &    \|\varphi -\!f\|^p_p  + 2 \sup _ {g \in  V _{\mathcal{N}} }   \, \big | {E}_{MC}(g ) \big | \, .
 \end{split}
\end{equation*}

Then it is known that the \emph{Rademacher complexity} of a function class $\mathcal G$
\begin{itemize}
	\item    called $R_ {N, \mathcal G } $ 
$$R_ {N, \mathcal G } := \EE \sup_{g \in \mathcal G}
\frac 1N  \sum _{i} 
\epsilon _i g(X_i ),$$
where $\epsilon _i  \in \{\pm 1 \}$ are independent and uniformly distributed, is known to control $\sup _ {\mathcal G }   \, \big | {E}_{MC}(g ) \big |\, .$ 
	\item Certain structural assumptions needed to control $R_ {N, \mathcal G } $ on   sets  $\mathcal G\ .$

\end{itemize}
There are similar alternative techniques in the \emph{statistical learning} literature that aim to control this error, typically uniformly with respect to the class $\mathcal G$. These approaches can sometimes be quite technical, but the theory is relatively mature by now and provides reasonable   asymptotic rates with respect to the dimension $d$.

However, the backbone of this approach relies on uniform, and hence structural, assumptions on the neural network spaces $V _{\mathcal{N}}$. In applications to PDEs and differential equations, the function $g$ typically involves derivatives arising from the differential operator. Consequently, applying such an approach can be particularly restrictive, as it relies on structural assumptions on the discrete spaces that may be difficult to satisfy in practice.

\subsubsection*{Estimates and the optimal matching problem}	
	
	We suggest an alternative approach that aims to avoid uniform bounds over a class. The key idea is to represent the loss in terms of an empirical measure and estimate the corresponding error using the distance between the empirical measure and the measure $\mu$. This, in turn, requires estimates for the \emph{optimal matching problem} in the setting under consideration. The resulting framework provides an alternative approach to error control, with the potential to yield new results and provide further insight into the behaviour of the algorithm.

	 To this end, as before, we rewrite the loss functional
\begin{equation*}\label{prob_E-d3}
\mathcal{E}_{N, \omega}( g )  =  \frac 1N  \sum _{i}   \Big  | g(X_i(\omega)\, ) - f(X_i (\omega) \, ) \Big  | ^p  \, 
\end{equation*}
using  the empirical measures on $\Omega,$
$	\mu _{N, \bm{X} (\omega)} =   \frac 1 N \sum _{i} \   \delta _{X_i(\omega)} \, ,
$
as
\begin{equation*}\label{energy_emp_meas-d3}
	\mathcal{E}_{N, \omega}( g )  = \int _\Omega |g (x)-f(x)|^p\, \d \mu _{N, \bm{X} (\omega)} (x) 
\end{equation*}
then, for $\omega \in \bbOmega , $
\begin{equation*}\label{prob_E_error-d6}
 {E}_{MC} (u_{\bm{X} (\omega)    }) =  \int_{\Omega}    |  u_{\bm{X}  (\omega)   } - f|^p    \d x -  \int _\Omega | u_{\bm{X}  (\omega)   }-f|^p\, \d \mu _{N, \bm{X} (\omega)}   \, . 
\end{equation*}

Recall that the Wasserstein distance of order $p$ between 
two probability measures $\mu$ and $\nu$ 
is defined by
\begin{equation*}\label {W-dist-d2}
	W_p(\mu, \nu) = \Big \{ \min_{\pi \in A (\mu, \nu)} \int _{\Omega\times \Omega} |x-y|^p\, \d \pi(x,y) \Big \} ^{1/p} ,
\end{equation*}
where $A(\mu,\nu)$ denotes the set of \emph{transport plans} between $\mu$ and $\nu.$ Then, utilising the push-forward measure and the formulas \eqref{pf-int} and \eqref{W-dist-pf} we observe:  
%
%
%
%
%
%
%
%
%
%
%
The triangular inequality in $L^p(\pi)$, for any $\pi\in A(\mu, \nu)$ we have
\begin{align*}
\big|\ \|\xi\|_{L^p(\mu)}-&\|\xi\|_{L^p(\nu)}\big|\\
&=\big| \big\{\int _{\Omega \times \Omega} |\xi (x) |^p\,  \d \pi(x,y)\big\}^{1/p}-\{\int _{\Omega\times \Omega} |\xi (y) |^p\,  \d \pi(x,y)\big\}^{1/p}\big|\\
&\leq \{\int _{\Omega\times \Omega} |\xi (x)-\xi(y) |^p\,  \d \pi(x,y)\big\}^{1/p}.
\end{align*}
Optimizing over $\pi\in A(\mu,\nu)$, we get that 
\begin{equation*} 
\big|\ \|\xi\|_{L^p(\mu)}-\|\xi\|_{L^p(\nu)}\big| \leq W_p(\xi _\sharp  \mu, \xi _\sharp \nu)
\end{equation*}
Now, set $\xi=u_{\bm {X}(\omega)}-f$ and $\nu=\mu_{N,\bm{X}(\omega)}$ to get
\begin{equation*} 
\begin{split}
 \|u_{\bm {X}(\omega)}\!-\!f\|_p \leq \Big \{  &\frac 1N  \sum _{i}   \Big  | u_{\bm {X}(\omega)} (X_i(\omega)\, ) - f(X_i (\omega) \, ) \Big  | ^p  \Big \}^{1/p}\\
 &+
 W_p\big([u_{\bm {X}(\omega)}\!-\!f\,]_\sharp\mu, [u_{\bm {X}(\omega)}\!-\!f\,]_\sharp\mu _{N, \bm{X} (\omega)}\big).	
\end{split}
\end{equation*}
Then, since 
$u_{\bm {X}(\omega)}$ is a minimiser of the discrete problem we obtain for any $\varphi   \in  V _{\mathcal{N}},$
\begin{equation*}\label {W-dist-d_0_3}
\begin{split}
\|u_{\bm {X}(\omega)}\!-\!f\|_p
 \leq &\Big \{  \frac 1N  \sum _{i}   \Big  | \varphi  (X_i(\omega)\, ) - f(X_i (\omega) \, ) \Big  | ^p  \Big \}^{1/p}\\
 &+  W_p(\big [u_{\bm {X}(\omega)}-f\, \big ]_\sharp\mu, \big [u_{\bm {X}(\omega)}-f\, \big ]_\sharp\mu _{N, \bm{X} (\omega)})
  \, .
\end{split}
\end{equation*}
In addition 
\begin{equation*}\label{prob_E_min_prop2}\begin{split} 
\int _\Omega \,  \frac 1N  \sum _{i}   \Big  | \varphi  (X_i(\omega)\, ) - f(X_i (\omega) \, ) \Big  | ^p   \, \d \P (\omega) \\
  = \int_\Omega|\varphi(x)-f(x)|^p\, \d \mu (x)\, . 	\end{split}
\end{equation*}
{Since $x\mapsto x^{1/p}$ is concave, Jensen's inequality gives $\EE[\mathcal{E}_{N,\cdot}(\varphi)^{1/p}]\le \|\varphi-f\|_p$; together with the measurability of $\omega\mapsto u_{\bm X(\omega)}$, assumed in Section~\ref{SubS:ModP:function-learning}, this yields the estimate in expectation below. We conclude therefore to the estimate:}

	 \begin{theorem}[\sc estimate of the generalisation error]
Consider for each $\omega \in \bbOmega$ the problem 
\begin{equation*}
	\min  _ {v \in  V _{\mathcal{N}} } \E _{N, \omega} (v)
\end{equation*}
and its solution     $ u_{\bm{X} (\omega)   },    $ the  $V_{\mathcal{N}}-$\emph{probabilistic deep neural network interpolant of $f$ at $ \bm{X}  = (X_1, \dots, X_N  )  \, .$}
Then for each $f\in L^p(\mu),\ \omega\in \bbOmega, $ and $\varphi   \in  V _{\mathcal{N}},$
\begin{equation*}
 \begin{split}
  \|u_{\bm {X}(\omega)}\!-\!f\|_p \leq &\Big \{  \frac 1N  \sum _{i}   \Big  | \varphi  (X_i(\omega)\, ) - f(X_i (\omega) \, ) \Big  | ^p  \Big \}^{1/p}\\
  & + 
W_p\big([u_{\bm {X}(\omega)}\!-\!f\, ]_\sharp\mu, [u_{\bm {X}(\omega)}\!-\!f\, ]_\sharp\mu_{N, \bm{X} (\omega)}\big).
 \end{split}
\end{equation*}
In particular,
\begin{equation*}
 \begin{split}
 \EE\big[ \|u_{\bm {X}(\cdot)}-f\|_p \big]  \leq & \inf_{\varphi   \in  V _{\mathcal{N}}}  \|\varphi-f\|_p\\
 & + 
 \EE\big[  W_p\big([u_{\bm {X}(\cdot)}\!-\!f\, ]_\sharp\mu, [u_{\bm {X}(\cdot)}\!-\!f\, ]_\sharp\mu_{N, \bm{X} (\cdot)}\big)  \big]\, .
 \end{split}
\end{equation*}
\end{theorem}

  This simple estimate involves two terms: the first is a pure approximability term reflecting the expressive capacity of the $V_{\mathcal{N}}$ space, while the second measures the Monte-Carlo integration error in this non-trivial setting. The latter is closely related to \emph{optimal matching problems}, a well-studied class of random problems in probability and mathematical physics concerning the distance between empirical measures and their underlying law, or between sets of randomly sampled points in $\R^d$. See, for example,  \cite{FournierGuillin:2015}; \cite{AjtaiKT:1984};\cite{AmbrosioST:2019,AGS:2005,AmbrosioGT:2019,AmbrosioGT:2022,Parisi:matching:2014,Talagrand:1994,Talagrand:1992,Talagrand:AKTgeneral:1992,Talagrand:book}, see also \cite{PanaretosWR:2019} for the relevance of these bounds in statistics. 
  
  The asymptotic behaviour of $\EE [W_p(\nu,\nu_N)]$, where $\nu_N$ is the empirical measure associated with $N$ independent samples drawn from $\nu$, has been studied extensively. Recent PDE-based approaches provide upper and lower bounds, see \cite{AmbrosioST:2019,AmbrosioGT:2019,AmbrosioGT:2022}, suggesting that, for $d\geq 2$,
\begin{equation}\label{emp_W_rates0}
\EE [W_p^p(\nu,\nu_N)]\asymp N^{-p/d},
\end{equation}
up to a logarithmic factor in the critical case $d=2$. 
In particular, \cite{FournierGuillin:2015} yields, for $d>2p$ and under suitable moment assumptions on $\nu$,
\begin{equation}\label{emp_W_rates}
\EE [W_p(\nu,\nu_N)]\leq C(\nu)N^{-1/d},
\end{equation}
with a logarithmic correction in the critical case $d=2p$. The dimension that enters here is in fact the intrinsic dimension of $\nu$ rather than that of the ambient space: by \cite{Weed_Bach:2019}, the rate is $N^{-1/d^*}$ with $d^*$ the upper Wasserstein dimension of $\nu$, which is smaller than $d$ when $\nu$ is concentrated on a lower-dimensional set; see also \cite{hundrieser2022empirical} for related results on empirical optimal transport. For a measure $\nu$ chosen by us, such as the uniform distribution on $\Omega$, $d^*=d$ and there is no gain; for data-driven problems, in which $\nu$ is the data distribution, the gain can be substantial.

The second term in the estimate can be bounded through the inputs: for a Lipschitz map $\xi$ and any plan $\pi\in A(\mu,\nu)$, $\int|\xi(x)-\xi(y)|^p\,\d\pi\le \mathrm{Lip}(\xi)^p\int|x-y|^p\,\d\pi$, so that, by \eqref{W-dist-pf},
\begin{equation}\label{Lip_pushforward}
W_p(\xi_\sharp\mu,\xi_\sharp\nu)\le \mathrm{Lip}(\xi)\,W_p(\mu,\nu)\, .
\end{equation}

{The next result essentially states that under a uniform bound on the Lipschitz constants of the computed minimisers, the expected generalisation error is bounded by the best approximation error of the network space plus the rate at which empirical measures approach their law.} We thus have:
\begin{proposition}
[\sc asymptotic bound of the generalisation error]
\label{Thm:estimate}
Let    $ u_{\bm{X} (\omega)   }    $ be the  $V_{\mathcal{N}}-$\emph{probabilistic deep neural network interpolant of $f$ at $ \bm{X}  = $$(X_1, \dots,$ $ X_N  )  \, .$}
Assume that $f$ is Lipschitz and let us denote the Lipschitz constant of $ u_{\bm{X} (\omega)   }  -f$ by $L_{\bm{X} (\omega)   }\, .$ 
We denote by $E_N (\nu),$
$$ E_N (\nu)= \EE [   W_p(\nu, \nu _{N})] \, .$$
Then, if, $\mathbb{P}$-a.s.,\begin{equation}\label {W-dist-ass}
L_{\bm{X}}\,  \leq L_N
\end{equation}
we have
\begin{equation}\label {W-dist-d_0}
 \EE\big[ \|u_{\bm {X}(\cdot)}-f\|_p \big]  \leq  \inf_{\varphi   \in  V _{\mathcal{N}}}  \|\varphi-f\|_p
  + 
L_N\, E_N (\mu) \, .
\end{equation}
Furthermore, if 
\begin{equation}\label {W-dist-d_01}
   \inf_{\varphi   \in  V _{\mathcal{N}}}  \|\varphi-f\|_p \to 0 \quad \text{and} \quad
L_N\,   E_N (\mu) \to 0 \, , 
\end{equation}
as the network space is enriched and $N \to \infty$, respectively,
we conclude that 
\begin{equation}\label {W-dist-d_02}
 \EE\big[ \|u_{\bm {X}(\cdot)}-f\|_p \big] \to 0\, .
\end{equation}	
\end{proposition}
Two remarks are in order. The uniform bound \eqref{W-dist-ass} on the Lipschitz constants of the computed interpolants is an assumption of the kind discussed in Section~\ref{SubS:aPaP:Error_analysis}: it is not derived from the scheme, and it may be restrictive. The per-instance estimate of the theorem does not require it, and in its Lipschitz form \eqref{Lip_pushforward} the constant $L_{\bm X(\omega)}$ is that of the \emph{computed} network, a quantity that can be measured a posteriori.
\bigskip

In view of \eqref{emp_W_rates}, the bound \eqref{W-dist-d_0} has an unfavourable dependence on the dimension $d$. However, this results from a crude estimate of the term

$$
W_p(\big [u_{\bm {X}(\omega)}-f\, \big ]_\sharp\mu, \big [u_{\bm {X}(\omega)}-f\, \big ]_\sharp\mu _{N, \bm{X} (\omega)}).
$$
Indeed, the measure $\nu=\big [u_{\bm {X}(\omega)}-f, \big ]_\sharp\mu$ is one-dimensional, while
$$
\big [u_{\bm {X}(\omega)}-f\, \big ]_\sharp\mu _{N, \bm{X} (\omega)}
=\frac{1}{N}\sum_i\delta_{[u_{\bm {X}(\omega)}-f](X_i(\omega))}
$$
Note that it is \emph{not} the empirical measure of i.i.d.\ samples of $\nu$, since $u_{\bm X(\omega)}$ depends on all the $X_i$. For genuine empirical measures in one dimension, \cite{FournierGuillin:2015} yields
\begin{equation}\label{emp_W_rates_d=1}
\EE [W_p^p(\nu,\nu_N)]\leq C(\nu)N^{-1/2},\qquad\text{i.e.}\qquad \EE [W_p(\nu,\nu_N)]\leq C(\nu)N^{-1/(2p)},
\end{equation}
and the sharper rate $\EE [W_p(\nu,\nu_N)]\leq C(\nu)N^{-1/2}$ holds under additional conditions on $\nu$, for instance for measures with a density bounded away from zero on an interval; the one-dimensional case is studied in detail in \cite{Bobkov_Ledoux:2019}. One may therefore expect, under suitable conditions on the law of $u_{\bm X(\omega)}-f$,
\begin{equation}\label{emp_meas_bound}
	\EE [\,  W_p(\big [u_{\bm {X}(\omega)}-f\, \big ]_\sharp\mu, \big [u_{\bm {X}(\omega)}-f\, \big ]_\sharp\mu _{N, \bm{X} (\omega)}) \, ]  \leq C \,    N^{-1/2}\, .
\end{equation}
Establishing \eqref{emp_meas_bound} is an interesting open problem; preliminary numerical evidence suggests that the observed behaviour is indeed compatible with this rate.

An important practical consequence is that the Wasserstein distances appearing in these estimates are distances between measures on the real line, for which $W_p$ between two empirical measures is computed exactly by sorting, through the quantile representation of one-dimensional optimal transport, at a cost of $O(N\log N)$; no approximate transport algorithm is needed. Thus, beyond providing theoretical error bounds, our framework offers a route to the \emph{empirical estimation of the generalisation error}. The term $W_p([u_{\bm X}-f]_\sharp\mu,[u_{\bm X}-f]_\sharp\mu_{N,\bm X})$ itself is not computable, since it involves $f$ everywhere, but the same argument with an independent validation sample $Y_1,\dots,Y_M$ of $\mu$, independent of $\bm X$, with empirical measure $\mu'_M$, gives
\begin{equation}\label{validation_estimate}
\|u_{\bm X}-f\|_p\le\Big\{\frac1M\sum_{j=1}^M\big|u_{\bm X}(Y_j)-f(Y_j)\big|^p\Big\}^{1/p}
+W_p\big([u_{\bm X}-f]_\sharp\mu,[u_{\bm X}-f]_\sharp\mu'_M\big)\, .
\end{equation}
Conditionally on the training sample $\bm X$, the function $u_{\bm X}-f$ is fixed and $\mu'_M$ is the empirical measure of genuine i.i.d.\ samples, so that the last term is a one-dimensional empirical term to which \eqref{emp_W_rates_d=1} applies, with a rate that does not depend on $d$. The first term is the validation loss used in practice; what the one-dimensional push-forward adds is an explicit, dimension-free control of its error. The open problem \eqref{emp_meas_bound} concerns the training term only.
This idea is particularly relevant to the PDE learning methods discussed next.

\subsection{A posteriori estimates for PINNs}
Now we turn to a  linear elliptic   PDE, i.e.,    
   boundary value problems  of the form,
\begin{equation}\label{ellipticPDE}
\left \{
\begin{alignedat}{3}
&L\, u = f\quad &&\text{in}\,\,&&\Omega\\[2pt]
&u =0 \quad &&\text{on}\,\,&&\partial \Omega \\[2pt]
\end{alignedat}
\right .
\end{equation}
where $ u : \Omega \subset \mathbb{R}^d \rightarrow \mathbb{R} \;, \: \Omega  $ is an open, bounded set with smooth enough boundary, 
$f \in L^2 (\Omega )  $ and $ L $ a self-adjoint elliptic operator as defined in Section~\ref{SubS:Ex:EllipticPDE}. 
   
As before, 
  we consider a collection $Y, X_1,X_2,\ldots$ of i.i.d.\ $\Omega$-valued random variables    with common law $\mu.$ 
Further we shall need boundary points through $  X^b_1,X^b_2,\ldots$ of i.i.d.\ $\partial \Omega$-valued random variables, defined on a  corresponding probability space 
with common law $\mu_S, $ and  
for $g\in L^p(\partial \Omega,\mu_S)$  
\begin{equation}\label{repr_one-d2}
	\int _{\partial \Omega} \, g (x) \, \d \mu_S  (x) = \int _{\tilde \bbOmega} \, g ( X^b_i (\tilde \omega)) \, \d \tilde \PP (\tilde \omega) \, .\end{equation} 
Monte Carlo integration is formulated  as follows: let $\omega \in \varOmega ,\ \tilde \omega \in \tilde \varOmega $ be fixed. Consider the discrete energy, 
\begin{equation}\label{prob_E-d4}
\mathcal{E}_{N, \omega, {\tilde \omega}}( g )  =  \frac 1N  \sum _{i}   \Big  | \, L\, g(X_i(\omega)\, ) - f(X_i (\omega) \, ) \Big  | ^2
+\frac 1{N _b} \sum _{j}   \Big  | g(X^b_j(\tilde \omega)\, ) \Big  | ^2  \, . 
\end{equation}

The method considered,   in this section is:   for each instance  
 $\omega \in \bbOmega, \tilde \omega \in \tilde \bbOmega$ we seek solutions to the minimisation problem%
\begin{equation}\label{prob_E_def-d2}
	\min  _ {v \in  V _{\mathcal{N}} } \mathcal{E}_{N, \omega, {\tilde \omega}} (v)\, . 
\end{equation}
 We denote by   $ u_{\bm{X} (\omega, {\tilde \omega} )   }  = u ^\star(\omega, \tilde \omega, \cdot)  \,  $ the  $V_{\mathcal{N}}-$\emph{ probabilistic PINN solution  of $Lu = f$ at $ \bm{X}  = (\, (X_1, \dots, X_N  ), (X^b_1, \dots, X^b_{N_b})\,   )  \, .$}
%

As in the previous section we introduce corresponding empirical measures. Now we shall need empirical measures to represent the boundary terms as well. 
With a slight abuse of notation, for fixed $\omega, \tilde \omega,$ consider the empirical measures, 
\begin{equation}\label{emp_meas-d4}
	\mu _{N, \bm{X} (\omega, \tilde \omega)} =   \frac 1 N \sum _{i} \   \delta _{X_i(\omega)} \, ,
\end{equation}
and
\begin{equation}\label{emp_meas-d5}
	\mu _{S, N^b, \bm{X} (\omega, \tilde \omega)} =   \frac 1 {N _b} \sum _{i} \   \delta _{X^b_i(\tilde \omega)} \, . 
\end{equation}
Then, it is useful to observe that the discrete sums  can be written as
\begin{equation}\label{energy_emp_meas-d4}
\frac 1N  \sum _{i}   \Big  | \, L\, g(X_i(\omega)\, ) - f(X_i (\omega) \, ) \Big  | ^2  = \int _\Omega | L\, g(x)-f(x)|^2\, \d \mu _{N, \bm{X} (\omega)} (x) \, .
\end{equation}
Similarly, 
\begin{equation}\label{energy_emp_meas_b}
\frac 1{N _b} \sum _{j}   \Big  | g(X^b_j(\tilde \omega)\, ) \Big  | ^2  = \int _{\partial \Omega}  |  g  (x)  |^2\, \d \mu _{S, N^b, \bm{X} (\omega, \tilde \omega)} (x) \, .
\end{equation}
Notice that we have, 
\begin{equation}\label{energy_emp_meas0}
	\begin{split}
	  \frac 1N  \sum _{i=1}^N  \EE \Big[  \,   \big  | L g(X_i \, ) - f(X_i   \, ) \big  | ^2 \Big] &=\|L g-f\|_{L^2(\Omega )}^2\\
	  \frac{1}{N_b}  \sum _{i=1}^{N_b}  \EE \Big[  \,   \big  | g(X_i^b \, ) \big  | ^2 \Big] &=\|g\|_{L^2(\partial \Omega )}^2 .
	  \end{split}
\end{equation}\\
{These identities require $\mu$ and $\mu_S$ to be the normalised Lebesgue and surface measures on $\Omega$ and $\partial\Omega$, which we assume in this section.}

{The following theorem states, in essence, that for every training set, the residual and the boundary values of the computed PINN solution are bounded by the computed loss plus one-dimensional Wasserstein terms of a posteriori nature, with no assumption on the network space.}
\begin{theorem}[\sc A posteriori error estimate for each instance]
\label{Thm:apost-PINN}
Consider for each $\omega \in \varOmega, \tilde \omega \in \tilde \varOmega$ the problem 
\begin{equation}\label{prob_E_def-d3}
	\min  _ {v \in  V _{\mathcal{N}} } \mathcal{E}_{N, \omega, {\tilde \omega}} (v)\, .
\end{equation}
 We denote by   $ u_{\bm{X} (\omega, {\tilde \omega} )   }  = u ^\star(\omega, {\tilde \omega} , \cdot)  \,  $ the  $V_{\mathcal{N}}-$\emph{ probabilistic PINN solution  of $Lu = f$ at $ \bm{X}  = (\, (X_1, \dots, X_N  ), (X^b_1, \dots, X^b_{N_b})\,   )  \, .$}
Then for each $\omega \in \bbOmega, \tilde \omega \in \tilde \bbOmega$  
\begin{equation}\label {W-dist-d_apost_1}
\begin{split}
 \Big \{ \  \int _{\Omega } \big |&\big [\, L \, u_{\bm {X}(\omega, \tilde \omega)}-\, L \, u \, \big ](y)\big|^2\, \d \mu(y) \Big \} ^{1/2}
 \leq \Big \{  \, \mathcal{E}_{N, \omega, {\tilde \omega}} (   u_{\bm {X}(\omega, \tilde \omega)} ) \Big \}^{1/2}\\
 &+  W_2(\big [\, L \, u_{\bm {X}(\omega, \tilde \omega)}-f\, \big ]_\sharp\mu, \big [\, L \,  u_{\bm {X}(\omega, \tilde \omega)}-f\, \big ]_\sharp\mu _{N, \bm{X} (\omega, \tilde \omega)})\, .\\
\end{split}
\end{equation}
In addition, 
\begin{equation}\label {W-dist-d_apost_1_bound}
\begin{split}
 \Big \{ \  \int _{\partial \Omega } \big |   u_{\bm {X}(\omega, \tilde \omega)} (s)\big|^2\, &\d \mu_S(s) \Big \} ^{1/2}
 \leq \Big \{  \, \mathcal{E}_{N, \omega, {\tilde \omega}} (   u_{\bm {X}(\omega, \tilde \omega)} ) \Big \}^{1/2}\\
 &+  W_2(\big [ u_{\bm {X}(\omega, \tilde \omega)}  \big ]_\sharp\mu_S, \big [ \,  u_{\bm {X}(\omega, \tilde \omega)}  \big ]_\sharp\mu _{S, N^b, \bm{X} (\omega, \tilde \omega)})\, .\\
\end{split}
\end{equation}
%
%
\end{theorem}

The proof is a simple adaptation of the analysis of the previous section. Notice that 
\begin{equation}\label {W-dist-d_0_2}
\begin{split}
\Big \{ \  \int _{\Omega } \big |\big [\, L \, u_{\bm {X}(\omega, \tilde \omega)}-f\, \big ](y)\big|^2\, \d \mu(y) \Big \} ^{1/2} =	W_2(\big [\, L \, u_{\bm {X}(\omega, \tilde \omega)}-f\, \big ]_\sharp\mu, \delta_0) \, .
\end{split}
\end{equation}
Then, since $W_2$ is a distance, 
\begin{equation*}\label {W-dist-d_0_3-d2}
\begin{split}
\Big \{ \  \int _{\Omega } \big |\big [\, L \, u_{\bm {X}(\omega, \tilde \omega)}-f\, \big ](y)\big|^2\, &\d \mu(y) \Big \} ^{1/2}
 \leq	 W_2(\big [\, L \, u_{\bm {X}(\omega, \tilde \omega)} -f\, \big ]_\sharp\mu _{N, \bm{X} (\omega)}, \delta_0)\\
  &+  W_2(\big [\, L \, u_{\bm {X}(\omega, \tilde \omega)}-f\, \big ]_\sharp\mu, \big [\, L \, u_{\bm {X}(\omega, \tilde \omega)}-f\, \big ]_\sharp\mu _{N, \bm{X} (\omega, \tilde \omega)})\\
 = &\Big \{  \frac 1N  \sum _{i}   \Big  | \, L \, u_{\bm {X}(\omega, \tilde \omega)} (X_i(\omega)\, ) - f(X_i (\omega) \, ) \Big  | ^2  \Big \}^{1/2}\\
 &+   W_2(\big [\, L \, u_{\bm {X}(\omega, \tilde \omega)}-f\, \big ]_\sharp\mu, \big [\, L \, u_{\bm {X}(\omega, \tilde \omega)}-f\, \big ]_\sharp\mu _{N, \bm{X} (\omega, \tilde \omega)})\\
 \leq & \Big \{  \, \mathcal{E}_{N, \omega, {\tilde \omega}} (   u_{\bm {X}(\omega, \tilde \omega)} ) \Big \}^{1/2}\\
 &+  W_2(\big [\, L \, u_{\bm {X}(\omega, \tilde \omega)}-f\, \big ]_\sharp\mu, \big [\, L \,  u_{\bm {X}(\omega, \tilde \omega)}-f\, \big ]_\sharp\mu _{N, \bm{X}(\omega, \tilde \omega)})\, .\\
  \end{split}
\end{equation*}
The boundary terms are handled similarly,  
\begin{equation}\label {W-dist-d_apost_1_bound-d2}
\begin{split}
 \Big \{ \  \int _{\partial \Omega } \big |   u_{\bm {X}(\omega, \tilde \omega)} &(s)\big|^2\, \d \mu_S(s) \Big \} ^{1/2} =  
 W_2(\big [ u_{\bm {X}(\omega, \tilde \omega)}  \big ]_\sharp\mu_S, \delta_0) \\
 \leq &\Big \{  \frac 1{N _b} \sum _{j}   \Big  |  u_{\bm {X}(\omega, \tilde \omega)}(X^b_j(\tilde \omega)\, ) \Big  | ^2  \Big \}^{1/2}\\
 &+  W_2(\big [ u_{\bm {X}(\omega, \tilde \omega)}  \big ]_\sharp\mu_S, \big [ \,  u_{\bm {X}(\omega, \tilde \omega)}  \big ]_\sharp\mu _{S, N^b, \bm{X} (\omega, \tilde \omega)})\, .\\
 \leq &\Big \{  \, \mathcal{E}_{N, \omega, {\tilde \omega}} (   u_{\bm {X}(\omega, \tilde \omega)} ) \Big \}^{1/2}\\
 &+  W_2(\big [ u_{\bm {X}(\omega, \tilde \omega)}  \big ]_\sharp\mu_S, \big [ \,  u_{\bm {X}(\omega, \tilde \omega)}  \big ]_\sharp\mu _{S, N^b, \bm{X} (\omega, \tilde \omega)})\, .\\
\end{split}
\end{equation}
The proof is  complete.\\

The theorem bounds the residual and the boundary values of $u_{\bm X(\omega,\tilde\omega)}$, not its error. The step from the one to the other is the stability of the continuous problem, as in Section~\ref{SubS:aPaP:Error_analysis}: for the very weak formulation of \eqref{ellipticPDE} on a smooth domain, see \cite{Lions_Magenes_I}, one has
\begin{equation}\label{very_weak_stability}
\|v\|_{L^2(\Omega)}\le C\big(\|Lv\|_{L^2(\Omega)}+\|v\|_{L^2(\partial\Omega)}\big)\qquad\text{for all } v\in H^2(\Omega),
\end{equation}
so that, applied to $v=u_{\bm X(\omega,\tilde\omega)}-u$, the two estimates of the theorem yield a computable a posteriori bound for $\|u_{\bm X(\omega,\tilde\omega)}-u\|_{L^2(\Omega)}$ in terms of the loss and the two one-dimensional Wasserstein terms. Notice that the regularity estimate \eqref{very_weak_stability} is an indicative bound, and in fact, if the $L^2(\partial\Omega)$ norm is used in the loss, the terms on the right hand side of this bound are not balanced  always. The way the boundary data enters the loss is a delicate issue not discussed here in detail, \cite{GGM_IMA_2025}.

	 In    \cite{Mishra_gen_err_pinn:2023}   the  cumulative error, i.e.,  an average over training sets, of generalisation error is estimated by quantities involving integrals over $D^N$ where $N$ is the number of samples. 
	 Probabilistic training is considered as well in \cite {Karniadakis:pinn:conv:2019, shin2020error, hong2022priori} using in the analysis, as mentioned above, 
	 Rademacher complexity. 
	 


\subsection{Estimates for the generation of measures and for operator learning}\label{SubS:Rel:ext}

\subsubsection*{Operator learning} The estimates of this section extend to the
operator learning problem of Section~\ref{SubS:ModP:function-learning}, and
it is instructive to state them in that setting, since the comparison with
the function learning case shows exactly where the dimension of the input
space enters and where it does not. Let $X$ be a separable metric space of
inputs with a probability measure $\mu\in\mathcal{P}_p(X)$, $Y$ a separable
Hilbert space, $\mathcal{G}:X\to Y$ measurable with
$\int_X\|\mathcal{G}\|_Y^p\,\d\mu<\infty$, and for i.i.d.\ samples
$a_i=X_i(\omega)$ of $\mu$ let $G_{\bm X(\omega)}\in V_{\mathcal{N}}$ be a
minimiser of
\begin{equation*}
\mathcal{E}_{N,\omega}(G)=\frac1N\sum_{i=1}^N\big\|G(X_i(\omega))-\mathcal{G}(X_i(\omega))\big\|_Y^p
=\int_X\|G(a)-\mathcal{G}(a)\|_Y^p\,\d\mu_{N,\bm X(\omega)}(a)
\end{equation*}
over a class $V_{\mathcal{N}}$ of neural operators, where
$\mu_{N,\bm X(\omega)}$ is the empirical measure of the inputs. The error is
measured in $L^p(X,\mu;Y)$, and the key observation is that both the error
and the loss depend on $G_{\bm X(\omega)}$ only through the scalar function
\begin{equation*}
\xi_\omega(a):=\big\|G_{\bm X(\omega)}(a)-\mathcal{G}(a)\big\|_Y,\qquad a\in X,
\end{equation*}
since $\|G_{\bm X(\omega)}-\mathcal{G}\|_{L^p(X,\mu;Y)}=\|\xi_\omega\|_{L^p(X,\mu)}$ and
$\mathcal{E}_{N,\omega}(G_{\bm X(\omega)})=\|\xi_\omega\|^p_{L^p(X,\mu_{N,\bm X(\omega)})}$.
The triangle inequality in $L^p(\pi)$ used above applies verbatim to
$\xi_\omega$, with $X$ in place of $\Omega$, and the push-forwards
$(\xi_\omega)_\sharp\mu$ and $(\xi_\omega)_\sharp\mu_{N,\bm X(\omega)}$ are
measures on $\R$ whatever $X$ and $Y$ are. Hence, for each $\omega$ and each
$\Phi\in V_{\mathcal{N}}$,
\begin{equation}\label{OL_apost}
\big\|G_{\bm X(\omega)}-\mathcal{G}\big\|_{L^p(X,\mu;Y)}
\le \mathcal{E}_{N,\omega}(\Phi)^{1/p}
+W_p\big((\xi_\omega)_\sharp\mu,(\xi_\omega)_\sharp\mu_{N,\bm X(\omega)}\big)\, ,
\end{equation}
and, with the measurability of $\omega\mapsto G_{\bm X(\omega)}$ and Jensen's
inequality as before,
\begin{equation}\label{OL_expect}
\EE\big[\|G_{\bm X(\cdot)}-\mathcal{G}\|_{L^p(X,\mu;Y)}\big]
\le \inf_{\Phi\in V_{\mathcal{N}}}\|\Phi-\mathcal{G}\|_{L^p(X,\mu;Y)}
+\EE\big[W_p\big((\xi_\cdot)_\sharp\mu,(\xi_\cdot)_\sharp\mu_{N,\bm X(\cdot)}\big)\big]\, .
\end{equation}
The validation estimate \eqref{validation_estimate} holds in the same form,
with the validation loss\\
$\frac1M\sum_j\|G_{\bm X}(Y_j)-\mathcal{G}(Y_j)\|_Y^p$
and a one-dimensional Wasserstein term for which \eqref{emp_W_rates_d=1}
applies conditionally on the training sample.

If, in addition, $\xi_\omega$ is Lipschitz continuous on $X$ with constant
$L_{\bm X(\omega)}$, which holds when $G_{\bm X(\omega)}$ and $\mathcal{G}$ are
Lipschitz from $X$ to $Y$, then \eqref{Lip_pushforward} gives
\begin{equation}\label{OL_Lip}
W_p\big((\xi_\omega)_\sharp\mu,(\xi_\omega)_\sharp\mu_{N,\bm X(\omega)}\big)
\le L_{\bm X(\omega)}\,W_p\big(\mu,\mu_{N,\bm X(\omega)}\big)\, ,
\end{equation}
and, under a uniform bound $L_{\bm X}\le L_N$ as in \eqref{W-dist-ass},
\begin{equation}\label{OL_apriori}
\EE\big[\|G_{\bm X(\cdot)}-\mathcal{G}\|_{L^p(X,\mu;Y)}\big]
\le \inf_{\Phi\in V_{\mathcal{N}}}\|\Phi-\mathcal{G}\|_{L^p(X,\mu;Y)}
+L_N\,\EE\big[W_p(\mu,\mu_{N,\bm X(\cdot)})\big]\, ,
\end{equation}
where now $W_p(\mu,\mu_{N,\bm X})$ is the distance of the empirical measure of
the inputs from $\mu$ on the space $X$.

The comparison with the function learning case is as follows. The
estimates \eqref{OL_apost} and \eqref{OL_expect}, and the validation
estimate, are parallel in form to those of
Section~\ref{SubS:Rel:MCInterpolation}. As there, the reduction to the
\emph{scalar} $\xi_\omega$ is an advantage of the approach, since the
Wasserstein terms compare measures on the real line, whatever the spaces
$X$ and $Y$ are. The output space $Y$ enters the analysis only in
Section~\ref{SubS:Conv:MCInterpolation}, where the $Y$-valued push-forward
and the Hilbert-space setting of \cite[Chapter 5]{AGS:2005} are needed to
identify the limit function.

The input space enters only through \eqref{OL_Lip}--\eqref{OL_apriori}, in
two ways. First, the Lipschitz continuity of $\xi_\omega$ is relative to the
metric of $X$, and it is well known that the choice of $X$ and $Y$
influences the continuity of the solution operators of PDEs, and of the
neural operators that approximate them. Second, the rate of
$\EE[W_p(\mu,\mu_{N,\bm X})]$ on $X$ is no longer given by
\eqref{emp_W_rates}, which has no analogue in an infinite-dimensional
space; it is governed by the intrinsic dimension of the input measure $\mu$
in the sense of \cite{Weed_Bach:2019}. For parametric inputs, which depend
on finitely many parameters, this dimension is finite and
\eqref{OL_apriori} yields corresponding rates. For infinite-dimensional
input measures no rate is available in general, and the Lipschitz route
provides convergence without a rate. In contrast, the a posteriori
estimates \eqref{OL_apost} and \eqref{validation_estimate} retain their
one-dimensional character in all cases, which is where the present approach
differs most from the error analyses of operator learning based on covering
numbers of the operator class, see Section~\ref{SubS:Conv:bib}: those
require uniform Lipschitz bounds on the networks and encode the input
discretisation in the complexity of the class, whereas here the only
network-dependent quantity is the computed $L_{\bm X(\omega)}$, and it
enters only the a priori form \eqref{OL_apriori}.

\subsubsection*{Generation of measures} In the static Wasserstein case of Section~\ref{SubS:Conv:measures}, the estimates \eqref{Wp:Cea} and \eqref{Wp:apost} contain the term $W_p(\nu_{\, \text{data}},\nu_{\, \text{true}})$, the distance of a genuine empirical measure from its law on $\R^k$, to which \eqref{emp_W_rates} applies with the dimension $k$ of the data, or rather with the intrinsic dimension of $\nu_{\, \text{true}}$; this is where the dimension of the data enters the analysis of generative models. When the generator is sampled, the additional term $W_p(\widehat\mu_{\theta,M},\mu_\theta)$ in \eqref{Wp:apost:sampled} is bounded by \eqref{Lip_pushforward} applied to $T_\theta$,
\begin{equation*}
W_p\big((T_\theta)_\sharp\widehat\gamma_M,(T_\theta)_\sharp\gamma_0\big)\le \mathrm{Lip}(T_\theta)\,W_p(\widehat\gamma_M,\gamma_0)\, ,
\end{equation*}
which is \eqref{Wp:pushforward}; here $\widehat\gamma_M$ is the empirical measure of $M$ i.i.d.\ samples of the reference measure $\gamma_0$, so that \eqref{emp_W_rates} applies with the latent dimension $m$, and the dependence on the computed generator is confined to its Lipschitz constant, exactly as in the proposition above.

 \subsection{Bibliography and comments}\label{SubS:Rel:bib}
The results of this section are based on the works \cite{loulakis2023newap} and \cite{KLMTsiouvras}.

	Rademacher complexity methods are among the most widely used tools for estimating generalisation error. Related approaches include VC dimension, local Rademacher complexity, and kernel methods; see, for example, \cite{Oneto_local_Rademacher:2018}, \cite{Yousefi_local_Rademacher:2018}, \cite{Bartlett_VCdimension:2019}, \cite{Jegelka_gen_opt_transport:2021}, \cite{Bartlett_Benignoverfitting:2020}, \cite{WE_Rademacher:2020}, \cite{Montanari_neural_kernel:2021}, and \cite{chenVanden-Eijnden_feature:2022}, and the references therein.

Optimal transport is an important tool in data science and machine learning; see, for example, \cite{Peyre:CompOT:2019}, \cite{PanaretosWR:2019}, and the references therein. More recent developments include \cite{Jegelka_gen_opt_transport:2021}, \cite{Kats_Pant_2022optimizing}, and \cite{Peyre2023unbalanced} and applications to generative models.

In \cite{Jegelka_gen_opt_transport:2021}, the Rademacher complexity of classes of neural network functions is controlled using Wasserstein distances and their relation to correlation. As in many such approaches, the analysis assumes that the functions in the class satisfy certain structural properties. The resulting bounds  provide useful theoretical and empirical tools for understanding and controlling generalisation error.

The estimates of this section provide an alternative to the generalisation bounds of statistical learning theory, which are obtained for empirical risk minimisers over classes that are uniformly bounded and whose size is tied to the number of samples, see e.g.\ \cite{Chakraborty_Bartlett:2025} for recent rates in terms of the intrinsic dimension of the data: they are a posteriori in nature, hold for any computed approximation with bounded loss, and require no uniform structural assumption on the neural network space, the only network-dependent quantity being the Lipschitz constant of the computed approximation. The rates they yield through \eqref{Lip_pushforward} depend, by \cite{Weed_Bach:2019}, on the intrinsic dimension of the sampling measure rather than on that of the ambient space, so that for data-driven problems they improve in the same way as the rates of statistical learning theory, while requiring no structural assumption; whether dimension-free rates can be derived in this framework remains open. The one-dimensional empirical measures that appear through the push-forward are the subject of \cite{Bobkov_Ledoux:2019}. For the stability of the  weak formulation of elliptic problems used in \eqref{very_weak_stability} see \cite{Lions_Magenes_I}.
 \\

\newpage

\section {From mathematical analysis to algorithmic design}\label{S:design}

The framework and analytical findings developed in the previous sections not only provide a better understanding of existing algorithms, but also motivate the design of new methods with potentially improved properties. Several directions naturally emerge from this perspective. These include methods specifically designed to preserve important qualitative properties at the discrete level, training strategies that are more efficient and better adapted to the underlying problem, hybrid approaches that combine well-established numerical techniques such as finite elements with neural-network minimisation, and, more fundamentally, the exploitation of alternative variational principles as a basis for the design of new algorithms.

A central finding of our analysis is that the stability and efficiency of neural-network methods may be strongly influenced by the properties of the underlying continuous variational problem. In particular, strong coercivity of the parent continuous functional can provide a robust foundation for the corresponding neural-network minimisation problem. This observation suggests a natural strategy for algorithmic design: rather than starting from a given differential equation and applying a standard neural-network approximation, we seek alternative variational formulations of the problem that possess particularly favourable analytical properties, and then use these formulations to construct the corresponding discrete neural-network methods.

This viewpoint opens up a two-way interaction between mathematical analysis and algorithmic design. On the one hand, identifying strongly coercive or otherwise stable variational principles can lead to neural-network algorithms with improved stability, robustness, and potentially more efficient training. On the other hand, the search for such formulations raises new mathematical questions concerning the underlying differential equations and their equivalent variational structures. In particular, alternative variational principles may reveal analytical properties that are not apparent from the original formulation and may provide new tools for understanding the behaviour of both the continuous problem and its numerical approximations.

It is important to emphasise that neural-network minimisation provides a particularly flexible framework for exploiting this idea. Once a suitable variational principle has been identified, the same general minimisation mechanism can be adapted to a wide range of formulations and problems. In this sense, neural-network spaces and optimisation procedures can serve as a common computational framework within which different continuous variational principles can be discretised and compared. This flexibility allows sequences of neural-network minimisation problems to approximate a broad class of continuous problems while retaining, at least in principle, the stability properties inherited from the underlying variational formulation.

The resulting philosophy is therefore not simply to apply neural networks to existing numerical formulations, but to use the flexibility of neural-network minimisation as a tool for exploring and exploiting new variational structures. In the following sections, we pursue this approach by considering alternative variational principles and investigating their analytical properties, with particular emphasis on coercivity, stability, and their implications for the design of reliable numerical methods.
The four examples that follow illustrate the two observations of Section~\ref{sec:1}: Runge--Kutta PINNs and training with finite elements concern the second, the design of time and space discretisations that retain the coercivity of the continuous functional, while the variational principles for the Helmholtz equation and the deep Uzawa methods concern the first, the choice of a continuous variational principle with   strong  coercivity.

\subsection{Runge-Kutta PINNs}\label{SubS:RK-PINN}
 In this section, we consider time-dependent PDEs discretised by a class of Physics-Informed Neural Networks based on Runge--Kutta and related time-Galerkin methods. The motivation is not only to achieve higher-order accuracy, but also to control the \emph{qualitative behaviour} of the discrete solutions. Such structure-preserving properties are central in numerical analysis, with well-established applications to wave and Schrödinger equations, diffusion, geometric and structure-preserving methods, and the Navier--Stokes equations; see, e.g., \cite{Butcher_1964}, \cite{Crouzeix_1}, \cite{Crouzeix_2}, \cite{BUTCHER1996113}, \cite{Glow_theta1987}, \cite{HairerNW1993}, \cite{Hairer2002}, \cite{Thomee_book_2006}.
 
 Physics-Informed Neural Networks are based on minimising the $L^2$ norm of the residual of an evolution PDE over a prescribed neural network space. The resulting loss requires numerical quadrature of the space-time integrals, typically using Monte Carlo or Quasi-Monte Carlo methods, particularly in high dimensions. However, the effect of these discrete losses on the qualitative behaviour of the resulting approximations remains largely unexplored. Even for the scalar wave equation, it is unclear whether, and under what conditions, such methods preserve fundamental conservation properties.

We propose a training approach based on two principles: (a) \emph{the discrete loss is designed within a time-discrete framework}, and (b) \emph{the final loss incorporates a Runge--Kutta or time-Galerkin discretisation} in a suitably structured form. We show that the resulting methods inherit the stability properties of the underlying time-discrete schemes and exhibit computational behaviour consistent with them. Since the time discretisation affects only the temporal variable, the approach is particularly well suited to high-dimensional evolution problems. 


\subsubsection*{A linear model problem}
To motivate our approach we consider as a model problem a linear  parabolic equation, as before, 
%
\begin{equation}\label{ParabolicPDE-I}
\left\{
\begin{alignedat}{3}
&u_t + Lu = f  \qquad&&\text{in }  &&\varOmega_T, \\
& u =  0 \qquad&&\text{on }  &&\partial \varOmega \times (0,T] , \\
&  u =  u^0\qquad&&\textrm{in }  &&\varOmega ,
\end{alignedat}
\right.
\end{equation}
%
where $ f \in L^2( \varOmega_T) ,\; u^0 \in H^1_0(\varOmega) $, and $ L $ is a coercive,
self-adjoint, second order elliptic operator. 
The associated energies used will be the $L^2$-residuals  
\begin{equation}\label{Functional}
\mathcal{E}(v) = \int_0^T \int_{\varOmega} |   v_t + Lv  - f |^2 \d  x\, \d t +\,       |  v {(0)}  - u^0 |_ {H^1} ^2 
\end{equation}
defined over smooth enough functions and  domains  $\varOmega .$ The choice of the $H^1$ seminorm for the initial  
condition,  $  |  v  - u^0 |_ {H^1},$ is done in order to obtain a balanced energy which is particularly convenient 
in the analysis; obviously other choices are possible.  

As usual we consider neural network spaces as basis for our discretisation.  To avoid extra technical problems,   it will be useful to consider, following  \cite{sukumar2022exact}, the set of functions  
which exactly satisfy the boundary conditions through appropriate distance functions depending only on the domain 
$\varOmega  . $ If $\varPhi$ is such a function, we define 
\begin{equation}\label {V_N0}
	V _{\mathcal{N}, 0}= \{ u_\theta : \varOmega  \to \R ,  \ \text{where }  u_\theta (x, t) =  \varPhi (x)\mathcal{C}_L (x, t), \ \text{for some } \mathcal{C}_L\in  \mathcal{N} \,  \} .
\end{equation}
Then, $V _{\mathcal{N}, 0} \subset H^1((0,T) ; L^2 (\varOmega)) \cap L^2((0,T) ; H^2(\varOmega) \cap H^1_0(\varOmega)),$ for smooth enough activation function $\sigma.$

 \subsubsection*{Deterministic and probabilistic training} 
 
 Approximating  integrals 
using  probabilistic   quadrature rules are based on  
  a collection $ X_1,X_2,\ldots$ of i.i.d.\ $\varOmega _T$-valued random variables, defined on an appropriate probability space 
corresponding to sample points in $\varOmega_T\, .$
Let  $\omega $   be a fixed instance, and $X_i(\omega )\in \varOmega_T$ be the corresponding values of the random variables. 
Monte Carlo approximation of the space-time integral yields the discrete sum 
\begin{equation}\label{prob_E-d5}
\mathcal{E}_{N, \omega }( v )  =  \frac 1N  \sum _{i=1 } ^N    \,  | v_t (X_i(\omega) ) + L v (X_i(\omega) ) - f(X_i(\omega) )| ^2 .     
 \end{equation}
 The discrete minimisation problem for each instance is 
\begin{equation}\label{ieE-minimize_NNQ-d6}
	\min  _ {v \in  V _{\mathcal{N}, 0} } \E _{N, \omega }(v)\,  +  \E ^0 _{N, \omega }(v) ,  
\end{equation}
where $\E ^0 _{N, \omega }(v)$ is a Monte Carlo approximation of the initial condition. One of the main advantages 
of these discretisations is that they scale reasonably with the dimension, and are thus preferable for high-dimensional operators $L .$

  In contrast to the above approach, Runge--Kutta PINN algorithms are based on \emph{time-discrete training}, where the time variable is discretised first. {To fix notation, let, as in Section~\ref{SubS:Ex:parabolic}, $0=t^0<t^1<\cdots<t^N=T$ define a partition of $[0,T]$, and set $I_n:=(t^{n-1},t^n]$ and $k_n:=t^n-t^{n-1}$, $n=1,\dots,N$. We denote by $v^m(\cdot)$ and $f^m(\cdot)$ the values $v(\cdot,t^m)$ and $f(\cdot,t^m)$, respectively.}
 A simple time-discrete version of the functional \eqref{Functional} was introduced in Section \ref{SubS:Ex:parabolic}:
\begin{equation}\label{Functional_k}
\mathcal{G}_{k, IE} (v) = \sum _{n=1}^{N} \,  k_n \, \int_{\varOmega }\big | \frac{v^{n}-  v^{n-1}}{k_n} + L v^{n}  - f^{n} \big |^2 \, \, \d x 
+\,    |  v^0  - u^0 |_{H^1(\varOmega)} ^2 \, . 
\end{equation}
The analysis of Section \ref{SubS:PINN-Parabolic} suggests a deeper connection between standard stability notions for time discretisations and this class of neural network algorithms. A key step in establishing this connection is an alternative representation of the time-discrete functional in terms of suitable reconstructions: 
\begin{equation}\label{def_hatU-d3}
\widehat v(t):=\ell_0^n(t) v^{n-1}+ \ell_1^n(t) v^{n}, \quad t\in I_n,
\end{equation}
with
$\ell_0^n(t) := (t^{n}-t)/k_n$  and  $\ell_1^n(t) := (t-t^{n-1})/k_n.$ 
If  the piecewise constant interpolant of $v_{j}$  is denoted by $\overline v, $  
\begin{equation}\label{def_overU-d4}
\overline v(t):= v^{n}, \quad t\in I_n\, ,\end{equation}
the time discrete energy  $ \mathcal{G}_{IE, k}   $ becomes
\begin{equation}\label{ParabFunctional_k-d2}
\begin{split}
	\mathcal{G}_{IE, k} (v) &= \| \widehat v_t + L \overline v - \overline f \|^2_{L^2((0,T); L^2(\varOmega))}  + |   \widehat v ^0  - u^0 |_{H^1(\varOmega)}^2 \\
&= \int_0^T \|  \widehat v_t + L \overline v - \overline f \|^2_{L^2(\varOmega)}\, \d t  +   |\widehat v^0  - u^0 |_{H^1(\varOmega)}^2  .
\end{split}
\end{equation}
This representation of the loss   can motivate   the design of the Runge--Kutta PINN methods next.

\subsubsection*{Runge--Kutta Physics Informed Neural Networks} 
To introduce the Runge--Kutta PINNs, we first recall the connection of Runge--Kutta methods and collocation time discretisations. 

Recall first the definition of a Runge-Kutta method: For $q\in \mathbb{N},$ let $0\leqslant c_1<\dotsb<c_q\leqslant 1$ denote the intermediate nodes of a Runge--Kutta method or collocation nodes.
{With starting value $U^0=u^0,$ we consider the discretisation of problem \eqref{ParabolicPDE-I} by a $q$-stage Runge--Kutta method: we recursively define approximations $U^{n}\in V$ to the nodal values $u(t^{n})$, as well as internal approximations $U_{n i}\in V$ to the intermediate values $u(t_{n i}),$ $t_{n i}=t^{n-1}+c_ik_n,$ $n=1,\dots,N$, by}
%
\begin{equation}\label{eq:RK1-I}
\left \{
\begin{alignedat}{2}
&U_{ni}=U^{n-1}- k _n \sum_{j=1}^q a_{ij} \big (L U_{nj} -f(t_{nj})\big ),\quad &&i=1,\dotsc,q,\\
&U^{n} =U^{n-1}- k _n \sum_{i=1}^q b_i \big (L U_{ni} -f(t_{n i})\big ), \quad &&
\end{alignedat}
\right . 
\end{equation}
$n=1,\dotsc,N.$ Here, without being very precise, $V$ is a functional space where our approximations are sought for every $t\in [0, T];$ 
typically $V$ is included in the domain of the operator $L.$   
An important class of Runge--Kutta methods is given by the \emph{collocation Runge--Kutta} methods, see \cite{GuillouSoule_1969}, \cite{Wright_1970}, \cite{HairerNW1993}. These methods are equivalent to a collocation approximation $\widehat U$, which is a continuous piecewise polynomial of local degree $q$ satisfying $\widehat U(0)=u^0$ and the collocation conditions
\begin{equation}
\label{eq:coll-1-I}
\widehat U'(t_{ni})+L\widehat U(t_{ni})=f(t_{ni}), \quad i=1,\dotsc,q, \quad n=1,\dotsc,N.
\end{equation}

These methods admit a crucial pointwise formulation that will be the basis for the loss functional defined below. Indeed, following \cite{AMN2} and \cite{AMN3}, let $I_{q-1}$ denote the piecewise interpolation operator by polynomials of degree $q-1$ at the collocation nodes $t_{ni}$, $i=1,\dotsc,q$, $n=1,\dotsc,N$. Since $\widehat U'$ is itself piecewise polynomial of degree $q-1$, \eqref{eq:coll-1-I} can be written in the \emph{pointwise form}
\begin{equation}
\label{eq:coll-2-I}
\widehat U'(t)+I_{q-1} L\widehat U(t)=I_{q-1} f(t), \quad t\in I_n, \quad n=1,\dotsc,N.
\end{equation}
The interpolants $U:=I_{q-1}\widehat U$ and $I_{q-1}f$ are piecewise polynomials of degree $q-1$ and, in general, are discontinuous at the nodes $t^0,\dotsc,t^{N-1}$. This pointwise formulation of the numerical method will provide the basis for defining the Runge--Kutta discrete loss.

To this end, and without being very precise at this stage, we introduce appropriate interpolation operators associated with the chosen Runge--Kutta method, acting on sufficiently regular functions. We shall also introduce some notation related to the representation of piecewise polynomial functions:
\begin{equation}
\label{eq:Lagr6}
\widehat v (t)=\widehat I_{q} v (t)= \sum_{i=0}^q {\tilde \ell}_{ni}(t) v(\tilde t_{ni}),\quad 
\overline v (t) = I_{q-1}\widehat v (t)=\sum_{j=1}^q \ell_{nj}(t) \widehat v(t_{nj}),
\quad t\in I_n.
\end{equation}
Notice that, in principle, $\tilde t_{ni}$ need not coincide with $t_{ni}$. These nodes can be chosen for convenience, and in most important cases $\tilde t_{ni}$ can be chosen as an extension of the collocation points $t_{ni}$ by including one additional node.

We are ready to define the  Runge--Kutta discrete loss by 
\begin{equation}\label{RK_ParabFunctional}
\begin{split}
	\mathcal{G}_{RK, k} (v) &= \| \widehat v_t(t)+I_{q-1} L\widehat v(t)-I_{q-1} f(t) \|^2_{L^2((0,T); L^2(\varOmega))}  
	+ |   \widehat v ^0  - u^0 |_{H^1(\varOmega)}^2 \\
&= \int_0^T \| \widehat v_t (t)+I_{q-1} L\widehat v(t)-I_{q-1} f(t)\|^2_{L^2(\varOmega)}\, \d t  + |   \widehat v^0  - u^0 |_{H^1(\varOmega)}^2 .
\end{split}
\end{equation}
The Runge--Kutta Physics Informed Neural Network method is based on the discrete minimisation problem for the loss $\mathcal{G}_{RK, k} (v), $
\begin{equation}\label{ieE-minimize_NNQ-d7}
	\min  _ {v \in  V _{\mathcal{N}, 0} } \mathcal{G}_{RK, k} (v) . 
\end{equation}
The formulation of the method for nonlinear evolution problems is straightforward by using similar interpolation operators; see \cite{AMN2}.

\subsubsection*{Discrete variational stability}

The analysis in this case follows the main steps discussed in Section \ref{SubS:PINN-Parabolic} and Section \ref{S:Analysis-abstract}. 
  We begin by assuming that the neural network spaces satisfy appropriate approximability properties, consistent with established results in approximation theory; see, e.g., \cite{Grohs_space_time_approx:2022, Mishra:appr:rough:2022, Xu} and the references therein. The precise assumptions depend on several factors that we will not address here.
As discussed in Section \ref{S:Analysis-abstract}, we consider a sequence of neural network spaces $V_\ell$ with increasingly better approximation properties as the parameter $\ell$ increases. We shall use the following notation for the corresponding discrete energies:
\begin{equation}\label{GdeltaEnergies}
\mathcal{G}_{\ell}(v ) = 
\left\{
\begin{alignedat}{2} 
&  \mathcal{G}_{RK, k(\ell)} (v ) ,\quad && v  \in V_\ell, 
\\ 
&+ \infty ,\quad && \textrm{otherwise.}
\end{alignedat}
\right.
\end{equation}
Here $k=k(\ell)$ are selected just to satisfy  $k=k(\ell)\to 0 $
as $\ell\to \infty\, , $ and $\mathcal{G}_{RK, k(\ell)} (v ) =\mathcal{G}_{RK, k} (v )$ is defined by \eqref{RK_ParabFunctional}.

 The next proposition establishes that the functional $ \mathcal{G}_{RK, k}$   is equicoercive with respect to $ \widehat U , \overline U $. 

\begin{proposition}\label{TD:EquicoercivityofG} Assume that the following maximal regularity estimate 
is satisfied
\begin{equation}
\label{eq:MR_Se:3}
\begin{split}
\|\overline U\|_{L^2((0,T); H^2(\varOmega))} &+ \|\widehat U'\|_{L^2((0,T); L^2(\varOmega))} \\&\leqslant C\Big [ \| \widehat U(0)\|_{H^1(\varOmega)}
+\| {\widehat U' + I_{q-1} L \widehat U} \|_{L^2((0,T); L^2(\varOmega))}\Big ];
\end{split}
\end{equation} 
then, the functional $ \mathcal{G}_{ {RK,}k}$ defined in \eqref{RK_ParabFunctional} is stable with respect to $ \widehat U , \overline U  {= I_{q-1}\widehat U}$ in the following sense: 
\begin{equation}\label{GEquicoercivity_RK}
\begin{gathered}
\textrm{If} \;\; \mathcal{G}_{ {RK,}k}(U) \leqslant C \;\; \textrm{for some} \;\: C >0 \;,\; \textrm{we have} \\
\|\overline U\|_{L^2((0,T); H^2(\varOmega))} + \|\widehat U'\|_{L^2((0,T); L^2(\varOmega))} \leqslant C_{MR}.
\end{gathered}
\end{equation}
\end{proposition}
 
 Although the \emph{discrete maximal regularity} bound, \eqref{eq:MR_Se:3}, for $\mathcal{G}_{IE, k}$ in Section \ref{SubS:PINN-Parabolic} follows from a direct calculation based on the stability properties of the Backward Euler scheme, the assumption \eqref{eq:MR_Se:3} in the present setting can also be established through novel energy-based proofs of maximal regularity estimates previously obtained by Kov\'acs, Li, and Lubich \cite{kovacs2016stable}, now in the Hilbert space setting and allowing for variable time steps. For detailed proofs see \cite{akrivis2025rungekuttaphysicsinformedneural}.

 The  {$\liminf$ inequality} is proved following the steps of 
 Section \ref{SubS:PINN-Parabolic} adapted to the present more technical setting:  We assume there is a sequence denoted by $ U_\ell $, such that $ \mathcal{G}_{\ell} (U_\ell) \leq C $ uniformly in $ \ell .$ Then up to sub-sequences 
 we get  the existence of $u_{(1)}$ and $u_{(2)}$ such that 
\begin{equation*}\label{proofGlowersemicontinuTDiscr-d2}
\begin{gathered}
L \overline  U_\ell \rightharpoonup L u_{(1)}\;\: \textrm{and} \;\: \widehat U_\ell' \rightharpoonup u_{(2)} '\;\: \textrm{weakly in} \;\: L^2(0,T; L^2(\Omega))\, .
\end{gathered}
\end{equation*}
we can finally prove   that  $u_{(1)} = u_{(2)} = u$ and thus,  
\begin{equation*}\label{proofGlowersemicontinuityeq2-d3}
\int_{\Omega_T} | u' + L u -f |^2  \leq \liminf_{\ell \rightarrow \infty}  \int_{\Omega_T} |  \widehat U_\ell' + L \overline  U_\ell -  \overline  f|^2\, .
\end{equation*}


 \subsubsection*{$\limsup$ inequality and  convergence of minimisers}\label{TD:limsup-d2}
 Let $ w \in H^1(0,T; L^2(\Omega))\cap L^2(0,T; H^2(\Omega)  \cap  H^1_0(\Omega))    $ The aim is to show  the existence of a recovery sequence $(w_\ell)$ such that $w_\ell \to w$ and 
   $$ \mathcal{G}(w) = \lim_{\ell \rightarrow \infty} \mathcal{G}_{\ell} (w_\ell) .$$
  This can indeed be achieved by exploiting the approximation properties of the neural network spaces, combined with suitable smoothing arguments and detailed estimates. The technical details of the proof can be found in \cite{akrivis2025rungekuttaphysicsinformedneural}.
As a result we can finally prove, following the steps described in Sections \ref{S:Analysis-abstract}, \ref{SubS:PINN-Parabolic}, 

{Roughly speaking, the next result asserts that minimisers of the Runge--Kutta loss converge to the exact solution, under the maximal regularity of the time-discrete scheme and the approximability of the network spaces.}
\begin{theorem}[Convergence]\label{Thrm:TimeD} 
Let $ \mathcal{G},\; \mathcal{G}_{ \ell} $ be the energy functionals 
 {defined in \eqref{Functional} and \eqref{GdeltaEnergies}, respectively.}
Let $u$ be the exact solution of  {\eqref{ParabolicPDE-I}} and let $ (u_\ell) , $ $u_\ell \in V_{\ell}, $ be a sequence 
of minimisers of $  {\mathcal{G}_{\ell}}  $, i.e.,
\begin{equation}\label{minofTDEnergy-d2}
\mathcal{G}_{\ell} (u_\ell) = \inf_{v_\ell \in  {V_\ell}} \mathcal{G}_{ \ell} (v_\ell)\, .
\end{equation}
Then, 
\begin{equation}\label{TDConvOfDiscrMinE-d2}
\hat u_\ell \rightarrow u,  \;\; \; \textrm{in} \;\: L^2((0,T); H^1(\varOmega))  ,\;\ \  
\end{equation}
where $\hat u_\ell ={\widehat I}_{q}u_\ell \, .$  
%
%
\end{theorem}

 \subsection{Training with Finite Elements
}\label{SubS:FE-training}
\smallskip

  \subsubsection{The Deep Ritz method  
}
We consider the  boundary value problem   
\begin{equation}\label{ivp-AC}
\left \{
\begin{alignedat}{3}
&-\varDelta u  = f\quad &&\text{in}\,\,&&\varOmega \\[2pt]
&\quad  u =0 \quad &&\text{on}\,\,&&\partial\varOmega \\[2pt]
\end{alignedat}
\right .
\end{equation}
in a polygonal   domain
$\varOmega\subset\R^d, 1\leqslant d\leqslant 3$.
The natural   energy functional associated to this problem is 
\begin{equation}\label{en-func} 
\mathcal{E} (u)=\int_{\varOmega} \Big( \frac12 |\nabla u|^2  - f\, u \Big) \d x  \, ;
\end{equation}
then, the   solution of \eqref{ivp-AC} is the unique minimiser of the problem 
\begin{equation}\label{mm}
	\min  _ {v \in  H^1 _0 (\varOmega ) } \E (v)\, .
\end{equation}
The boundary conditions can be also imposed weakly within $\mathcal{E} .$ The deep Ritz method is based on the minimisation of the above functional on discrete neural network spaces appropriately discretised with collocation type methods to yield computable  approximations.  

The semidiscrete method  considers  minimisation of $ \E (v)$ on discrete neural network spaces without further discretisation of the functional. 
As we discussed $\E \,  $ should be further discretised. Applying  a Monte-Carlo integration 
will yield the \emph{fully discrete}  Deep Ritz method of  \cite{e2017deep}, see Section \ref{SubS:Ex:EllipticPDE}.  
%
%
%
 An alternative is {training through quadrature/collocation,} Section\ref{SubS:Ex:EllipticPDE} where we may use deterministic quadrature for integrals over $\varOmega.$  
%
 %
The deterministic quadrature   results in the discrete energy \eqref{E_DQh}, which necessitates evaluating $\frac 12 |\nabla g(z)| ^2  - f(z) g(z)$ at the quadrature points. This can be computationally accomplished using the back-propagation algorithm implemented in standard neural network software packages.
%
{Monte Carlo quadrature, on the other hand, approximates $\mathcal{E}(v)$ well for a sufficiently large number of samples $N$, scales reasonably with the dimension, and is the most popular training approach for neural network discretisations of PDEs; in low dimensions, however, it is computationally demanding.}

\subsubsection*{Training through Finite Element quadrature}
 Let $T_h$ be a shape regular triangulation of a polygonal domain $\Omega$ with mesh size $h=h(x),$ see  \cite{BrennerScott} for  standard finite element notation and terminology.
For $K\in T_h$, $K$ is an element of the triangulation, and  $\mathbb{P}_q(K)$ denotes the set of polynomials of degree less or equal to $q$. 
We define the standard space of continuous  piecewise polynomial functions as  
\begin{equation*}\label{Sh_def}
\tilde \S_h (\Omega):=\left\{v\in C^0(\bar \Omega)\,:\, {v}|_{K}  \in \mathbb{P}_q(K), K \in T_h \right\}.
\end{equation*}
We also consider the finite element space where zero boundary conditions are enforced
\begin{equation*}\label{Sh_def-d2}
 \S_h(\Omega):=\left\{v\in C_0^0(\bar \Omega)\,:\, {v}|_{K}  \in \mathbb{P}_q(K), K \in T_h \right\}.
\end{equation*}
  Then, define
  \begin{equation*}
\label{E_hFE}
\mathcal{E}_{\S _h}( g )  =\int_{\Omega} \Big( \frac12 |\nabla I_{\S _h} ( g)|^2  -  I_{\S _h} ( f\, g) \Big) \d x  \, . 
\end{equation*}
The discrete minimisation problem   is 
\begin{equation*}\label{ieE-minimize_NNQ-d8}
	\min  _ {v \in  V _{\mathcal{N}} } \E _{\S _h}(v)
\end{equation*}
The finite element spaces are used solely to define the discrete energy. Specifically, for $g\in V_{\mathcal{N}}$, we have $I_{\S_h}g\in\S_h$. Both $I_{\S_h}g$ and its gradient $\nabla I_{\S_h}g$ are piecewise polynomial functions determined entirely by the values of $g$ at the Lagrangian degrees of freedom. Consequently, the corresponding integrals can be evaluated using standard finite element tools and software, which is a key advantage of the method.
In particular, the computation of
$$
\int_{\varOmega}|\nabla I_{\S_h}(g)|^2dx
$$
does not require backpropagation. Backpropagation is needed only to compute the derivative of $\mathcal{E}_{\S_h}$ with respect to the neural-network parameters, which is required for the iterative computation of its minimisers.
 
The point of this choice, in the terms of these notes, is twofold. First,
$\mathcal{E}_{\S_h}(v)$ is the \emph{exact} energy of the finite element
interpolant of $v$: it is coercive on $\S_h\subset H^1_0(\Omega)$ and
controls $I_{\S_h}v$ between the quadrature points, which is precisely what
the sampled energies of Sections~\ref{SubS:Ex:EllipticPDE}
and~\ref{SubS:Conv:MCInterpolation} fail to do; training with finite
elements is thus a training step designed to retain the coercivity of the
continuous functional, in the sense of discrete variational stability.
Second, the discrete energy is evaluated with the standard tools of finite
elements, which is both efficient, since no backpropagation is needed for
the energy itself, and an opening to the whole finite element toolbox:
boundary conditions can be imposed in various ways, even exactly through $\S_h$, the mesh can be
adapted, and the a posteriori error control and adaptivity of finite
element theory can prove useful to design and control the neural network method.

\subsubsection* {Discrete variational stability}
 
We consider a sequence of problems as follows: 
  For each   $\ell$ we associate a finite element space $\S _{h(\ell)}, $ with maximum diameter $h(\ell)$ such that 
$h(\ell)  \to 0,  \  \ \ell\to \infty\, .$ The corresponding   minimisation problems are
%
\begin{equation}\label{E_ell}
	\min_ {v \in  V _{\ell }} \mathcal{E}_{\ell}(v), \qquad \text{where } \mathcal{E}_{\ell}(v) := \mathcal{E}_{\S _{h(\ell)} }(v)\, .
\end{equation}
The corresponding $\S _{h(\ell)} $-finite element-deep Ritz minimisers are denoted by $u_\ell. $ 
%
To verify the stability of the method we observe that  
 the functional $ \mathcal{E} _\ell $   is stable with respect to the $ H^1 $-norm, in the following sense: Let $ (v_\ell) $ be a sequence of functions in $ V_\ell $ such that for a constant $ C>0 $ independent of $ \ell $, it holds that
\begin{equation*}\label{EquicoercivityofEdelta1(2)}
\mathcal{E}_\ell (v_\ell) \leq C.
\end{equation*} 
Then there exists a constant $ C_1>0 $ such that
\begin{equation*}\label{EquicoercivityofEdelta2(2)}
\| I_{\S _{h(\ell)} }v_\ell \|_{H^1(\Omega)} \leq C_1\, .
\end{equation*}
Furthermore, by   convexity, 
the  weak lower semicontinuity follows:
 $$\liminf\limits_{\ell\rightarrow \infty}
\int_\Omega |\nabla I_{\S _{h(\ell)} }v_\ell |^2 \ge \int_\Omega |\nabla     v|^2.$$

 \subsubsection*{Convergence}

\smallskip

{Roughly speaking, the next result essentially states that the finite element interpolants of the Deep Ritz minimisers trained with finite element quadrature converge to the exact solution, provided the mesh and the network spaces are refined together in a compatible way.}
\begin{theorem}[Convergence of the discrete minimisers]\label{Thrm:Gamma_funct} Let $ \mathcal{E},\; \mathcal{E}_{\ell} $ be the energy functionals defined in \eqref{en-func} and 
\eqref{E_ell}
 respectively, and $f\in C^0(\bar \Omega).$
 Let $ (u_\ell) , $ $u_\ell \in V_{\ell}, $ be a sequence of minimisers of $ \mathcal{E}_{\ell} $ and $$\widehat  u_\ell: =I_{\S _{h(\ell)} } u_\ell \, .$$
{Assume that the spaces $V_\ell$ have the approximability property \eqref{w_ell_7_hSm} with $[m,s,p]=[2,0,\infty]$, and that the finite element spaces are chosen such that $ \underline {h}_{E, \ell} ^ {-1/2}\,  (\tilde  \beta ^{[2,0,  \infty]} _\ell) ^{1-2\epsilon} \leq C$ for some $\epsilon\in(0,1/2)$, where $E_{h(\ell)}$ is the set of edges (faces in three dimensions) of $T_{h(\ell)}$, $h_e$ the diameter of $e\in E_{h(\ell)}$, and $\underline {h}_{E, \ell} = \min_ {e\in E_{h(\ell)}} h_e.$ Then we have }
\begin{equation*}\label{ConvOfDiscrMinE}
\widehat  u_\ell \rightarrow u ,  \;\; \; \textrm{in} \;\: L^2 (\Omega), 
\quad 
\widehat   u_\ell \rightharpoonup u  \, ,  \;\; \; \textrm{in} \;\: H^1 (\Omega),
\qquad \ell \to \infty\, .
\end{equation*}
where $u $ is the exact solution of the  problem.
\end{theorem}

{The approximability property of the spaces $V_\ell$ assumed in the theorem is the following: for each $w\in W^{m,p}(\Omega)\cap H^1_0(\Omega)$ there exist $w_\ell\in V_\ell$ such that}
%
 \begin{equation*}\label{w_ell_7_hSm}\begin{split} 
 \|w_{\ell}-w\|_{W^{s, p} (\Omega)}  \leq \ \tilde  \beta ^{[m, s,  p]} _\ell \,  | w | _ {W^{m, p} (\Omega)} , \qquad 
 \text{and } \ \tilde \beta _\ell ^{[m, s,  p]} \, \to 0,  \  \ \ell\to \infty\, .	\end{split} \end{equation*}
Then, since  $ \widehat v_\ell  \rightarrow    v$ in $L^2(\Omega)$  we can show   that    
$\lim _{\ell\to  \infty}
\int_\Omega  I_{\S _{h(\ell)} } ( v_\ell  f )\d x  = \int_ \Omega     v f \d x . $
Furthermore for any $ w \in H^1_0 (\Omega)  $  we can show the existence of a recovery sequence  $(w_\ell)$, such that  $ \mathcal{E}(w) = \lim_{\ell \rightarrow \infty} \mathcal{E}_\ell (w_\ell)  .$  
 Details are  rather technical and the last step requires smooth activation function.   
The condition of the theorem couples the smallest edge of the mesh with the rate at which the network spaces approximate, in $L^\infty$, functions bounded in $W^{2,\infty}$: the finite element interpolant $\widehat w_\ell$ of the network function $w_\ell$ is controlled through inverse estimates on the edges, and the condition ensures that $\widehat w_\ell\to w$ in $H^1(\Omega)$, from which the limsup inequality follows by the continuity of the energy. The liminf inequality and the compactness follow from the stability of $\mathcal{E}_\ell$ in the $H^1$-norm established above. In addition, $\mathcal{E}_\ell(u_\ell)\to\mathcal{E}(u)$, as the liminf--limsup argument gives.

{For ReLU activation functions the network functions are not in $W^{2,\infty}(\Omega)$, and the last step of the proof requires a different argument, based on a Cl\'ement-type interpolant applied to a mollified recovery sequence; we refer to \cite{grekas2024deepritzfinite} for the details.}

\newpage
\subsection{Variational Principles for the Helmholtz equation}\label{SubS:Helmholtz}
\begingroup 
\def\dimx{\nu}
 {Consider the Helmholtz equation with an impedance 
boundary condition }
\begin{equation}\label{Helm_PDE}
  -\Delta u - k^2 u = f\quad \text{ in }  \varOmega 
  \qquad \partial_{\mathbf{n}} u = i k  u \quad \text{ on }  \partial\varOmega\, .
\end{equation}
This is the time-harmonic analogue of the wave equation with first-order absorbing boundary conditions. The function $u$ is complex-valued, $u:\varOmega\to\mathbb{C}$, with $\varOmega\subset\mathbb{R}^{\dimx}$ a bounded Lipschitz domain, $f\in L^2(\varOmega)$ complex-valued, and $k>0$ the wavenumber. We use the compact notation $\L u:=-\Delta u-k^2u$. Absorbing boundary conditions allow wave propagation problems on large or unbounded domains to be reduced to computationally tractable bounded domains. Despite its classical nature, accurately approximating this model at high wavenumbers remains a significant numerical challenge.

In this section, we investigate whether \emph{variational principles} can be associated with the Helmholtz equation, and in particular with \eqref{Helm_PDE}. We design new strongly coercive variational principles while 
our derivations are guided throughout by the aim of preserving the underlying physical structure of the problem.
Beyond their theoretical interest, these variational principles have important implications for computational methods. Neural-network-based approaches typically rely on discrete minimisation problems based either on continuous energies or on the $L^2$ residual of the differential operator \cite{Karniadakis:pinn:orig:2019}. For highly oscillatory problems such as \eqref{Helm_PDE}, variational formulations with strong coercivity offer the prospect of more robust and reliable discretisations. 

We show that the proposed principle leads to an effective neural-network formulation of \eqref{Helm_PDE}. Strongly coercive bilinear forms for the Helmholtz equation were previously developed in \cite{Moiola_Spence2014}, leading in general to finite element methods with non-symmetric matrices. Although we employ similar analytical tools, notably the Rellich and Morawetz identities, our motivation and design strategy are substantially different and prioritise the derivation of novel variational principles.

\def\div{\text{\rm div}}
\def\union{\cup} 
\def\cal \mathcal{ } 
\def\ah{A_h} 
\def\bh{B_h} 
\def\nv{\|_ V} 
\def\nstar{\|_*} 
\def\nst2{\| _*} 
\def\ninfty{\| _ \mathcal{ B} } 
\def\a12{A_h ^{1/2} } 
\def\d{{\mathrm d}}
\def\tr|{|\! |\! |} 
\def\vp{\varpi }

\def\i{\text {\rm i} }
\def\e{\text {\rm e} }
\def\C{{\mathbb C}}
\def\R {{\mathbb R}}
\def\P{{\mathbb P}}
\def\N{{\mathbb N}}
\def\Z{{\mathbb Z}}
\def\N{{\mathbb N}}
\def\E{{\mathcal{E}}}
\def\A{{\mathcal{A}}}
\def\B{{\mathcal{B}}}
\def\D{\mathscr{D}}
\def\L{{\mathscr{L}}}
\def\K{{\mathscr{K}}}
\def\AA{\text{${\mathcal O}\hskip-3pt\iota$}}
\def\rr{ {\text {\bf r} }}
\def\vv{ {\text {\bf v} }}
\def\uu{ {\text {\bf u} }}
\def\UU{ {\text {\bf U} }}

\def\wU{\widehat U}
\def\wR{\widehat R}
\def\he{\hat e}

\def \x{\overline x}
\def \l{\ell }
\def \f{\varphi }
\def \th{\vartheta }
\def \r{\rho }
\def \s{\sigma }
\def \z{\zeta }
\def \a{\alpha }
\def \b{\beta }
\def \g{\gamma }
\def\wt{\widetilde} 
\def\co{\overline} 
\def\pp{ {\wt p} } 
\def\qq{ {\wt q} }
\def\cq{ {\co q} } 
\def\cqq{ \co {\wt q} } 
\def\cp{ \co {\wt p} } 
\def\cv{ {\co v} } 
\def\Dh{{\varDelta}_{h}^{n}} 
\def\Dhj{{\varDelta}_{h} ^J}
\def\Ik#1{\frac 1{k_ {#1} } }

\def\IL{{{I_{\hbox{\srm L}}}}}
\def\Th{{\varTheta}}
\def\T_h{{{\mathcal T}_h}}
\def\La{{\Lambda}}
\def\<{{\langle }}
\def\>{{\rangle }}
\def\V{{\mathcal{V}}}
\def\S{{\mathcal{S}}}
\def\U{{\mathcal{U}}}

\def\eps{\varepsilon}
\def\ds{\displaystyle}

\def\MR#1{\href{http://www.ams.org/mathscinet-getitem?mr=#1}{MR#1}}


 \newcommand{\norm}[1]{\ensuremath{\left|#1\right|}}
  \newcommand{\Norm}[1]{\ensuremath{\left\|#1\right\|}}
  
  \renewcommand{\vec}[1]{\boldsymbol{#1}}


\def\eps{\varepsilon}
\def\normH#1{\|#1\|_{\mathscr {H}}}
\def\normV#1{\|#1\|_{\mathscr {V}}}
\def\normVB#1{\|#1\|_{\mathscr {V}_{BC}}}

\def\Vhnod{{\stackrel{\!\scriptscriptstyle 0}{V}}_h}
\def\Hcnod{\stackrel{\,\scriptscriptstyle 0}{\Hc}}
\def\Hnod1{\stackrel{\!\scriptscriptstyle 0}{H^1}\!\!}
\def\Hnods{\stackrel{\scriptscriptstyle 0}{H^s}}
\def\Hdots{\dot{H}^s}
\def\Hdotone{\dot{H}^1}
\def\Hdotnod{\dot{H}^0}

\newcommand{\Dbar}{\mbox{$\overline{D}$}}
\newcommand{\Omegabar}{\mbox{$\overline{\Omega}$}}
\def\ubar{\overline{u}}
\def\vbar{\overline{v}}
\def\phibar{\overline{\phi}}
\def\kappat{\tilde{\kappa}}
\def\xt{\tilde{x}}
\def\Rt{\widetilde{R}}
\def\St{\widetilde{S}}
\def\Tt{\widetilde{T}}

\def\xn{x_{\sN}}
\def\xf{x_{\sF}}

\def\Ac{{\mathcal A}}
\def\Bc{{\mathcal B}}
\def\Cc{{\mathcal C}}
\def\Gc{{\mathcal G}}
\def\Hc{{\mathscr H}}
\def\Nc{{\mathcal N}}
\def\Qc{{\mathcal Q}}
\def\Rc{{\mathcal R}}
\def\Sc{{\mathcal S}}
\def\Tc{{\mathcal T}}
\def\Uc{{\mathcal U}}
\def\Vc{{\mathscr V}}
\def\Vcbc{{\mathscr V}_{BC}}
\def\Vctil{\widetilde {\mathscr V}}
\def\Hcnod{\stackrel{\,\scriptscriptstyle 0}{\Hc}}

\newcommand{\Pd}{\mathbb{P}}
\def\abold{\boldsymbol{\alpha}}
\def\ebold{\boldsymbol{e}}
\def\nubold{\boldsymbol{\nu}}
\def\fbold{\boldsymbol{f}}
\def\sbold{\boldsymbol{s}}
\def\ubold{\boldsymbol{u}}
\def\sF{\mbox{\tiny F}}
\newcommand{\sL}{\mbox{\tiny L}}
\def\sE{\mbox{\tiny E}}
\def\sN{\mbox{\tiny N}}
\newcommand{\sR}{\mbox{\tiny R}}
\newcommand{\cis}{\mbox{\scriptsize i}}
\newcommand{\smin}{\mbox{\scriptsize min}}
\newcommand{\smax}{\mbox{\scriptsize max}}
\def\astt{\mbox{\tiny{E}}}

\newcommand{\tu}{ \tilde{u} }
\newcommand{\tmbu}{ \tilde{\bm u} }
\newcommand{\hmbu}{ \hat{\bm u} }
\newcommand{\tV}{ \tilde{V} }
\newcommand{\mbr}{\bm r}
\newcommand{\mbR}{\bm R}
\newcommand{\mba}{\bm a}
\newcommand{\mbb}{\bm b}
\newcommand{\mbw}{\bm w}
\newcommand{\mbz}{\bm z}
\newcommand{\mbph}{\bm \phi}
\newcommand{\mbx}{\bm x}
\newcommand{\mbxi}{\bm \xi}
\newcommand{\mby}{\bm y}
\newcommand{\mbu}{\bm u }
\newcommand{\mbU}{\bm U }
\newcommand{\mbv}{\bm v}
\newcommand{\mbg}{\bm g }
\newcommand{\mbn}{\bm n}
\newcommand{\mbe}{\bm e }
\newcommand{\mbc}{\bm c }
\newcommand{\mbF}{\bm F }
\newcommand{\mbf}{\bm f }
\newcommand{\mbC}{\bm C }
\newcommand{\mbO}{\bm 0 }
\newcommand{\tW}{\tilde{W} }
\newcommand{\lm}{\tilde{\lambda}_m}
\newcommand{\lM}{\tilde{\lambda}_M}
\newcommand{\mbS}{\bm S }
\newcommand{\mbphi}{\bm \phi }
\newcommand{\mE}{\mathcal E}

\def\D{\mathscr{D}}
\def\L{{\mathscr{L}}}
\def\K{{\mathscr{K}}}

 \subsubsection*{Notation and preliminaries}

We introduce some standard notation. {Throughout, $\Omega\subset\mathbb{R}^{\dimx}$ is a bounded Lipschitz domain with outward unit normal $\mathbf{n},$ \cite{Grisvard_book}.}  {On the boundary we use the trace spaces $H^s(\partial\Omega)$, $s\in\{\tfrac32,\tfrac12,-\tfrac12,-\tfrac32\}$, and the trace operator $\tau :H^1(\Omega)\to H^{1/2}(\partial\Omega)$ \cite{Lions_Magenes_I}. When no confusion may arise we will use just $u$ to denote $ \tau u\, .$ }

 {We denote  by $\Hc$ the space  $H^1(\Omega)$ equipped with the wavenumber-dependent norm 
\[
\|v\|_{H^1_k(\Omega)} := \Bigl(\|\nabla v\|_\Omega^2 + k^2 \|v\|_\Omega^2\Bigr)^{1/2}\, .
\]}
{For energies that involve the bulk residual $ \L v:=-\Delta v-k^2 v $ in $L^2(\Omega)$ we use the  space}
\[
\Vc := \Bigl\{  v\in H^1(\Omega) :  \L v \in L^2(\Omega),\ v\in H^1(\partial \Omega),  \text{and}\ \partial_{\mathbf{n}} v - i k  v \in  L^2(\partial \Omega ) 
 \Bigr\},
\]
{endowed with the norm}
\[
\|v\|_{\Vc}^2 := \|\nabla v\|_\Omega^2 + k^2 \|v\|_\Omega^2 + \|\L v\|_\Omega^2  + k ^2 \|v \|_{\partial \Omega} ^2  + \| \nabla v \|_{\partial \Omega} ^2 + \| \partial_{\mathbf{n}} v - i k  v \|_{\partial \Omega} ^2 .
\]
{The homogeneous-impedance subspace is}
\[
\Vcbc := \Bigl\{  v\in \Vc  :  \partial_{\mathbf{n}} v - i k  v = 0 \ \text{in } L^2(\partial\Omega)  \Bigr\},
\]
and in the sequel we shall use the notation
$$
  \normVB{v}^2:=\|\nabla v\|_\Omega^2 + k^2 \|v\|_\Omega^2 + \|\L v\|_\Omega^2 + k ^2 \|v \|_{\partial \Omega} ^2 +  \| \nabla v \|_{\partial \Omega} ^2 .
$$

The standard variational formulation of the impedance problem is based on  the sesquilinear form $\Bc:$    
\begin{equation}   \label{eq:bf}
  \Bc(u,v) := \int_\Omega \nabla u\cdot \nabla \vbar \, \d x - k^2 \int_\Omega  u\vbar\, \d x
           - i\, k\int_{\partial \Omega} u \vbar\, \d S \, ,
\end{equation}
and is: 
\emph{Find $u\in \Hc$ such that}
\begin{equation}\label{eq:wf-gen}
\Bc(u,v) = (f,v)_\Omega   
\qquad \text{for all } v\in \Hc .
\end{equation}

\subsubsection*{Variational Principles}
The physical energy associated with the Helmholtz equation  {in unrestricted domains} can be formulated as 
\begin{equation}\label{Energy_H_main_s2}
\mathcal{E}_P(v) = \int_{\Omega }  \frac12  \Bigl(    |\nabla v |^2 - k^2|v|^2 \Bigr)   \d x  -     \Re \int_{\Omega } f   \overline{v}   \d  x   .
\end{equation}
Then the map  $\mathcal{E}_P  : \Hc \to \R $ is clearly differentiable. 
{We take variations with respect to the \emph{real} inner product on $H^1(\Omega;\mathbb{C})$, i.e. $\langle \phi,\psi\rangle_{\mathbb{R}}:=\Re \int_{\Omega}\phi \overline{\psi} \d x$, so that $D\mathcal{E}_P(u)[v]\in\mathbb{R}$ for all $u,v$.}
To find its derivative we first notice 
\begin{equation}\label{Energy_H_main_Re}
\mathcal{E}_P(v) =  \Re \Bigl\{  \int_{\Omega }  \frac12  \bigl( |\nabla v |^2 - k^2|v|^2 \bigr)  \d x  -  \int_{\Omega } f   \overline{v}   \d  x \Bigr\}. 
\end{equation}
Then, 
\begin{equation}\label{Energy_H_main_Re2}\begin{split}
\< D\mathcal{E}_P(u)  , v  \>   =&  \Re \Bigl\{  \int_{\Omega }  \tfrac12 \Bigl(   \nabla u\cdot \nabla \vbar  
+   \nabla v\cdot \nabla \ubar  - k^2( u \vbar + v \ubar ) \Bigr) \d x  -  \int_{\Omega } f   \overline{v}   \d  x \Bigr\}\\[2pt]
 = &  \Re \Bigl\{  \int_{\Omega }   \bigl( \nabla u\cdot \nabla \vbar  - k^2 u \vbar \bigr)   \d x  -  \int_{\Omega } f   \overline{v}   \d  x \Bigr\}.
\end{split}
\end{equation}
 Consider {  $u= u_R + i u_I $, $v= v_R + i v_I$, and $f = f_R + i f_I$}. 
Then stationary points of $ \mathcal{E}_P, $ i.e., $u\in \Hc $ satisfying,
\begin{equation}\label{Energy_H_main_stat_p}\begin{split}
\< D\mathcal{E}_P(u)  , v  \>  =   0, \quad \text{for all } \ v \in \Hc   , \end{split}
\end{equation}
are such that, for any $v\in \Hc , $
\begin{equation}\label{Energy_H_main_Re_system}
\begin{split}
  &    \int_{\Omega }   \Bigl(     \nabla u _R \cdot \nabla v_R 
 \ - k^2  u_R   v_R \Bigr) \d x  -  \int_{\Omega }{   f_R}    v_R   \d  x   =  0  , \\
 &    \int_{\Omega }   \Bigl(     \nabla u _I \cdot \nabla v_I 
 \ - k^2  u_I   v_I \Bigr) \d x  -   {    \int_{\Omega } f_I    v_I   \d  x }  =  0  . \\ 
\end{split}
\end{equation}
{Since $C_c^\infty(\Omega)\subset \Hc$ is dense, testing \eqref{Energy_H_main_stat_p} with $v\in C_c^\infty(\Omega)$ yields the Euler-Lagrange equation in distributional form.}
Considering test functions $v\in \Hc \cap C ^\infty _ c (\Omega)$ we conclude that stationary points 
of $ \mathcal{E}_P  $  satisfy $- \Delta u - k^2 u = f$ in $ \Omega,$ 
in the sense of distributions.
\begin{proposition}
 Consider the energy $\mathcal{E}_P$ defined in \eqref{Energy_H_main_s2}. The stationary points of 
 $\mathcal{E}_P$ defined on $\Hc$ satisfy the Helmholtz equation in the sense of distributions.
\end{proposition}
{However,  we cannot incorporate the absorbing boundary conditions into the energy when we restrict the computational domain. One way to view this is to split the domain as $\Omega = \Omega_{\text{Restr}} \cup \Omega_2$. The corresponding energies can then be written as follows, assuming that $f$ has compact support in $\Omega_{\text{Restr}}$:}
\begin{equation}\label{Energy_H_main_s2_2}
\begin{split}
\mathcal{E}_P(v) = &\int_{ \Omega _{\text{Restr}}  }  \frac12  \Bigl(    |\nabla v |^2 - k^2|v|^2 \Bigr)   \d x  -     \Re \int_{\Omega _{\text{Restr}}  } f   \overline{v}   \d  x  +  \int_{\Omega _{2} }  \frac12  \Bigl(    |\nabla v |^2 - k^2|v|^2 \Bigr)   \d x \\
=&: \mathcal{E}_{P, \Omega _{\text{Restr}} }(v) + \mathcal{E}_{P, \Omega _{2}}(v)
\end{split}\end{equation}
An energetically consistent approach requires representing the contribution of the complementary domain, $\mathcal{E}_{P, \Omega _{2}}(v)$,   through a surface integral incorporating appropriate boundary conditions.  However, it remains unclear whether impedance conditions can be used to represent this energy contribution. To preserve the variational structure associated with the boundary conditions in \eqref{Helm_PDE}, we therefore leverage complex Lagrangians, that is functionals where $v $ and $ \vbar$ are considered independent variables. 
To this end, we  define 
\begin{equation}\label{Energy_H_complex_1}
\begin{split}
\mathcal{E}_{P, \C } (v) = \mathcal{E}_{P, \C } (v,  \overline{v}) =  &\int_{ \Omega    }  \frac12  \Bigl(    |\nabla v |^2 - k^2|v|^2 \Bigr)   \d x  -     \Re \int_{\Omega   } f   \overline{v}   \d  x   - \frac{ik}{2}\int_{\partial \Omega}|v|^2 \d x  \\
 =  &\int_{ \Omega    }  \frac12  \Bigl(    \nabla v\cdot \nabla \vbar  - k^2v \vbar \Bigr)   \d x  -     \Re \int_{\Omega   } f   \overline{v}   \d  x   - \frac{ik}{2}\int_{\partial \Omega}v \vbar \d x  
\end{split}\end{equation}

We can then show the following

\begin{proposition}
 Consider the Lagrangian $\mathcal{E}_{P, \C }$ given by \eqref{Energy_H_complex_1} and defined on $\Hc$.
If $ \Bc(\cdot ,\cdot) $ is the  form defined in \eqref{eq:bf}, the stationary points of $\mathcal{E}_{P, \C }$ with respect to $\vbar,$  
 satisfy, 
 \begin{equation*}   \label{eq:bf_new}
  \Bc(u,v) = \int_{\Omega   } f   \overline{v}   \d  x, \quad \text{for all } v\in \Hc\, , 
\end{equation*}
 i.e., they are weak solutions of    
 Helmholtz problem \eqref{Helm_PDE}.   
\end{proposition}

%

\subsubsection*{Regularised Lagrangians for Helmholtz}

The Lagrangians discussed  above provide a variational characterisation of the Helmholtz equation originating from first principles. However, they are indefinite and lack coercivity. To address this, we construct regularised functionals that penalise the Helmholtz residual in a least-squares sense. This allows us to gain further insight into the analytical properties of the reduced wave equation and provides a framework for designing computational  methods with potentially improved stability and robustness.

The introduction of the following real Lagrangian will be instrumental 
\begin{equation}
\begin{split}
\mathcal{F}_{\gamma} (v)
&= \int_{\Omega }  \frac12  \Bigl( |\nabla v |^2 - k^2|v|^2 \Bigr)  \d x
 -  {\Re \int_{\Omega } f  \overline{v}  \d  x} \\
&\quad + \gamma_1 \int_{\Omega} \bigl| \L v - f \bigr|^2   \d x
 +  \gamma_2 \int_{\partial \Omega} \Bigl| \partial_{\mathbf{n}} v - i k  v \Bigr|^2   \d S.
\end{split}
\end{equation}
We introduce the corresponding complex Lagrangian to incorporate the influence of the impedance boundary conditions into the principal part
\begin{equation}\label{Complex-Energy_H_gamma_wbc}
\begin{split}
\mathcal{F}_{\gamma,\,  \C \, } (v)
&= \int_{\Omega }  \frac12  \Bigl( |\nabla v |^2 - k^2|v|^2 \Bigr)  \d x
 -  {\Re \int_{\Omega } f  \overline{v}  \d  x} \\
&\quad + \gamma_1 \int_{\Omega} \bigl| \L v - f \bigr|^2   \d x
 +  \gamma_2 \int_{\partial \Omega} \Bigl| \partial_{\mathbf{n}} v - i k  v \Bigr|^2   \d S - \frac{ik}{2}\int_{\partial \Omega}|v|^2 \d x .
\end{split}
\end{equation}
To emphasise the role of $v$ and $\vbar$ as independent variables when differentiating complex Lagrangians,   and using the obvious notation, we denote 
\begin{equation}\label{Complex-Energy_v_vbar}
\begin{split}
 \mathcal{F}_{\gamma,\,  \C \, } (v,  \overline{v}) = \,  &\mathcal{F}_{\gamma,\,  \C \, } (v)\, .
\end{split}\end{equation} 
One can prove the  following key result,  stating that the homogeneous part of   $F_\gamma(u), $ $F_\gamma(u):=\mathcal{F}_\gamma(u)\big|_{f=0}$ is coercive in $\Vc.$ This requires  the following geometric assumption on $\Omega$, 
\begin{center}
	\emph{There exists a constant $L_0>0$ such that
$\vec{x}\cdot\vec{n}\ge L_0$ for all $\vec{x}\in\partial\Omega$,
i.e. $\Omega$ is strictly star-shaped with respect to the origin.}
\end{center}

The following coercivity result is important for the following.

\begin{theorem} 
\label{theorem_coerc_wbc}
Assume that $\Omega$ is strictly star-shaped with respect to the origin. 
Then, there exist positive constants
$c_0,c_1,c_2,c_3$ and thresholds $\gamma_{1,0},\gamma_{2,0}$,
all independent of $k$, such that, for
$\gamma_1\ge\gamma_{1,0}$ and $\gamma_2\ge\gamma_{2,0}$,
\begin{equation}
\label{coerc_d}
\begin{split}
F_\gamma(u) \ge
c_0\big(\|\L u\|^2+\|\nabla u\|^2+k^2\|u\|^2\big)
+ c_1\|\nabla u\|^2_{\partial\Omega}
+ c_2 k^2\|u\|^2_{\partial\Omega}
+ c_3\|\partial_{\mathbf n}u-iku\|^2_{\partial\Omega}.
\end{split}
\end{equation}
\end{theorem}

As a consequence of Theorem~\ref{theorem_coerc_wbc}, we have the following result.
\begin{theorem}\label{Theorem_SCVP_2}
Consider the Lagrangian  $\mathcal{F}_{\gamma,\,  \C \, } $ defined in \eqref{Complex-Energy_H_gamma_wbc}. 
Then, under the assumptions of Theorem \ref{theorem_coerc_wbc}, the following hold:
\begin{enumerate}
\item[(i)] {There exists $\boldsymbol{\gamma}_0=(\gamma_{1,0},\gamma_{2,0})$ such that for all $\boldsymbol{\gamma}=(\gamma_1,\gamma_2)$ with $\gamma_j\ge \gamma_{j,0}$, $j=1,2$, the  real part of $\mathcal{F}_{\gamma,\,  \C \, } \big|_{f=0}$   is strongly coercive in the norm $\|v\|_{\Vc}$.}
\item[(ii)] The stationary points of $\mathcal{F}_{\gamma,\,  \C \, } $ with respect to $\vbar, $    restricted to $\Vc, $ i.e., the elements   $u\in \Vc$ satisfying 
$$   D _ {\vbar } \, \mathcal{F}_{\gamma,\,  \C \, } ( u ) =0 \, $$
 are {exactly} the solutions of \eqref{Helm_PDE}.
\end{enumerate}
\end{theorem}

\noindent
Statement (i) is precisely Theorem \ref{theorem_coerc_wbc}.
Assuming this coercivity estimate, statement (ii) follows from the
complex version of the Lax--Milgram theorem \cite[Chapter VII]{Dautray_Lions_v2_1988}. To this end, consider the sesquilinear form $\mathcal{A}_{\gamma,\,  \C \, }$ defined by 
\begin{equation}\label{Aform_reg_wbc}
\begin{split}
\mathcal{A}_{\gamma,\,  \C \, }(u,v)
&:= \int_\Omega \nabla u\cdot \nabla \overline{v}  \d x  -  k^2 \int_\Omega u \overline{v}  \d x \\
&\quad +  2\gamma_1 \int_{\Omega} (\L u)  \overline{\L v}  \d x
 +  2\gamma_2 \int_{\partial\Omega} 
\Bigl(\partial_{\mathbf n}u - i k u \Bigr) 
\overline{\Bigl(\partial_{\mathbf n}v - i k v \Bigr)}  \d S \\
&\quad   -  {ik} \int_{\partial \Omega}u \, \vbar  \d x .
\end{split}
\end{equation}
{This form is linear in its first argument and antilinear in its second, and it arises as the first variation of $\mathcal{F}_{\gamma,\,  \C \, }$ with respect to $\vbar .$}
We conclude that
\begin{equation}\label{Energy_H_gamma_Re_2_wbc}
\begin{split}
\langle D_{\vbar\, }  \mathcal{F}_{\gamma,\,  \C \, }(u), v \rangle
=  \Bigl\{ \mathcal{A}_{\gamma,\,  \C \, }(u,v)  -  \int_{\Omega} f  \overline{\bigl(v + 2\gamma_1 \L v\bigr)}  \d x \Bigr\}.
\end{split}
\end{equation}
Thus, the stationary points of $\mathcal{F}_\gamma$ satisfy
\begin{equation}\label{Energy_H_gamma_Re_3_wbc}
\mathcal{A}_{\gamma,\,  \C \, } (u,v)  =  \int_{\Omega} f  \overline{\bigl(v + 2\gamma_1 \L v\bigr)}  \d x ,
\qquad \text{for all } v \in \Vc .
\end{equation}
Theorem \ref{theorem_coerc_wbc} implies
\begin{equation}
\Re \  \frac{1}{2} \mathcal{A}_{\gamma,\,  \C \, }(v,v)
 = \Re \  \mathcal{F}_{\gamma,\,  \C \, } (v)\Big|_{f=0}
 \ge  \tilde c  \|v\|_{\Vc}^{2},
\qquad \text{for all } v\in \Vc .
\end{equation}
{Thus $\mathcal{A}_{\gamma,\,  \C \, }$ is bounded on $\Vc \times\Vc $ and $  \mathcal{A}_{\gamma,\,  \C \, }$ is $\Vc$-elliptic, while the right-hand side $v\mapsto \int_\Omega f \overline{(v+2\gamma \L v)}$ defines a bounded antilinear functional on $\Vc$.}
Therefore, by the Lax-Milgram theorem \cite[Section~1, Chapter~VII]{Dautray_Lions_v2_1988} the variational problem 
\eqref{Energy_H_gamma_Re_3_wbc} has a unique solution. 
To complete the proof, consider $u$ being the unique weak solution of \eqref{Helm_PDE}. Then $u$ solves 
 \begin{equation*}
  \Bc(u,v) = \int_{\Omega   } f   \overline{v}   \d  x, \quad \text{for all } v\in \Hc\, . 
\end{equation*}
Furthermore, since $\Omega $ is a Lipschitz domain, $u \in \Vc, $ see \cite{MakrakisMMPP_2025variationalHelm}. In addition, 
$\L u =f$ in $L^2 (\Omega) $ and the impedance boundary conditions are satisfied in the $L^2 ( \partial \Omega) $
sense. Since $\Vc \subset \Hc $ we conclude that  $u$ satisfies \eqref{Energy_H_gamma_Re_3_wbc} and it is exactly the unique stationary point of $ \mathcal{F}_{\gamma,\,  \C \, }$ with respect to $\vbar\, .$

\subsubsection{Related minimisation problems}\label{SubS:Helm-minimisation}
Our analysis may lead to alternative approaches for approximating \eqref{Helm_PDE}. In particular, it provides a mathematical foundation for iterative methods that incorporate boundary conditions either as Lagrange multipliers or through fixed-point iterations. To fix ideas, consider the following problem: assume that an approximation of $u, $ denoted by $\tilde w$, is given, and we seek to update it by solving the following problem:
  Find $w \in \Vc$ such that 
\begin{equation}\label{fixed_p_1}
\begin{split}
  \int_\Omega \nabla w\cdot \nabla \overline{v}  \d x  &-  k^2 \int_\Omega w \overline{v}  \d x
  +
  2\gamma_1 \int_{\Omega} (\L w)  \overline{\L v}  \d x
 +  2\gamma_2 \int_{\partial\Omega} 
\Bigl(\partial_{\mathbf n}w - i k w \Bigr) 
\overline{\Bigl(\partial_{\mathbf n}v - i k v \Bigr)}  \d S \\
&\quad  =   {ik} \int_{\partial \Omega}\widetilde w \, \vbar \,   \d S  +  \int_{\Omega} f  \overline{\bigl(v + 2\gamma_1 \L v\bigr)}  \d x ,
\qquad \text{for all } v \in \Vc .
\end{split}
\end{equation}
Consider this problem with all data fixed except for  
$\widetilde w.$ We shall show that it admits a unique solution, and denote by $T$ the solution operator that maps   $\widetilde w \mapsto w :$   
$$ T : \Vc \to \Vc \, , \qquad w = T ( \widetilde w )\, .$$
We 
consider the sesquilinear form $ \mathcal{A}_{\text{wbc}} $ defined by 
\begin{equation}
\begin{split}
\mathcal{A}_{\text{wbc}}(u,v)
&:= \int_\Omega \nabla u\cdot \nabla \overline{v}  \d x  -  k^2 \int_\Omega u \overline{v}  \d x \\
&\quad +  2\gamma_1 \int_{\Omega} (\L u)  \overline{\L v}  \d x
 +  2\gamma_2 \int_{\partial\Omega} 
\Bigl(\partial_{\mathbf n}u - i k u \Bigr) 
\overline{\Bigl(\partial_{\mathbf n}v - i k v \Bigr)}  \d S .
\end{split}
\end{equation}
Then $w$ solves
\begin{equation}\label{fixed_p_2}
\mathcal{A}_{\text{wbc}} (w,v)  = {ik} \int_{\partial \Omega}\widetilde w \, \vbar \,   \d S   + \int_{\Omega} f  \overline{\bigl(v + 2\gamma_1 \L v\bigr)}  \d x    
\qquad \text{for all } v \in \Vc .
\end{equation}
Notice that  $\mathcal{A}_{\text{wbc}} $ {is Hermitian (conjugate-symmetric)}, i.e. $\mathcal{A}_{\text{wbc}} (u,v)=\overline{\mathcal{A}_{\text{wbc}} (v,u)}$.
The next theorem shows that $w$ can be obtained as a minimiser of the real Lagrangian
\begin{equation}\label{Energy_H_gamma_wbc}
\begin{split}
\widetilde {\mathcal{F}}_{\gamma} (v)
&= \int_{\Omega }  \frac12  \Bigl( |\nabla v |^2 - k^2|v|^2 \Bigr)  \d x
 -  {\Re \ {ik} \int_{\partial \Omega}\widetilde w \, \vbar \, \d S} -  {\Re \int_{\Omega } f  \overline{v}  \d  x}\\
&\quad + \gamma_1 \int_{\Omega} \bigl| \L v - f \bigr|^2   \d x
 +  \gamma_2 \int_{\partial \Omega} \Bigl| \partial_{\mathbf{n}} v - i k  v \Bigr|^2   \d S \, .
\end{split}
\end{equation}

\begin{theorem}
Consider the energy $\widetilde {\mathcal{F}}_{\gamma}$ defined in \eqref{Energy_H_gamma_wbc}, and the minimisation problem
\begin{equation}\label{min_V}
\min_{v \in \Vc} \widetilde {\mathcal{F}}_{\gamma} (v).
\end{equation}
Then, under the assumptions of Theorem \ref{theorem_coerc_wbc}, the following hold:
\begin{enumerate}
\item[(i)] {There exists $\boldsymbol{\gamma}_0=(\gamma_{1,0},\gamma_{2,0})$ such that for all $\boldsymbol{\gamma}=(\gamma_1,\gamma_2)$ with $\gamma_j\ge \gamma_{j,0}$, $j=1,2$, the quadratic part of $ \widetilde {\mathcal{F}}_{\gamma}$ (i.e. with $f=0$ and $\widetilde w=0$) is strongly coercive in the norm $\|v\|_{\Vc}$.}
\item[(ii)] The problem \eqref{min_V} has a unique solution {that coincides with the   solution of} \eqref{fixed_p_2}.
\item[(iii)] Both \eqref{min_V} and \eqref{fixed_p_2} admit a unique solution $w$  that satisfies: 
$$ \| w\| _{\Vc} \leq C _{T} \big ( \| \widetilde w \|  _{\Vc} + \| f\|_{L^2 (  \Omega)} \big )  \, .$$
Furthermore, if in addition $\gamma_1\ge \tilde \gamma_{1,0}\geq \gamma_{1,0}$ for an appropriately chosen $\tilde \gamma_{1,0},$ $T$ is a contraction in $ L^2 (\partial \Omega):$ 
$$  \|T( \widetilde w _ 1) - T(\widetilde w _2) \| _{L^2 (\partial \Omega)} 
\leq \, \eta \,  \| \widetilde w _ 1 - \widetilde w _2\|  _{L^2 (\partial \Omega)}  \, , \quad \text{for a constant } \ 0< \eta < 1\, .   $$
\end{enumerate}
\end{theorem}

As before (i) follows by Theorem \ref{theorem_coerc_wbc}.    Statements (ii) and (iii) then follow from the {complex} Lax-Milgram theorem. In fact,   the sesquilinear form $ \mathcal{A}_{\text{wbc}} $  arises as the first variation of the \emph{real} Lagrangian $\widetilde {\mathcal{F}}_{\gamma} $.
We conclude that
\begin{equation}\label{Energy_H_gamma_Re_2_wbc-d2}
\begin{split}
\langle D\widetilde {\mathcal{F}}_{\gamma}(u), v \rangle
= \Re \Bigl\{ \mathcal{A}_{\text{wbc}}(u,v)  -   \ {ik} \int_{\partial \Omega}\widetilde w \, \vbar \,  \d x - \int_{\Omega} f  \overline{\bigl(v + 2\gamma_1 \L v\bigr)}  \d x \Bigr\}.
\end{split}
\end{equation}
Thus, the stationary points of $\widetilde {\mathcal{F}}_{\gamma} $ are solutions of \eqref{fixed_p_2}. 
Theorem \ref{theorem_coerc_wbc} implies
\begin{equation}\label{Energy_H_gamma_Re_4_wbc}
\Re  \frac{1}{2} \mathcal{A}_{\text{wbc}}(v,v) 
 \ge  \tilde c  \|v\|_{\Vc}^{2},
\qquad \text{for all } v\in \Vc .
\end{equation}
In addition the right-hand side of  \eqref{fixed_p_2} defines a bounded antilinear functional on $\Vc$. {Therefore, by \cite[Section~1, Chapter~VII]{Dautray_Lions_v2_1988}, the variational problem \eqref{fixed_p_2} admits a unique solution, and statement (ii) follows.} 
The fact that $T$ is a contraction in $L^2(\partial \Omega)$ follows from the observation that the boundary term involving $\widetilde w$ in \eqref{fixed_p_2} is independent of $\gamma_1$ and $\gamma_2$, together with the coercivity bound established in \cite{MakrakisMMPP_2025variationalHelm}, provided that $\beta$ is chosen such that
in $\Big(2\beta - \alpha L - \frac{\alpha L^2}{\varepsilon_2}\Big)$
is sufficiently large, which in turn implies that $\gamma_1$ must be taken sufficiently large. The remainder of statement (iii) is then immediate.

\subsubsection*{A convergent fixed point iteration} Let $\rho >0$ and consider  the following iterative scheme:
  \begin{equation}\label{fixed_point_it_tilde}
  	\widetilde w ^{n+1} = (1-\rho)\widetilde w ^{n} + \rho\,  w^{n+1} = (1-\rho)\widetilde w ^{n} + \rho\,  T (\widetilde w ^{n}) =: G (\widetilde w ^{n} )\, .
  \end{equation}  
  Utilising (iii) of the above theorem we have
  $$  \|G( \widetilde w _ 1) - G(\widetilde w _2) \| _{L^2 (\partial \Omega)} 
\leq \, \Big [(1-\rho) + \eta \rho \, \Big ]\,  \| \widetilde w _ 1 - \widetilde w _2\|  _{L^2 (\partial \Omega)}  \, .$$
Since $ 0< \eta < 1    $ it follows that $G$ is a contraction for any $\rho ,$ $0< \rho < 1\, .$
Hence 
the sequence converges to the fixed point $u$ in  $L^2 (\partial \Omega)\, .$
By the stability bound, 
$$ \| u- T(\widetilde w ^{n}) \| _{\Vc} \leq C   \, k \, \| u - \widetilde w ^{n}  \|  _{L^2 (\partial \Omega)}   \, ,$$
we conclude that $T(\widetilde w ^{n}) \to u$ in $\Vc\, .$

\subsubsection{Rellich, Morawetz identities and proof design strategy}

 The key idea of the proof is to combine in an appropriate manner the identity known as Rellich/Morawetz Identity and a lower order version of Morawetz identity needed to handle the boundary terms.

  The first identity is well known, e.g.,  \cite{CummingsFeng2006,    Moiola_Spence2014},

\begin{proposition}[Rellich/Morawetz Identity]
  \label{cor:intmorawetz}
  {  Let $u\in \Vc,$ $\L u = -\Delta u - k^2 u$   and $\Omega$ be a   domain in $\mathbb{R}^\dimx $ 
  with Lipschitz boundary $\partial \Omega$.} Then the following  identity holds:
  \begin{multline}    \label{eq:IntRellich}
    \int_\Omega 2 \Re\qp{\qp{\vec{x} \cdot \nabla u} \overline{\L  u}} \,  \d x
    +
    \int_\Omega \qp{(\dimx  -2) \norm{\nabla u}^2 - \dimx k^2 \norm{u}^2} \,  \d x \\
    +
    \int_{\partial \Omega} \qp{
      2 \Re\qp{\qp{\vec{x} \cdot \overline{\nabla u}} \qp{\nabla u \cdot \vec{n}}}
      -
      \norm{\nabla u}^2 \qp{\vec{x} \cdot \vec{n}}
      +
      k^2 \norm{u}^2 \qp{\vec{x} \cdot \vec{n}}
    } \,  \d S =0,
  \end{multline}
  where $\vec{n}$ is the outward unit normal vector on $\partial
  \Omega$, and $ \d S$ is the surface element.
\end{proposition}

The following identity will be also useful. 

\begin{proposition}[Low-order Morawetz Identity]
  \label{the:boundarymorawetz}
  Let $u\in \Vc$, $\L u = -\Delta u - k^2 u$, and let $\Omega$ be {a bounded} domain in $\mathbb{R}^\dimx$ 
  with Lipschitz boundary $\partial \Omega$. 
  Define
  \begin{equation}  \label{eq:multM0}
    M_0 u :=  - i k \beta  u,
  \end{equation}
  where $\beta \in \mathbb{R}$ is a constant. Then the following integral identity holds:
 \begin{equation} \label {M0}
    \int_\Omega  2  \Re  \big( \overline{M_0 u}  \L  u\big)    \d x
     -    2 k ^2 \beta  \int_{\partial \Omega}   |u|^2    \d S       
     +  2  \Re   \int_{\partial \Omega} i k \beta   \overline{u}   \Big ( \nabla u \cdot \vec{n}
          -ik   u \Big )    \d S  =  0.
  \end{equation}
  {When traces are not $L^2$-valued, the boundary integrals are interpreted in the \\
  $H^{-1/2}(\partial\Omega)$-$H^{1/2}(\partial\Omega)$ duality sense.}
\end{proposition}

\noindent
\emph{Design idea:\/}
Our construction of the enhanced Lagrangians  is intrinsically connected to the analysis establishing their coercivity. 
The key idea is to observe that the energy $ \mathcal{E}_P (u)$ contains the indefinite term   
\begin{equation}\label{energy_wd}
    \frac{1}{2} \int_\Omega \qp{\norm{\nabla u}^2 - k^2 \norm{u}^2} \,  \d x \, , 
\end{equation}
and at the same time a standard Rellich identity, see \eqref{eq:IntRellich}, can be written as
 \begin{equation}
      \label{eq:IntRellich_wd}
    \int_\Omega \qp{(\dimx  -2) \norm{\nabla u}^2 - \dimx k^2 \norm{u}^2} \,  \d x \\
     + \int_\Omega 2 \Re\qp{\qp{\vec{x} \cdot \nabla u} \overline{\L  u}} \,  \d x
    + B_R(u)=0,
 	 \end{equation}
where $ B_R(u)$ accounts for the boundary terms in \eqref{eq:IntRellich}. 
Then, an appropriate linear combination of the principal part of the Lagrangian for $f=0$  and {  the identity} \eqref{eq:IntRellich_wd} 
has two {  key advantages}: the energy {  changes} only by a constant factor, since the contribution of \eqref{eq:IntRellich_wd}  is zero,
and the coefficients  of $\norm{\nabla u}^2$ and $k^2 \norm{u}^2 $ {  can both be made} constants of the same sign. 
The least-squares regularisation term in the Lagrangian  functional {  makes it possible} to control the second term in \eqref{eq:IntRellich_wd}. 
{  Moreover,} the boundary term in \eqref{eq:IntRellich_wd} can also be controlled, {  particularly} through the use of {  another} 
lower-order Morawetz identity---a fact that is not immediately apparent, see \eqref{M0}. It is important to emphasise that, within this framework, the structure 
of the regularised Largangian is preserved, and coercivity is ensured for all $\gamma \geq \gamma_0$, where the constant $ \gamma_0$
can be determined explicitly.  

\subsubsection{Neural network approximation}

Our aim is consider a neural network discretisation of the augmented variational formulation. To this end we employ the fixed point iteration described in Section \ref{SubS:Helm-minimisation}.
We introduce a complex boundary variable $\lambda$ and at a fixed point, we
recover the impedance term appearing in the variational formulation.
For a fixed boundary flux $\lambda$, the primal network is obtained by
approximately minimising
\begin{equation}
\label{eq:nn_deep_uzawa_functional}
\begin{split}
\mathcal{J}_{\gamma,\lambda}(v)
&=
\frac12\int_\Omega\qp{|\nabla v|^2-k^2|v|^2}\d x
-\Re\int_\Omega f\overline v\d x
+\gamma_1\int_\Omega|\L v-f|^2\d x
\\
&\quad
+\gamma_2\int_{\partial\Omega}
|\partial_{\mathbf n}v-ikv|^2\d S
-\Re\int_{\partial\Omega}\lambda\overline v\d S .
\end{split}
\end{equation}
As discussed in Section \ref{SubS:Helm-minimisation} we aim to approximate the solution of the Helmholtz problem by a sequence of 
\emph{strongly coercive, real minimisation problems}. Each one of these problems is posed as usual in a discrete neural network space of complex functions
\[
u_\theta=u_{\theta,R}+iu_{\theta,I},
\]
where $u_{\theta,R}$ and $u_{\theta,I}$ are the two real-valued outputs of a Sinusoidal Representation Network. Derivatives in $\L u_\theta$ and on $\partial\Omega$ are computed by automatic differentiation. On a fixed boundary collocation set $\{y_m\}_{m=1}^{M_\partial}\subset\partial\Omega$, the flux is updated by the damped {fixed point} step, compare to Theorem \ref{Theorem_SCVP_2} and 
\eqref{fixed_point_it_tilde},
\begin{equation}
\label{eq:nn_uzawa_update}
\lambda^{m+1}(y)
=
(1-\rho)\lambda^m(y)
+\rho\,ik\,u_{\theta^{m+1}}(y),
\qquad
y\in\partial\Omega,
\end{equation}
with $0<\rho\le 1$. Equivalently,
\[
\lambda_R^{m+1}
=
(1-\rho)\lambda_R^m-\rho k u_I^{m+1},
\qquad
\lambda_I^{m+1}
=
(1-\rho)\lambda_I^m+\rho k u_R^{m+1}.
\]

\newpage

\endgroup

\subsection{Deep Uzawa : From saddle point to minimisation}\label{SubS:DeepUzawa}

\newcommand{\ensemble}[2]{\ensuremath{\left\{ #1:\;#2 \right\}}}
\newcommand{\qa}[1]{\ensuremath{\left\langle{#1}\right\rangle}}
\newcommand{\ltwop}[2]{\ensuremath{\qa{#1,#2}}}

This section is devoted to a general  framework for analysing constrained optimisation problems solved using neural networks. 
\begin{equation}\label{eq:intro:problem}
\min_{u \in V} J(u)
\qquad
\text{subject to }
\mathcal{C}u = g,
\end{equation}
where $J$ is a convex energy functional, $V$ is a Hilbert space of admissible functions, and the constraint $\mathcal{C}u=g$ encodes structural conditions that may include boundary conditions, divergence constraints, or inequality admissibility. Inequality constraints can be handled through an equivalent constrained formulation in which feasibility is expressed by membership in a closed constraint set in the range of $\mathcal{C}$ and subsequently enforced through the same penalised and augmented machinery.
A standard approach would be to use residual-type losses, in which the constraint is incorporated directly into the functional to be minimised. Our key idea is instead to develop variations of the classical Uzawa algorithm based on the following principles:
\begin{center}
\begin{itemize}
\begin{enumerate}
\item [(a)] the solution of \eqref{eq:intro:problem} is approximated by a sequence of minimisation problems and constraint updates;
\item [(b)] each minimisation problem is solved over discrete neural network spaces;
\item [(c)] the variational structure of each minimisation problem is strongly coercive.
\end{enumerate}
\end{itemize}
	\end{center}

We refer to the resulting algorithms as \emph{deep Uzawa} methods. At each level of approximation, the energy $J$ and the constraint $\mathcal{C}$ are replaced by discrete approximations constructed using sampling or quadrature, while the function space $V$ is replaced by a nonlinear class of neural network functions $V_{\mathcal{N}}$.

\subsubsection*{Problem setup}
\label{sec:setup}

Let $V$ and $B$ be Hilbert spaces and let $\mathcal{C}:V\to B$ be a bounded linear operator encoding a constraint. Let $J:V\to\mathbb{R}$ be Fr\'echet differentiable, coercive and strongly convex. Given $g\in B$, we seek the solution of the constrained minimisation problem
\begin{equation}\label{eq:setup:primal}
u^* = \arg\min_{v\in V} J(v)
\qquad
\text{subject to }
\mathcal{C}v = g.
\end{equation}
Assume the feasible set $\ensemble{v\in V}{\mathcal{C}v=g}$ is nonempty. Then strong convexity of $J$ implies that the constrained minimiser $u^*$ is unique.

To enforce the constraint weakly let with $B'$ denote the continuous dual of $B$, introduce a Lagrange multiplier $\lambda\in B'$ and define the Lagrangian
\begin{equation}\label{eq:setup:L}
L(u,\lambda) := J(u) + \langle \mathcal{C}u-g,\lambda\rangle_{B\times B'} ,
\end{equation}
where $\langle\cdot,\cdot\rangle_{B\times B'}$ denotes the duality pairing. We assume there exists a saddle point $(u^*,\lambda^*)\in V\times B'$ satisfying
\begin{equation}\label{eq:setup:saddle}
L(u^*,\lambda)\le L(u^*,\lambda^*) \le L(u,\lambda^*)
\qquad
\text{for all }(u,\lambda)\in V\times B'.
\end{equation}
Equivalently, $(u^*,\lambda^*)$ satisfies the KKT conditions
\begin{equation}
  \begin{split}
    J'(u^*)+\mathcal{C}^*\lambda^*&=0 \text{ in } V'
    \\
    \mathcal{C}u^*&=g \text{ in } B,
  \end{split}
\end{equation}
where $\mathcal{C}^*:B'\to V'$ is the adjoint operator.

Since we will also use quadratic augmentation, it is convenient to identify residuals in $B$ with elements of $B'$ via the Riesz map $\mathcal{R}_B:B\to B'$ defined by
\begin{equation}\label{eq:setup:riesz}
\langle \psi,\mathcal{R}_B b\rangle_{B\times B'} = \ltwop{b}{\psi}_B
\qquad
\text{for all }b,\psi\in B,
\end{equation}
where $\ltwop{\cdot}{\cdot}_B$ denotes the inner product on $B$. In particular, $\mathcal{R}_B$ is an isometric isomorphism and $\|\mathcal{R}_B b\|_{B'}=\|b\|_B$.

This is the only place where Hilbert structure on $B$ is used. If one wishes to work in a Banach constraint space, the same updates can be formulated using the duality map in place of $\mathcal{R}_B$, but we do not pursue that here.

\subsubsection*{Uzawa iteration}

To compute a saddle point, one can apply the Uzawa algorithm \cite{uzawa1958iterative}. Given $\lambda^0\in B'$ and a step size $\rho>0$, define $\{(u^k,\lambda^k)\}_{k\ge 0}$ by
\begin{equation}\label{eq:setup:uzawa}
\begin{aligned}
u^k &\in \arg\min_{u\in V} \Bigl(J(u) + \langle \mathcal{C}u-g,\lambda^k\rangle_{B\times B'}\Bigr),\\
\lambda^{k+1} &= \lambda^k + \rho \mathcal{R}_B\qp{\mathcal{C}u^k-g}.
\end{aligned}
\end{equation}
Equivalently, the multiplier update may be read weakly as: for all $\psi\in B$,
\begin{equation}\label{eq:setup:uzawa-weak}
\langle \psi,\lambda^{k+1}-\lambda^k\rangle_{B\times B'}
= \rho \ltwop{\mathcal{C}u^k-g}{\psi}_B .
\end{equation}
The scheme alternates between exact minimisation in $u$ for fixed $\lambda^k$ and an explicit ascent step in $\lambda$ driven by the constraint residual.

\subsubsection*{Augmented Uzawa iteration}

One way to enforce strong coercivity which is capable to control the Lagrange multiplier term in the functional is to    incorporate a quadratic penalty into the Lagrangian, \cite{Fortin_Glowinski:AugmentedL:1983}, giving the augmented Lagrangian, 
\begin{equation}\label{eq:setup:augLag}
L_\gamma(u,\lambda) := J(u)
+ \langle \mathcal{C}u-g,\lambda\rangle_{B\times B'}
+ \frac{\gamma}{2} \|\mathcal{C}u-g\|_B^2,
\end{equation}
where $\gamma>0$ and $\|\cdot\|_B$ is the norm induced by $\ltwop{\cdot}{\cdot}_B$. Applying Uzawa iteration to $L_\gamma$ yields $\{(u^k,\lambda^k)\}_{k\ge 0}$ defined by
\begin{equation}\label{eq:setup:auguzawa}
\begin{aligned}
u^k &\in \arg\min_{u\in V} \Bigl(J(u) + \langle \mathcal{C}u-g,\lambda^k\rangle_{B\times B'}
+ \frac{\gamma}{2} \|\mathcal{C}u-g\|_B^2\Bigr),\\
\lambda^{k+1} &= \lambda^k + \rho \mathcal{R}_B\qp{\mathcal{C}u^k-g}.
\end{aligned}
\end{equation}
The augmented term improves conditioning of the primal subproblem and stabilises the iteration in stiff regimes. It is also compatible with sampling-based surrogates used in neural discretisations, where $\mathcal{C}$ and $J$ are replaced by discrete approximations and the multiplier update is performed in a finite-dimensional surrogate dual space. This is the setting analysed in the sequel. One can obtain the convergence  at the continuous level by adopting classical arguments, \cite {Glow_Lions_Trem:1981}, \cite{Fortin_Glowinski:AugmentedL:1983}:

\begin{theorem}[Augmented Uzawa convergence]\label{thm:auguzawa}
Let $V$ and $B$ be Hilbert spaces and let $\mathcal{C}:V\to B$ be a bounded linear operator. Let $J:V\to\mathbb{R}$ be convex, coercive and G\^ateaux differentiable. Assume the constrained problem
\begin{equation}\label{eq:thm:primal}
\min\{J(v):v\in V,\ \mathcal{C}v=g\}
\end{equation}
admits a unique minimiser $u^*\in V$ and assume there exists at least one $\lambda^*\in B'$ such that $(u^*,\lambda^*)$ is a saddle point of the augmented Lagrangian
\begin{equation}\label{eq:thm:augLag}
L_\gamma(u,\lambda)\ :=\ J(u)+\ltwop{\mathcal{C}u-g}{\lambda}_{B\times B'}
+\frac{\gamma}{2}\|\mathcal{C}u-g\|_B^2,
\end{equation}
with some $\gamma>0$. Let $\mathcal{R}_B:B\to B'$ denote the Riesz map and define $(u^k,\lambda^k)$ by the augmented Uzawa iteration
\begin{equation}\label{eq:aug-uzawa-iter}
\begin{aligned}
  u^k  &\in  \arg\min_{u\in V}\ \Bigl[J(u)+\ltwop{\mathcal{C}u-g}{\lambda^k}_{B\times B'}
  +\frac{\gamma}{2}\|\mathcal{C}u-g\|_B^2\Bigr],\\
  \lambda^{k+1}  &=  \lambda^k + \rho \mathcal{R}_B\qp{\mathcal{C}u^k-g}\quad\text{in }B',
\end{aligned}
\end{equation}
with $\rho\in(0,2\gamma)$. Then the constraint residuals $r_k:=\mathcal{C}u^k-g$ satisfy $\sum_{k\ge 0}\|r_k\|_B^2<\infty$ and hence $r_k\to 0$ strongly in $B$. Moreover, $u^k\rightharpoonup u^*$ weakly in $V$. Finally, the multiplier sequence $\{\lambda^k\}$ is bounded in $B'$ and converges weakly in $B'$ to some $\bar\lambda\in\Lambda^*$, where
\begin{equation}\label{eq:thm:Lambda}
\Lambda^*:=\ensemble{\lambda\in B'}{u^*\in \arg\min_{u\in V} L_\gamma(u,\lambda)}.
\end{equation}
If $\Lambda^*$ is a singleton then $\lambda^k\rightharpoonup \lambda^*$ weakly in $B'$.
\end{theorem}

The arguments are based heavily on the comparison of the optimality conditions and of the Euler-Lagrange equations for the minimisation step. Preliminary results suggest that 
alternative, neural network compatible  proofs
utilising  the framework of Section~\ref{S:Analysis} can be derived. 

 \subsection{Bibliography and comments}\label{SubS:an-to-ad:bib}

In Section \ref{S:Analysis-abstract}  we established that the training process---that is, the evaluation of the loss functional through quadrature or other discrete approximations---can play a decisive role in determining the stability, and in particular the coercivity, of the resulting method. 
It is natural therefore to seek alternative training methods with emphasis on qualitative behaviour  and improved performance.   In Section~\ref{SubS:RK-PINN} based on \cite{akrivis2025rungekuttaphysicsinformedneural}, we establish a systematic connection between Runge--Kutta (RK) time-discrete schemes and neural solvers through the training mechanism. The resulting methods inherit key stability properties and computational features of the underlying RK schemes. This framework leads, in particular, to energy-conservative neural network solvers for the wave equation. From a theoretical perspective, stability in the parabolic case is established using new maximal-regularity estimates for classical RK schemes, also developed in \cite{akrivis2025rungekuttaphysicsinformedneural}.
In Section~\ref{SubS:FE-training}, based on \cite{grekas2024deepritzfinite}, we propose incorporating finite element techniques into the training of neural solvers for elliptic problems in settings where meshes can be constructed efficiently. In such cases, the resulting methods exhibit strong stability properties and can offer improved memory and computational efficiency. In \cite{grekas2025pinndgresidualneuralnetwork}, this approach is extended to  residual methods 
by merging Discontinuous Galerkin finite elements and PINNs. The need of specially designed loss functionals satisfying the Discrete Variational Stability is linked to well known coercivity of DG formulations.

Although the minimisation problem over neural network functions is inherently non-convex, our work suggests that algorithms based on strongly coercive functionals, formulated at the continuous infinite-dimensional level, can inherit favourable stability and performance properties at the discrete level. This observation motivates the development, already at the PDE level, of alternative variational principles with enhanced analytical and computational properties. 
In Section~\ref{SubS:Helmholtz} based on  \cite{MakrakisMMPP_2025variationalHelm},   Morawetz identities are used to derive novel strongly coercive variational principles for the Helmholtz equation, leading through fixed point iterations to efficient neural network solvers.

Deep Uzawa methods discussed in Section~\ref{SubS:DeepUzawa} are based on a related principle: saddle-point problems are reformulated and approximated through a sequence of strongly coercive minimisation problems combined with constraint updates. This approach has been applied to PDE-constrained optimisation in \cite{makr_pryer2024deepUz}, to the weak enforcement of boundary conditions through Lagrange multipliers in \cite{makr_pryer2024uzLagr}, and to stationary kinetic equations in \cite{makridakis2025deepuzawakinetictransport}.

\subsubsection*{Acknowledgements}
\addcontentsline{toc}{section}{Acknowledgements}
These notes are based on the lectures I gave at the mathematics summer
school and conference held within the 2025 edition of the festival Festum-$\pi$, led by
D.~Dervis-Bournias and C.~Villani.
  I am grateful to Professor C\'edric Villani for his invitation
and his encouragement, and I would also like to warmly thank Professors
N.~Alikakos, P.~Giannotis and G.~Fournodavlos. The material presented here
is based on joint work and discussions with several collaborators and
friends, among them G.~Akrivis, G.~Arampatzis, L.~Banjai,  M.~Dashti, D.~Gazoulis, E.~Georgoulis,
I.~Gkannis, G.~Grekas, P.~Herbert, M.~Jensen, T.~Katsaounis, S.~Katsarakis, M.~Katsoulakis, K.~Koumatos, O.~Lakkis, B.~Li, G.~Makrakis, D.~Mitsoudis, I.~Pantazis, M.~Plexousakis,
N.~Rekatsinas, P.~Rosakis, G.~Savar\'e, K.~Smaragdakis, D.~Tsagkarogiannis, A.~Tzavaras, A.~Veeser, A.~Vikelis, J.~Xu, and K.~Zygalakis to whom I owe much of my understanding of the subject; most of
the ideas were developed during my long-standing collaboration with
Michalis Loulakis and Tristan Pryer, to whom I am especially grateful.

\bibliographystyle{plain}

\bibliography{proposal_bib_2025}
\end{document}